\documentclass[12pt]{article}

\usepackage{epic}
\usepackage{eepic}
\usepackage{epsf} 
\usepackage{ amssymb, amscd}
\usepackage{amsmath}
\usepackage{color}
\usepackage{epsfig}
\usepackage{graphicx}

\numberwithin{equation}{section}

\def\cA            {{\mathcal{A}}}

\def\cD            {{\mathcal{D}}}
\def\cE            {{\mathcal{E}}}

\def\cM            {{\mathcal{M}}}

\def\cO            {{\mathcal{O}}}
\def\cP            {{\mathcal{P}}}

\def\cV            {{\mathcal{V}}}

\def\cX            {{\mathcal{X}}}

\def\cZ            {{\mathcal{Z}}}

\def\bbC           {\mathbb{C}}

\def\bbQ           {\mathbb{Q}}
\def\bbR           {\mathbb{R}}
\def\bbT           {\mathbb{T}}
\def\bbZ           {\mathbb{Z}}

\def\MXN           {{}_M {\mathcal{X}}_N}
\def\MXM           {{}_M {\mathcal{X}}_M}

\def\NXN           {{}_N {\mathcal{X}}_N}

\def\NXM           {{}_N {\mathcal{X}}_M}
\def\Tr            {\hbox{Tr}}

\def\a             {\alpha}
\def\la            {\lambda}

\def\Tr            {{\mathrm{Tr}}}

\def\i{{\rm i}}
\def\sdprod{{\times\!\vrule height5pt depth0pt width0.4pt\,}}

\def\AXB           {{}_A {\mathcal{X}}_B}
\def\BXB           {{}_B {\mathcal{X}}_B}
\def\AXA           {{}_A {\mathcal{X}}_A}

\def\id            {\mathrm{id}}

\def\Hg          {\mathrm{Hg}}

\title{The exoticness and realisability of twisted Haagerup--Izumi modular data}
\author{
{\sc David E.\ Evans}\\
 {\footnotesize School of Mathematics, Cardiff University,}\\
 {\footnotesize Senghennydd Road, Cardiff CF24 4AG, Wales, U.K.}\\
 {\footnotesize e-mail: {\tt EvansDE@cf.ac.uk}}\\ \\
 {\sc Terry  Gannon }\\
 {\footnotesize Department of Mathematics, University of Alberta,}\\
{\footnotesize Edmonton, Alberta, Canada T6G 2G1}\\
{\footnotesize e-mail: {\tt tgannon@math.ualberta.ca}} }

\begin{document}
\maketitle

\begin{abstract} The quantum double of the Haagerup subfactor, the first irreducible finite depth subfactor with index
above 4, is the most obvious candidate for exotic modular data. We show that its modular data $\cD \Hg$
fits into a family
$\cD^\omega \Hg_{2n+1}$, where  $n\ge 0$ and  $\omega\in \bbZ_{2n+1}$.  
We show $\cD^0 \Hg_{2n+1}$ is related to the subfactors Izumi hypothetically associates to
the cyclic groups $\bbZ_{2n+1}$. Their modular data comes equipped with 
canonical and dual canonical modular invariants; we compute the corresponding
alpha-inductions etc.
In addition, we show there are (respectively) 1, 2, 0  subfactors  of
Izumi type $\bbZ_7$, $\bbZ_9$ and $\bbZ_3^2$, and find numerical evidence for 2, 1, 1, 1, 2 
subfactors of Izumi type $\bbZ_{11},\bbZ_{13},\bbZ_{15},\bbZ_{17},\bbZ_{19}$ 
(previously, Izumi had shown uniqueness for $\bbZ_3$ and
$\bbZ_5$), and we identify their modular data.
 We explain how  $\cD \Hg$ (more generally
$\cD^\omega \Hg_{2n+1}$) is a \textit{graft} of the quantum double $\cD Sym(3)$ (resp. the twisted double
$\cD^\omega D_{2n+1}$) by affine so(13) (resp. so$(4n^2+4n+5)$) at level 2.
We discuss the vertex operator algebra (or conformal field theory)
realisation of  the modular data 
$\cD^\omega \Hg_{2n+1}$. For example we show there are exactly 2 possible character
vectors (giving graded dimensions of all modules) for the Haagerup VOA at
central charge $c=8$. It seems unlikely that any of this twisted Haagerup-Izumi  modular data
{can be regarded as} exotic, in any reasonable sense.

\end{abstract}

{\footnotesize
\tableofcontents
}

\section{Introduction}

From the early years of conformal field theory (CFT) (see e.g. \cite{MoSe,Wi}), we find:

\medskip\noindent{\textbf{Speculation.}}  \textit {The standard
constructions (orbifolds, cosets, simple-current extensions, ...) applied to the basic theories (lattice compactifications, 
affine Kac-Moody algebras, ...) exhaust all rational theories.}\medskip

Perhaps this should be applied to the modular tensor category rather than
the full CFT, in which case this would constitute a sort of generalised Tannaka-Krein duality 
holding for modular tensor categories, {where the dual to the category (i.e. the
analogue of  the 
compact group) is (say) a vertex operator algebra constructed
in a standard way from the basic examples.}

Possible counterexamples have been around just as long. For instance 
Walton \cite{Wal} explained that rank-level duality applied to conformal embeddings of affine 
algebra CFTs  (the simplest being $A_{1,10}\subset C_{2,1}$, yielding an extension of 
$A_{9,2}$) yield what seem to be new CFTs. The possibility  that rank-level
duality is somehow inherently sick was eliminated by Xu \cite{Xu}, who realised many
of these examples with completely rational nets of subfactors (these should correspond to
rational CFT). Perhaps all that this accomplishes is to insist
that  rank-level duality (or in subfactor language the \textit{mirror}) should be included as 
one of  those standard constructions. More recently,
\cite{DoGe} proposed four modular data (each with 6 or 7 primaries) as possible 
counterexamples to this speculation, but didn't show they could be realised by a rational CFT
(or vertex operator algebra or net of subfactors).

Perhaps the most obvious place to look for a truly  exotic example is the Haagerup subfactor 
\cite{ahaag}. It is the first irreducible finite depth subfactor with
index greater than 4 --- it has index $(5+\sqrt{13})/2\approx  4.30278$. It is generally regarded as exotic, since
(so far) it can only be constructed by hand without any natural algebraic symmetries.
However it is nonbraided and so to get a braided system (which has a chance
to correspond to the fusion ring and modular data of a rational CFT) we should take its quantum 
double, or equivalently asymptotic inclusion (see
Ocneanu as in \cite[Chapter 12]{EK}) or Longo-Rehren inclusion  \cite{LR,masuda}. Its modular 
data $\cD \Hg$ was computed in \cite{iz3},
and subfactor realisations of its modular invariants  were studied in \cite{EP}.

\medskip\noindent\textbf{Question 1.}  \textit{Is this Haagerup quantum double realized by a rational 
CFT (or  rational vertex operator algebra --- see Definition 3 below ---  or completely rational net of subfactors \cite{KLM}), in the sense that they 
share the same modular data $\cD\Hg$?}\medskip

 We would 
want a deeper relation between the Haagerup double and the corresponding vertex operator 
algebra (VOA) than merely that their modular data coincide, but for now that can suffice.

A more direct realisation of the Haagerup subfactor $N\subset M$ in a rational CFT could be that its $N$-$N$ 
sectors (say) form the algebra of defect lines (or full system) of a rational CFT. After all this will
in general fail to be commutative. However, even this cannot happen: the modular data of
the quantum double
of a full system must be in factorised form ${\mathcal{X}} \otimes {\mathcal{X}}^{\mathrm{opp}}$.
{This is because the quantum double of the full system is equivalent to the double 
of the fusion algebra ${\mathcal{X}}$ by \cite{m1}, \cite[Cor.\ 2.2]{ost}, and
this double factors as ${\mathcal{X}} \otimes {\mathcal{X}}^{\mathrm{opp}}$
by \cite[Prop.\ 2.2]{EK1} as long as the braiding on ${\mathcal{X}}$ is non-degenerate. }
However,  the modular data of the double of the Haagerup is easily seen
to not be in factorised form \cite{HRW}. {There remains the possibility, which
we won't explore in this paper, that the 
Haagerup systems are the full system of a \textit{degenerately} braided system. Note that much of the 
theory of modular invariants, including alpha-induction and the  full system, holds for degenerately braided systems \cite{BEK1, BEK2, BE4, BE5, BE6}.}

\medskip\noindent\textbf{Question 2.} \textit{If $\cD\Hg$ is realised by a
rational VOA, say, then is that VOA
exotic, in the sense that it cannot be constructed from standard methods and
examples?}\medskip

{The three known constructions of the Haagerup subfactor (namely Haagerup's
connection computation \cite{Haag,ahaag}, Izumi's Cuntz algebra construction \cite{iz3} and
Peters' planar algebra construction \cite{pet}) are largely combinatorial tours de force (although
Izumi's suggests an underlying cyclic group), and for this reason
the Haagerup subfactor is generally regarded as exotic. Even given this,} it is conceivable  that 
its quantum double $\cD \Hg$ could be constructed directly in more standard ways. \cite{HRW} 
confirm that $\cD \Hg$ does not fall into the simplest possibilities, namely the modular data 
of an affine algebra, lattice, or  finite group.

We argue that the answer to Question 1 is yes, and that for Question 2 is no.

In particular, we explain why such a VOA should be a conformal subalgebra of the central charge $c=8$ 
VOA $\cV(E_6A_2)$ corresponding to (root) lattice $E_{6}\oplus A_{2}$
(by \textit{conformal subalgebra} we mean a subVOA of identical central charge and identical
conformal vector). A 
more familiar conformal subVOA of
$\cV(E_{6}A_{2})$ is an  order-2 orbifold,  realising the modular data $\cD S_3$ of the quantum double of 
the symmetric group $Sym(3)=S_3$.  An early indication for us that the  Haagerup is 
related to $S_3$ was the similarity of their tube algebras. Tube algebras were 
introduced by Ocneanu {(see e.g. \cite{EK})} as a way of computing the {irreducible objects and 
fusion rules of the quantum doubles of
subfactors, and developed further by Izumi \cite{iz1, iz3} in particular in analysing and determining 
their modular data}. The list of modular invariants and nimreps for the doubles of the Haagerup 
and $S_3$ \cite{CGR,EP,EP3}
are also strikingly similar, and in fact the modular invariants lead us to $E_{6}\oplus A_{2}$.
The connection between modular data and VOAs is made much more explicit using character
vectors. A second, equally promising possibility for a VOA realisation of
$\cD\Hg$, by GKO cosets using the affine algebra VOA $\cV(B_{6,2})$, is discussed
but not analysed in detail.

We generalize the modular data $\cD \Hg$ in two directions, by fiting it into
an infinite sequence, and showing the $n$th term in this sequence can be twisted by
$\bbZ_{2n+1}$. The role of $S_3$ is now played by the twisted quantum double of $D_{2n+1}$, and 
$\cV(E_{6} A_{2})$ (which realises the $\cD \bbZ_3$ modular data) by  a holomorphic orbifold
by $\bbZ_{2n+1}$.  Subfactor realizations of $\cD^0\Hg_{2n+1}$ should be provided by
Izumi's hypothetical generalisation of the Haagerup: we construct a unique subfactor of
Izumi type which realises $\cD^0\Hg_\nu$ for all $1\le \nu\le 19$ and conjecture that indeed this continues
for every $n$. For this reason we suggest calling $\cD^0\Hg_{2n+1}$ \textit{Haagerup-Izumi} 
modular data. von Neumann algebra realisations of \textit{twisted} Haagerup-Izumi modular data $\cD^\omega
\Hg_{2n+1}$ are still unclear.

\medskip\noindent\textbf{Question 3.} \textit{Can we orbifold a VOA by something more general
than a group?}\medskip

If so, this should provide a simple construction of VOAs realising $\cD^\omega \Hg_{2n+1}$,
starting from a holomorphic orbifold by $\bbZ_{2n+1}$. Perhaps this is not that unreasonable
--- after all, a subfactor itself is really a generalised orbifold, in the sense that the quantum double 
of a subfactor $M^G\subset M$ of fixed points recovers $\cD G$. {Indeed the whole motivation of subfactor theory (see e.g. \cite{EK}) is to understand non-grouplike  quantum symmetries 
through subfactors or inclusions of factors $N \subset M$, and associated notions of $N$-$M$,
$N$-$M$, $M$-$M$ irreducible bimodules or sectors, their fusion rules, paragroups, $\lambda$-lattices
or planar algebras. This notion of generalised orbifold construction also features in the framework of \cite{FFRS}.}

Our work suggests a further generalisation. Consider a quadruple 
$(K,H,\alpha,\omega)$ where $K$ and $H$ are finite groups, and
$H$ acts freely on $K$, i.e. the group homomorphism $\alpha:H\rightarrow
\mathrm{Aut}(K)$ obeys $\alpha(h)k=k$ iff $h=e$ or $k=e$. {Freeness  
implies the projection $h\mapsto\mathrm{Out}(K)$ is injective;} such a 
semi-direct product $D=K\sdprod_\alpha H$ is called a \textit{Frobenius group}.
The twist $\omega\in H^3(B_K;\bbT)$ should be compatible with $\alpha$ in the sense that
it lies in the image of the natural map $H^3(B_D;\bbT)\rightarrow H^3(B_K;\bbT)$.
The Haagerup subfactor is associated to 
$(\bbZ_3,\bbZ_2,\alpha(1)=-,[0])$ and $D=S_3$.

\medskip\noindent\textbf{Wildly Optimistic Guess.} \textit{Let $(K,H,\alpha,\omega)$ be any
quadruple defined above, corresponding to Frobenius group $D=K\sdprod_\alpha H$.
 There is an irreducible finite-depth subfactor associated to
the pair $(K,H)$; the triple $(K,H,\alpha)$ is realised by a $Q$-system  of
endomorphisms. To any such quadruple there is a rational VOA $\cV(K,H,\alpha,\omega)$, and
a completely rational net of subfactors realising the modular data of this twisted quantum double.
This VOA is a generalised orbifold (controlled in some sense by $H$) of a holomorphic
orbifold $\cV_K$ by $K$. This VOA $\cV_K$ also contains a 
holomorphic orbifold $\cV_D$; both $\cV(K,H,\alpha,\omega)$ and
$\cV_D$ contain a common rational VOA.
For $K$ an odd cyclic group and $H=\bbZ_2$ acting
by $n\mapsto-n$, this recovers the modular data given in Section 3.2 and also a subfactor
of Izumi type $K$, i.e. a solution to equations (7.1)-(7.5) of \cite{iz3}; the
$K=\bbZ_3$ special case
recovers the Haagerup subfactor.}\medskip

This guess suggests that the Haagerup
subfactor itself cannot be regarded as exotic. But this guess is in some control only for the odd cyclic case discussed
above, where it seems quite plausible.
The reason for requiring freeness is that this condition plays an
important role in the derivation of $\cD^\omega\Hg_\nu$ in Section 3.2, as well as in
\cite{iz3}; probably it can be dropped but at a cost of complexity.   

The possibility
that the Haagerup subfactor is related in some way to $S_3$ appears to have been first made 
in \cite{E2}. The present paper clarifies, deepens, and generalises that relation.
The remainder of this first section recalls modular data and modular invariants
(Section 1.1), reviews the basic theory of subfactors and the tube algebra (Section 1.2), and
discusses the modular data and character vectors of VOAs (Section 1.4). {We also explain in
Section 1.3 how to recover the original subfactor from its quantum double --- unlike the
rest of Section 1, this subsection contains original material.} Section 2 develops 
the relation between the Haagerup subfactor and $S_3$. Section 3 puts the Haagerup modular
data into a sequence and twists it. We also introduce the notion of
\textit{grafting} and discuss the role of affine so$(\nu^2+4)$ at level 2.
Section 4  relates our untwisted sequence to
Izumi's hypothetical family of 3-star, 5-star, 7-star, ... subfactors. Section 5
discusses VOA interpretations for our twisted sequence, and  explains the
Haagerup-dihedral diamond which generalises Section 2.

The Haagerup subfactor arose in Haagerup's classification {\cite{Haag}} of irreducible finite depth 
subfactors of index between 4 and 
$3+\sqrt{3}\approx  4.73205$. Two other subfactors {$N\subset M$} (both also 
regarded as exotic) appear there:
the \textit{Asaeda-Haagerup} subfactor with Jones index $(5+\sqrt{17})/2\approx  4.56155$ 
\cite{ahaag}, and the 
\textit{extended Haagerup} subfactor \cite{BMPS} with index  $\approx 4.37720$. Both seem unrelated to our 
family, and neither lies in a known sequence. {As with the Haagerup subfactor,
neither their $N$-$N$ nor $M$-$M$ systems are  braided}. Of course 
Questions 1 and 2 should also be asked of them, but no serious work can begin until
the modular data of their doubles has been computed.

\subsection{Modular data and modular invariants}

Modular data arises naturally in several contexts (for instance CFT, subfactors, and VOAs) ---
see \cite{Gan2} for a review. In CFT, modular invariants correspond
to the 1-loop closed string partition function.

\medskip\noindent\textbf{Definition 1.} Modular data $(\Phi,0,S,T)$ \textit{consists of a finite set $\Phi$,
an element $0\in\Phi$, and matrices $S=(S_{i,j})_{i,j\in\Phi}$ and 
$T=(T_{i,j})_{i,j\in\Phi}$, such that:} 

\noindent a) \textit{$S,T$ are unitary, $S$ is symmetric, $T$ is diagonal;}

\noindent b) \textit{$S^2=(ST)^3$, $S^4=I$, the identity matrix;}

\noindent c) \textit{$C:=S^2$ is a permutation matrix, and row 0 of $S$ consists of nonzero real numbers;}

\noindent d) \textit{for all $i,j,k\in\Phi$, the following quantities are
nonnegative integers:} 
\begin{equation}\label{verl}N_{i,j}^k:=\sum_{l\in\Phi}\frac{S_{i,l}S_{j,l}\overline{S_{k,l}}}{S_{0,l}}\,.\end{equation}

 We call modular data $(\Phi,0,S,T)$ and $(\Phi',0',S',T')$ 
\textit{equivalent} if there is a bijection $\phi:\Phi\rightarrow\Phi'$ such that $\phi(0)=0'$,
$S'_{\phi(i),\phi(j)}=S_{i,j}$ and $T'_{\phi(i),\phi(i)}=T_{i,i}$ for all $i,j\in\Phi$.
The  $i\in\Phi$ are called \textit{primaries} and $0\in\Phi$ is called the \textit{vacuum}.
The order 1 or 2  permutation $C$ is called \textit{charge-conjugation}. $T_{0,0}=e^{-\pi\i c/12}$ for a
real number $c$ (only defined mod 24 by the modular data) called the \textit{central charge}.
 Equation (\ref{verl}) is called the Verlinde formula; the 
quantities $N_{ij}^k$ are called \textit{fusion coefficients} and comprise the structure constants
of a commutative associative algebra called the \textit{fusion ring}. 
  The coefficients $N_{i,j,k}:=N_{i,j}^{Ci}$ are invariant under all 6 reorderings of indices.
  Also, $N_{0,i,j}=\delta_{i,Cj}$.

Condition b) says that modular data generates a representation $\rho$ of the modular group
SL$_2(\bbZ)$ through the assignment
\begin{equation}\label{modrep}
S=\rho\left(\begin{matrix} 0&-1\cr 1&0\end{matrix}\right)\ ,\
T=\rho\left(\begin{matrix} 1&1\cr 0&1\end{matrix}\right)\,.\end{equation}
Most examples considered in this paper have the additional property
that $C=I$, i.e. $S$ is real and the representation is factors to PSL$_2(\bbZ)$. 
We assume throughout that the first row of $S$ is strictly positive --- such
 modular data is called \textit{unitary}.

From Definition 1 we obtain
\begin{equation}\label{SNT} S_{i,j}=\overline{T}_{i,i}\overline{T}_{j,j}T_{0,0}\sum_{k\in\Phi}
T_{k,k}S_{k,0}N_{i,j}^k\,.\end{equation}
Thus if $T$ is known exactly but $S$ only approximately, \eqref{SNT} together with the
integrality of the fusion coefficients \eqref{verl} can be used to determine $S$ exactly (the
\textit{quantum-dimension} $S_{i,0}/S_{0,0}$ is the Perron-Frobenius  eigenvalue of the matrix $N_i=(N_{i,j}^k)$, and
$S_{0,0}>0$ is determined by the square 
$\sum_i(S_{i,0}/S_{0,0})^2$ of the \textit{global dimension} $1/S_{0,0}$). Equivalently, if $S,T$ 
and $S',T$ are both modular data sharing the same $T$, and $S'$ is sufficiently close to
$S$, then $S=S'$. This is how we'll
identify the modular data of the double of subfactors in Section 4.1, given numerical estimates
of $S$. 

One consequence of this definition is that the entries of $S$ and $T$  lie
in a cyclotomic extension $\bbQ[\xi_L]$ of the rationals (throughout this paper we
write $\xi_k$ for $\exp[2\pi \i/k]$). For any Galois automorphism
$\sigma\in\mathrm{Gal}(\bbQ[\xi_L]/\bbQ)\cong\bbZ^\times_L$, $(\Phi,0,\sigma S,\sigma T)$
will be (generally inequivalent) modular data with identical fusion coefficients, where $\sigma$
acts entry-wise on $S$ and $T$. It can be shown \cite{gancoste} that for any 
$\sigma\in \mathrm{Gal}(\bbQ[\xi_L]/\bbQ)$, there is a permutation 
$j\mapsto \sigma j$ of $\Phi$ and choice $s_\sigma(j)$ of signs such that 
\begin{equation}\label{galois} \sigma(S_{i,j})=s_\sigma(i)\,S_{\sigma i,j}
=s_\sigma(j)\,S_{i,\sigma j}\,.\end{equation}

An important class of examples of modular data comes from finite groups $G$.
Write $e$ for the identity and $g^h$ for $h^{-1}gh$.
Select a cocycle $\omega\in Z^3(B_G;\mathbb{T})$, where we write $\bbT$ for the unit circle in $\bbC$; 
the corresponding modular data depends
up to equivalence only on its class in $H^3(B_G;\mathbb{T})$. Define
\begin{equation}\label{theta}\theta_g(h,l)=\omega(h,g^h,l)\,\overline{\omega}(g,h,l)\,\overline{\omega}(h,l,g^{hl})\,.\end{equation}
Then the restriction of each $\theta_a$ to the centraliser $\mathcal{Z}_a:=\{h\in G:ha=ah\}$ in $G$ 
is a normalised 2-cocycle.
 In this paper we are only interested in the \textit{cohomologically trivial} case where all $\theta_a$
 are {coboundaries} (this will happen for instance whenever the Schur multiplier $H^2(B_{\mathcal{Z}_a};\bbT)$ 
 of each centraliser $\mathcal{Z}_a$ is trivial).
This means that there exist  1-cochains
$\epsilon_a : \mathcal{Z}_a \rightarrow \mathbb{T}$ for which $\epsilon_a(e)=1$ and both
\begin{align}
&&\theta_a(h,g) = (\delta\epsilon_a)(h,g) =
\epsilon_a(h)\,\epsilon_a(g)\,\overline{\epsilon}_a(hg)\,,\label{epsilon1}\\
&&\epsilon_{x^{-1}ax}(x^{-1}hx)=\theta_a(x,x^{-1}hx)\,\overline{\theta_a}(h,x)\,
\epsilon_a(h)\,,\label{epsilon2}
\end{align}
for all $g,h\in\mathcal{Z}_a$ and $x\in G$.  The primaries are all pairs $(a,\pi)$ where $a$
runs through representatives of conjugacy classes in $G$ and 
$\pi\in\widehat{\cZ_a}=\mathrm{Irr}(\cZ_a) $ are irreducible
representations (\textit{irreps}). The
vacuum is $(e,1)$. The modular matrices are \cite{DV3}
\begin{subequations}\label{moddataG}\begin{align}
&&S^{\omega}_{(a,\chi),(a',\chi')}={1\over |\mathcal{Z}_a| |\mathcal{Z}_{a'}|}
\sum_{g\in G(a,a')}\overline{\chi}(ga'g^{-1})\,\overline{\chi}'(a^g)\,\overline{\epsilon}_a 
(ga'g^{-1})\,\overline{\epsilon}_{ga'g^{-1}}(a)\,,
\\ &&T^{\omega}_{(a,\chi),(a',\chi')}
= \delta_{a,a'}\delta_{\chi,\chi'}{\chi(a) \over  \chi(e)}\,\epsilon_a(a)\,,\hskip 2truecm 
\end{align}\end{subequations}
where  $G(a,a')=\{g\in G:a^gb=ba^g\}$.
In the following we let $\cD^\omega G$ denote this modular data.
See e.g. \cite{DV3,DW,CGR} for more details; we recover this formula next subsection  using
tube algebras.

Another class of modular data is associated to even positive-definite lattices. Let $L$ be
such a lattice, and $n$ its dimension. Then $\Phi=L^*/L$ where $L^*$ is the dual of $L$. The 
vacuum 0 corresponds to the coset $[0]$. Then 
\begin{subequations}\label{moddataL}\begin{align}
&&S_{[u],[v]}=|L^*/L|^{-1/2}\exp[2\pi\i \,u\cdot v]\,,\\ &&T_{[u],[u]}=\exp[\pi\i\, u\cdot u-\pi\i\, n/12]\,.\end{align}
\end{subequations}

The final class we need is associated to affine nontwisted Lie  algebras 
$\mathfrak{g}^{(1)}$ and  positive integers $k$ (see e.g. \cite{kac}). 
Let $\mathfrak{g}_{,k}$ denote the corresponding modular data:
$\Phi$ is the set $P_+^k(
\mathfrak{g}^{(1)})$ of integrable highest-weights of level $k$, the vacuum 
is $k\Lambda_0$, and explicit formulas for the matrices $S$ and $T$ are given
for instance in \cite[Ch.13]{kac}  or \cite{Gan2}.   Finite group, lattice and affine
algebra modular data {(and as we'll see next subsection, the double of a subfactor)} is always unitary.

%All modular data considered in this paper has the \textit{congruence subgroup property}:
%$T$ has finite order $N$, and the kernel of $\rho$ in (\ref{modrep}) contains the principal 
%congruence subgroup $\Gamma(N):=\{A\in{\rm SL}_2(\bbZ)
%\,:\, A\equiv I\ ({\rm mod}\ N)\}$. Moreover, $S_{ij}\in\bbQ[\xi_N]$ and $T_{\sigma i,\sigma i}
%=T_{i,i}^{\ell^2}$ where $\sigma\xi_N=\xi_N^{\ell^2}$. This property makes the calculations
%sketched in Section 1.3 much more manageable.

\medskip\noindent\textbf{Definition 2.} \textit{A matrix $Z=(Z_{i,j})_{i,j\in\Phi}$ is called a} {modular invariant} \textit{provided}

\noindent i) \textit{all entries $Z_{i,j}$ are nonnegative integers and $Z_{0,0}=1$;}

\noindent ii) \textit{$ZS=SZ$ and $ZT=TZ$.}\medskip

Examples are $Z=I$ and $Z=C$.
A modular invariant is often written in equivalent form as the formal expression
\begin{equation} Z=\sum_{i,j\in\Phi}Z_{i,j}ch_i\overline{ch}_j\,.\end{equation}
Unitary modular data  has only finitely many modular invariants.
The Galois symmetry \eqref{galois} implies, for any $\sigma\in\mathrm{Gal}(\bbQ[\xi_L]/\bbQ)$,
\begin{align}
&M_{i,j}=s_\sigma(i)\,s_\sigma(j)\,M_{\sigma i,\sigma j}\,,&\label{galsym}\\
&M_{i,j}\ne 0\ \mathrm{implies}\  s_\sigma(i)=s_\sigma(j)\,.&\label{Mselrule}\end{align}

%%%%%%%%
\subsection{Subfactors, quantum doubles and tube algebras}

We refer to \cite{J1,EK} for the basic theory of subfactors, principal graphs etc, \cite{Ln5} for the
theory of sectors, and \cite{BEK1,BE5,BE6} for the theory of alpha-induction.

Given a type III factor $N$,  let $\NXN$ denote a finite system of endomorphisms on $N$ {\cite[Defn 2.1]{BEK1}.}
Write   $\Sigma(\NXN)$ for the endomorphisms which decompose into a finite number of irreducibles
from $\NXN$. For  $\la, \rho\in\Sigma(\NXN)$,   the  intertwiner space
Hom$(\la,\rho)$ is a finite-dimensional Hilbert space. Write $\langle \la,\rho\rangle$
for its dimension. The \textit{sector} $[\lambda]$ identifies all endomorphisms
Ad$(u)\lambda$ for unitaries $u$ in the target algebra.

Suppose $\NXN$ is nondegenerately braided {\cite[Section 2.2]{BEK1}.} Among other things this means  $\lambda,\mu\in
\NXN$ commute up to a unitary $\epsilon=\epsilon(\lambda,\mu)$, i.e. $\lambda\mu=\mathrm{Ad}
(\epsilon)\mu\lambda$, and these
unitaries $\{\epsilon(\lambda,\mu)\}$ can be chosen to satisfy the {braiding--fusion relations}. Unitary modular
data $S,T$ is obtained from this set-up by the intertwiners associated to the Hopf link and
twist \cite{R1,tur}; here $\Phi=\NXN$ and $1$ is the identity endomorphism, and the corresponding fusion ring 
is realised by composition in $\NXN$: $[\lambda][\mu]=
\sum N_{[\lambda],[\mu]}^{[\nu]}[\nu]$ where $N_{[\lambda],[\mu]}^{[\nu]}=\langle [\lambda][\mu],
[\nu]\rangle$. 

Modular invariants are recovered as follows. Suppose we have a subfactor $N\subset M$.
Let $\iota:N\rightarrow M$ be the inclusion and $\overline{\iota}:M\rightarrow N$ its conjugate.
Then $\theta=\overline{\iota}\iota$ is called the \textit{canonical endomorphism} and $\gamma=
\iota\overline{\iota}$ its \textit{dual canonical endomorphism}. {Suppose $\theta$ is in 
$\Sigma(\NXN)$.}
Using the braiding $\epsilon^+:=\epsilon$
or its opposite $\epsilon^-:=\epsilon^{-1}$, we can lift an endomorphism $\lambda\in\NXN$
of $N$ to one of $M$: $\alpha^{\pm }_\lambda:=\gamma^{-1}\mathrm{Ad}(\epsilon^\pm(\lambda,\theta))
\lambda\gamma$. Then $Z_{\lambda,\mu}:=\langle \alpha^+_\lambda,\alpha^-_\mu\rangle$
is a modular invariant \cite{BEK1, Ev2}. The induced $\alpha^\pm(\NXN)=\MXM^\pm$
generate  the \textit{full system} $\MXM$. The $N$-$M$ system $\NXM$ (resp. 
the $M$-$N$ system $\MXN$)
consists of all irreducibles in $\iota\theta^n$ (resp. $\overline{\iota}\gamma^n$).
 By the \textit{nimrep} we mean the $\NXN$ action on the
$N$-$M$ system $\NXM$. This is one of 8 (6 independent) natural 
products ${}_P\cX_Q\times {}_Q\cX_R\rightarrow {}_P\cX_R$ among the sectors, one for each 
triple $P,Q,R\in\{M,N\}$.

Many examples (e.g. the Haagerup and finite groups) arise naturally
as \textit{nonbraided} systems of endomorphisms. To get a braided system, one takes the
quantum double (asymptotic inclusion)  of $\NXN$. This can be realised by the Longo-Rehren 
inclusion  $A\subset
B$  where $B=N\otimes N^{{\rm opp}}$ (see e.g. \cite{LR, iz1}): we are interested in the double $\cD(\NXN)$ 
system on $A$ (which contains the $A$-$A$ system of this subfactor, cf. Remark
(i) p.146 of \cite{iz1}) {whereas the full system ${}_B\cX_B$ is simply $\NXN \otimes {\NXN}^{\rm opp}$ (see Theorem 1 below) with dual 
canonical endomorphism $\gamma_{LR} = \Sigma_{\xi \in \NXN}  \xi \otimes \xi^{\rm opp}$.}

{If $\MXM$ is an $M$-$M$ system of a subfactor $N\subset M$, then the dual
Longo-Rehren inclusion is $\widehat{\iota}:A_1\subset M\otimes M^{\mathrm{opp}}=:B_1$,
where $\widehat{\gamma}=\overline{\widehat{\iota}}\widehat{\iota}=\sum_{\eta\in\MXM}\eta
\otimes\eta^{\mathrm{opp}}$. 
Here ${}_{B_1}\cX_{B_1}$ is $\MXM\otimes\MXM^{\mathrm{opp}}$
but we can and will identify the double $\cD(\NXN)$ on $A$ with the double
$\cD(\MXM)$ on $A_1$, and  $A_1$ with $A$.

 The induction-restriction graph of the
quantum double system is constructed in the following way
\cite{iz1}. First, given $\rho\in\Sigma(\NXN)$, a system of unitaries $\{\cE_\rho(\xi)\}$
is called a \textit{half-braiding} of $\rho$ if
$\cE_\rho(\xi)\in{\rm Hom}(\rho\xi,\xi\rho)$ for some
$\xi\in\NXN,$ and
$x\cE_\rho(\zeta)=\xi(\cE_\rho(\eta))\cE_\rho(\xi)\rho(x)\
 {\rm for\ every\ } x\in{\rm Hom}
(\zeta,\xi\eta),$
as in \cite[Defn.\ 4.2]{iz1}. 
This can be described by a matrix representation 
for  an orthonormal basis
$\{w_\rho(\xi)_i\}$ of Hom($\xi,\rho$), if we set
\begin{equation}\cE_\rho^p(\xi)_{(\eta,i),(\zeta,j)}:=
\rho_\xi(w_\rho(\zeta)_j^\ast)\cE_\rho^p(\xi) w_\rho(\eta)_i\in
{\rm Hom}(\eta\xi, \xi\zeta)\, .\nonumber\end{equation}
Then the even vertices of the quantum double system are labelled by 
inequivalent half-braidings $\cE_\rho^p$, and the odd vertices are labelled by
$(\xi\otimes{\rm id}^{\rm opp})\iota$ with $\xi$ running
in $\NXN$. Finally vertices $\cE_\rho$ and $(\xi\otimes{\rm
id}^{\rm opp})\iota$ are joined by $\langle \rho,\xi\rangle$ edges, see
\cite[page 154]{iz1}. The principal graph of $A\subset B$ is the
connected component containing $({\rm id}\otimes{\rm id}^{\rm
opp})\iota$. 
%In turn we automatically have the dual canonical
%endomorphism $\theta_{{\rm LR}}$ of $A\subset B$.

Incidentally, the forgetful functor sending the half-braiding $\cE_\rho^p$ to the representation
$\rho$ is an algebra homomorphism, and indeed has a name: it is alpha-induction.
See also the discussion at the end of the next subsection.

The tube algebra is an effective way to compute the modular data $S,T$ of the quantum double.
In Section 2 we compare the tube algebras of the Haagerup subfactor and the group $S_3$. Their 
similarity is the first indication of a possible relation between them. 
The \textit{tube algebra} of $\NXN$ is the finite dimensional C$^*$-algebra
\begin{eqnarray}
{\rm Tube}(\NXN)=\bigoplus_{\xi,\eta,\zeta\in\NXN}
{\rm Hom}(\xi\cdot
\zeta,\zeta\cdot \eta)\,.
\end{eqnarray}
The simple summands of Tube($\NXN$) are labelled by inequivalent half-braidings $\cE_\rho^p$.
The modular data $S,T$ is obtained from the half-braidings by (2.8),(2.9) of \cite{iz3}.

Let's consider in more detail the special case of tube algebras of finite groups.
Let $G$ be a finite group with identity $e$. Given a type III factor $N$, write
Int$(N)$ for the group of inner automorphisms and
Out$(N):={\rm Aut}(N)/{\rm Int}(N)$. Let $\a: G\to \hbox{Out}(N)$ be a homomorphism. We take a 
lift $\a: G\to {\rm Aut}(N)$, as in \cite{suth}. 

Since $[\a_g][\a_h]=[\a_{gh}]$ as sectors, there exists a unitary
$u_{g,h}\in\ N$ satisfying
$\a_g\cdot \a_h=\hbox{Ad}(u_{g,h})\a_{gh}.$
In particular $u_{g,h}\in\hbox{Hom}(\a_{gh},\a_g\cdot\a_h)$.
Associativity $(\a_g\a_h)\a_k=\a_g(\a_h\a_k)$ implies
$\hbox{Ad}(u_{g,h})u_{gh,k}\cdot \a_{ghk}=
\hbox{Ad}(\a_g(u_{h,k})u_{g,hk})\cdot \a_{ghk}.$
Then since $N$
is a factor, there exists a scalar $\omega(g,h,k)\in\bbT$ satisfying
$u_{g,h}u_{gh,k}=\omega(g,h,k)\a_g(u_{h,k})u_{g,hk},$
i.e. $\omega$ is an element in $Z^3(B_G;\bbT)$. Conversely, every element in $Z^3(B_G;\bbT)$
arises as such an $\omega$ {\cite{suth, J2}.}

Set ${\bf c}(g,h):= u_{h,g^h}u_{g,h}^\ast$. The product and *-structure of the tube
algebra Tube($\NXN$) are:
\begin{align}
&{\bf c}(g,h)\,{\bf c}(k,l)=\delta_{k,g^h}\,\theta_g(h,l)\,{\bf c}(g,hl)\,,\label{ccdelta}&\\
&{\bf c}(g,h)^\ast=\omega(gh,h^{-1},h)\,\omega(h^{-1},gh,h^{-1})\,
\omega(h^{-1},h,g^h)\,{\bf c}(g^h,h^{-1})&
\end{align}
(recall (\ref{theta})). Consequently, 
$${\bf c}(e,g)\,{\bf c}(e,h)={\bf c}(e,gh)\ \hbox{ and }\ {\bf c}(e,g)^\ast
={\bf c}(e,g^{-1})\,,$$
and thus the group algebra $\bbC[G]$ of $G$ is a C$^\ast$-subalgebra of
Tube($\NXN$).
The identity of the tube algebra is
$\sum_{g\in G}{\bf c}(g,e)$.

Let $K_g=\{g^h: h\in G\}$ denote the conjugacy class of an element
$g\in G$.  Then we have the decomposition
${\rm Tube}(\NXN)=\bigoplus_{K_u} {\rm Tube}(K_u)$
where the sum ranges over the conjugacy classes of $G$,
and Tube$(K_u)=\bigoplus_{g\in K_u,l\in G}{\rm
Hom}(\a_g\a_l,\a_l\a_{g^l})$. Hence $\bbC[G]={\rm Tube}(K_e)$.

This decomposition can be further refined. Consider first a trivial twist $\omega$ \cite{EP4}.
The outer action $\alpha$ gives a subfactor
$N\subset  N\rtimes G=M$. The $N$-$N$
system is $\{\a_g\}\equiv G$ whereas the $M$-$M$ system is the irreps
$\widehat{G}$.  {The tube algebras  Tube$(G)\simeq{\rm
Tube}(\widehat{G})$ are Morita equivalent} as finite-dimensional
C$^\ast$-algebras \cite{iz1}, {so their centres, the quantum doubles of $G$ and 
$\widehat{G}$, are identified}.
The simple components of ${\rm Tube}(G)$, equivalently the even vertices of the quantum
double of $G$,  are labelled by pairs $(K_u,\pi)$ where 
$\pi\in\widehat{\mathcal{Z}_u}$ are irreps of the centralizers
${\cal Z}_u$, see e.g.\
\cite{KY1} or \cite[Sect.\ 4]{EP}. To every such 
pair $(K_u,\pi)$ is attached the endomorphism
$\rho_{(K_u,\pi)}\in\Sigma(G)$ defined as
\begin{eqnarray}
\rho_{(K_u,\pi)}={\rm dim}(\pi)\sum_{h\in K_u}\a_h\,.
\end{eqnarray}

In the dual system $\widehat{G}$,
the endomorphism $\hat{\rho}_{(K_u,\pi)}$ in $\Sigma(\widehat{G})$ attached to
such a pair $(K_u,\pi)$ is the Mackey induction ${\rm Ind}_{{\cal Z}_u}^{{G}}\pi$
from ${\cal Z}_u$ to $G$, decomposed into irreps of $G$.

 The number of inequivalent half-braidings
associated to $\rho_{(K_u,\pi)}$ equals the number of 
inequivalent irreps $\pi^\prime\in\widehat{{\cal Z}_u}$ such that
${\rm dim}(\pi^\prime)={\rm dim}(\pi)$. The half-braiding for $\rho=\rho_{(K_u,\pi)}$ {
in its matrix decomposition}  is
\begin{equation}
{{\cal E}_{\rho}^\pi(\xi)}_{(u^p,i),(u^q,j)}=
\pi_{ji}(p\xi^{-1} q^{-1}){\rm id}_N\,, \label{halfbraiding}
\end{equation}
{for all $\xi \in p^{-1}{\cal Z}_uq, u^p, u^q \in K _u$, where $i,j$ label the basis vectors in
 the representation space of $\pi$.}
%the irreducible representation $\tilde{\pi}$ of ${G}$ --- possibly null.

%Of course we mau well have different pairs $(K,\pi)$ and $(K,\pi^\prime)$ with
%the same endomorphism
%$\rho_{(K,\pi)}=\rho_{(K,\pi^\prime)}$, but this amounts to the number
%of inequevalent half braidings for $\rho_{(K,\pi)}$ \cite{iz1}: the
%number of inequivalent half braidings for the
%endomorphism $\rho_{(K,\pi)}$ equals the number of
%$\pi^\prime\in\widehat{{\cal Z}}_g$ whose dim($\pi^\prime$)=dim($\pi$).
%The same applies to the $\widehat{\rho}_{(K,\pi)}$.

%Moreover, the dual canonical sectors $[\theta_{{\rm LR}(G)}]$ and
%$[\theta_{\rm LR(\widehat{G})}]$ of the Longo-Rehren inclusions of
%the systems $G$ and $\widehat{G}$ are, respectively,
%\begin{eqnarray}\label{dualendgroup}
%[\theta_{{\rm LR}(G)}]=\bigoplus_{\pi\in \widehat{G}}{\rm
%dim}(\pi)[\a_e]\,,\qquad [\theta_{\rm LR(\widehat{G})}]={\rm
%Ind}_{id}^{\widehat G}\,. 
%\end{eqnarray}

%\end{remark}

All this extends to the twisted case.
The $\theta_a$'s of (\ref{theta}) are normalised twisted
cocycles
on $G$, namely from (\ref{ccdelta}) they satisfy
$$
\theta_a(x,y) \, \theta_a(xy,z) = \theta_a(x,yz) \, \theta_{x^{-1}ax}(y,z)\,,
\qquad \forall x,y,z \in G\,.
\label{cocy}
$$
They thus describe  projective representations of $\mathcal{Z}_a$.
As we know \cite{DW}, the primary fields  in the model twisted by a
given 3-cocycle $\omega$ consist of all pairs $(K_a,\tilde \pi)$ where
%$a \in \T$ (as before) and $\tilde \chi\in\beta_a$-Irr$(C_G(a))$.
now $\tilde \pi\in\theta_a$-Irr$(\mathcal{Z}_a)$.

In the cohomologically trivial case discussed in Section 1.1, this is immediate: we can twist
 the formula (\ref{halfbraiding}) for half-braidings by $\epsilon$; inserting into (2.8),(2.9)
 of \cite{iz3} recovers (\ref{moddataG}).

\subsection{The canonical and dual canonical modular invariants}

This subsection contains original material.

Suppose we have a subfactor $N\subset M$, with $N$-$N$ system denoted by
$\Delta$ and $M$-$M$
system $\widehat{\Delta}$. The double $\cD\Delta=LR(\Delta)$ of $\Delta$ is realised as sectors on a
factor $A$ in $N\otimes N^{\mathrm{opp}}=:B$. Alpha-induction reverses this. That is, if you
take the appropriate modular invariant $Z_\Delta$ in the double, and take alpha-induction
from $A$ to $B$, then the full system is $\Delta\otimes\Delta^{\mathrm{opp}}$. In some sense
then, the factor $A$ and the modular invariant $Z_\Delta$ remembers the original system
$\Delta$ on $A$. In fact it is a consequence of Theorem 1 below that the $A$-$B$ system is just
 $\Delta$ regarded as a nimrep. 

To get $B$, only the canonical endomorphism $\theta$ on $A$ is needed { 
(because of Theorem 1 below, 
we can read $\theta$ off from} the modular invariant), plus the associated Q-system on $\theta$.
That is, all the information about $B$ and the full system $\Delta\otimes\Delta^{\mathrm{opp}}$
is carried in the primary fields and modular data of the double, the modular invariant $Z_\Delta$
and its vacuum block $\theta$. We suggest this modular invariant $Z_\Delta$ be called the
\textit{canonical modular invariant}.

Dually, there is an inclusion $A\subset M\otimes M^{\mathrm{opp}}=:B$, where the dual
$\widehat{\Delta}$ acts on $M$, but the same double acts on $A$. This comes with another
modular invariant $Z_{\widehat{\Delta}}$, with a corresponding $\widehat{\theta}$ which
regains $\widehat{\Delta}\otimes\widehat{\Delta}^{\mathrm{opp}}$ as the full system. Here
(again from Theorem 1 below) the $A$-$B$ system is just $\widehat{\Delta}$. Again the modular invariant
 $Z_{\widehat{\Delta}}$ --- the \textit{dual canonical modular invariant} --- and $\widehat{\theta}$ 
 together encode the original system $\widehat{\Delta}$.
 
 This means we should regard the double $\cD \Delta\cong \cD\widehat{\Delta}$ as equipped
 with two canonical modular invariants $Z_\Delta$ and $Z_{\widehat{\Delta}}$, from
 which using say alpha-induction one can recover the full system $B$-$B$  and the nimrep 
 $A$-$B$ entirely in terms of $\Delta$ or $\widehat{\Delta}$. In particular the full system is
 $\Delta\otimes\Delta^{\mathrm{opp}}$ or $\widehat{\Delta}\otimes
 \widehat{\Delta}^{\mathrm{opp}}$, depending on the choice of modular invariant
$Z_\Delta$ or $Z_{\widehat{\Delta}}$, and
 the nimrep is $\Delta$ or $\widehat{\Delta}$. 
 
 Note that in both cases the $A$-$B$ system here is an algebra. 
 For most modular invariants the $A$-$B$ system is only a module, but when the inclusion
 is type I, which implies for instance that the modular invariant
 is a sum of squares, the $A$-$B$ system will necessarily be an algebra ($A\subset B$ is
 type I iff the equivalent conditions of Proposition 3.2 in \cite{BE4} hold; in the nets of subfactors 
 setting $A(I)\subset B(I)$, this means the extended net is local). This suggests
that both inclusions $A\subset N\otimes N^{\mathrm{opp}}$ and $A\subset M\otimes M^{\mathrm{opp}}$ are type I. In fact more is true:

\medskip\noindent\textbf{Theorem 1.} \textit{The inclusions $A\subset N\otimes N^{\mathrm{opp}}$ 
and $A\subset M\otimes M^{\mathrm{opp}}$ are both type I. Moreover, we have 
$(Z_\Delta)_{i,j}=(Z_\Delta)_{i,0}(Z_\Delta)_{j,0}\,$, i.e.} $Z_\Delta=\left|\sum_i Z_{i,0}ch_i\right|^2$
\textit{(similarly for  $Z_{\widehat{\Delta}}$). 
The vectors $u$ and $v$ with entries $u_i=(Z_\Delta)_{i,0}$ and 
$v_i=(Z_{\widehat{\Delta}})_{i,0}$ are eigenvectors with eigenvalue 1 of both $S$ and $T$.}

\medskip\noindent\textit{Proof.} {Corollaries 6.3 and 6.4 of \cite{iz1} show that, at the level
of subfactors, $\BXB^+:=\alpha^+(\BXB)\subseteq \Delta\otimes 1$ and $\BXB^-:=
\alpha^-(\AXA)
\subseteq 1\otimes{\Delta}^{\mathrm{opp}}$. The neutral system $\BXB^0=\BXB^+\cap\BXB^-$ is thus contained in
$\Delta\otimes 1\cap 1\otimes\Delta^{\mathrm{opp}}$ and hence equals 1.
But dim$\,\BXB^\pm$ is computed in \cite[Prop.3.1]{BE4},  for such a nondegenerate 
system like ours, to be 
\begin{equation}
\frac{1}{S_{0,0}}=\sum\frac{S_{i,0}}{S_{0,0}}Z_{i,0}=\sum\frac{S_{j,0}}{S_{0,0}}Z_{0,j}\,,
\label{dimtheta}\end{equation}
which matches the dimensions of $\Delta\otimes 1$ and $1\otimes \Delta^{\mathrm{opp}}$.
Therefore we get the equalities $\BXB^+=\Delta\otimes 1$ and $\BXB^-=1\otimes
\Delta^{\mathrm{opp}}$.

Any canonical endomorphism $\theta$ is bounded below by the vacuum sector:
$\theta\ge \sum_\lambda Z_{\lambda,0}\lambda$. This inequality is a special case of
$\langle\theta\lambda,\mu\rangle\ge \langle\alpha^\pm_\lambda,\alpha^\pm_\mu\rangle$ (\cite[eq.(37)]{BEK2}; see also
\cite[Thm.3.9]{BE1}). But dim$\,\theta$ is given by
\begin{equation}\mathrm{dim}\,\theta=\mathrm{dim}(\overline{\iota}\iota)=\mathrm{dim}(\iota\overline{\iota})
=\sum\mathrm{dim}(\xi\otimes\xi^{\mathrm{opp}})=\sum(\mathrm{dim}\,\xi)^2=\frac{1}{S_{0,0}}\,,
\nonumber\end{equation}
agreeing with 
\eqref{dimtheta}, and thus $\theta=\sum_\lambda Z_{\lambda,0}\lambda$. By 
\cite[Prop.3.2]{BE4}, this implies the inclusion $A\subset B$ is type I.

Because of alpha-induction, we know a modular invariant takes the form $Z=\sum_{\tau,\lambda,
\mu}b^+_{\tau,\lambda}b_{\tau,\mu}$, where $\tau$ runs over all sectors of the neutral system
and $b^\pm_{\tau,\lambda}=\langle \tau,\alpha^\pm_\lambda\rangle$ are the branching
coefficients. Because this system is type I, we have $b^+=b^-$.
Because the neutral system is trivial, we have only one $\tau$. Thus the modular invariant
$Z_\Delta$ corresponding to $\theta$, and $Z_{\widehat{\Delta}}$ corresponding to
$\widehat{\theta}$, take the desired forms. The statements about the vectors $u,v$ follow
from modular invariance.}
QED\medskip

From the theory of alpha-induction, we know that $\Delta$ (respectively $\widehat{\Delta}$)
is noncommutative iff some $(Z_\Delta)_{i,j}>1$ (resp. some $(Z_{\widehat{\Delta}})_{i,j}>1$)
--- see \cite[Cor.6.9]{BEK1},\cite[Thm.4.11]{BEK2}.
For type I inclusions, we know from Section 4.1 of \cite{BE3} that the $A$-$B$ system
$\AXB$ is isomorphic to $\BXB^\pm$.

We will call any modular invariant of the type $Z=|\sum_i Z_{i,0}\,ch_i|^2$ a \textit{monomial
modular invariant}.

Of course all this applies for a system $\Delta$ of endomorphisms not
necessarily coming from a subfactor --- we expect this is relevant for our
twisted Haagerup data. Again $\Delta$ is recovered from the
double $\cD\Delta$ and the canonical modular invariant $Z_\Delta$, which again 
is monomial. The only 
difference is that there is no $\widehat{\Delta}$.

Suppose a finite group  $G$ acts by outer automorphism on a type III
factor $N$. Then as mentioned last subsection, by taking a crossed product we have a subfactor $N\subset M$ where the
$N$-$N$ and $M$-$M$ systems are identified with $G$ and $\widehat{G}$ respectively.
Taking the Longo-Rehren inclusion $A\subset N\otimes N^{\mathrm{opp}}$, the doubled $A$-$A$
system can be identified with the untwisted quantum double $\cD G$. We will sometimes find it useful to use the
K-theory language developed in \cite{Ev,EG1}, which identifies $\cD G$ with
the equivariant K-group $K_{G}^0(G) \cong K^0_{\Delta-\Delta}(G\times G)$ where
$G$ acts on $G$ by the conjugate action and $\Delta=\Delta(G)\cong G$ acts 
diagonally on the left and right of $G\times G$.
The neutral system is trivial
($K^0(e,e)= \bbZ$), with sigma-restriction $\sigma_{\mathrm{id}}=\theta_{\cD G}\cong
\sum_{\chi\in\widehat{G}}\mathrm{dim}\,\chi\,(e,\chi)\,$. The full
system is $\Delta\times\Delta^{\mathrm{opp}}\cong K^0(G \times G)$ with canonical modular invariant
\begin{equation}\label{gpmodinv1}
Z_\Delta=\left|\sum \chi(e)\,ch_{(e,\chi)}\right|^2\,,\end{equation}
the sum over $\chi\in\widehat{G}$.  The dual LR-inclusion is
$A\subset M\otimes M^{\mathrm{opp}}$, where the doubled $A$-$A$ system is again the quantum
double of $G$, the full system is $\widehat{G}\times\widehat{G}^{\mathrm{opp}}
\cong K^0_{G \times G}(G/G \times G/G)$ and the neutral system $K^0(G/G \times G/G)$ is trivial. Here
sigma-restriction is given by  $\widehat{\sigma}_{\mathrm{id}}=\widehat{\theta}_{\cD G}
\cong \sum (g,\mathrm{id})$ where $g$ runs over representatives
of all conjugacy classes, and the dual canonical modular invariant is
\begin{equation} \label{gpmodinv2}
 Z_{\widehat{\Delta}}=\left|
\sum ch_{(a,1)}\right|^2\,,\end{equation}
where the sum runs over representatives $a$ of all conjugacy classes of $G$.
 We recover from \eqref{gpmodinv1},\eqref{gpmodinv2} that $\widehat{\Delta}=\widehat{G}$
is always commutative, whereas $\Delta=\bbC G$ is commutative iff all dimensions $\chi(e)=1$,
i.e. iff $G$ is abelian.

By the Galois theory of \cite{iz1}, a subgroup $K<G$ induces an intermediate
subfactor $A\subset C\subset N\otimes N^{\mathrm{opp}}$ where the doubled $C$-$C$ system
is the untwisted quantum double $\cD K\cong K^0_{K}(K)\cong K^0_{\Delta K-
\Delta K}(K \times K)$.
Then $A\subset C$ is a braided subfactor of index $|G/K|$. The full system ${}_C\mathcal{X}_C\cong K^0_{\Delta K
-\Delta K}(G\times G)$ and chiral systems ${}_C\mathcal{X}_C^\pm
\cong {}_C\cX_A\cong K^0_{\Delta K-\Delta G}(G \times G)$. The branching coefficients and
sigma-restriction
$K^0_{K}(K)\rightarrow K^0_{G}(G)$ are given by $(g,\pi)\rightarrow \Sigma(g',
\mathrm{Ind}^{Z_g}_{K \cap Z_g} \pi)$, where the summation is over equivalence classes of $g'$ 
in $G$ conjugate to $g$. The trivial example $K=1$ recovers the canonical modular invariant
\eqref{gpmodinv1}. In the case of the dihedral group $G=D_\nu$ and its cyclic subgroup
$K=\bbZ_\nu$, in which we are interested,
the canonical modular invariant arising from the inclusion is  given in \eqref{modinvDnuZnu}
below.

\subsection{VOAs and  vector-valued modular functions}

For the basic theory of VOAs, see \cite{LeLi}. Let $\cV$ be a VOA. Write Irr$(\cV)$ for the set of 
irreducible $\cV$-modules $M$. Among other things, $\cV$ and its modules carry a representation
of the Virasoro algebra $\mathrm{Vir}=\mathrm{Span}\{L_n,C\}_{n\in\bbZ}$, where the central term
$C$ acts as a scalar $c=c(\cV)$ called the \textit{central charge}. $c$ also
arose in Section 
2.1, where it was defined only mod 24. The action of 
$L_0\in Vir$ on each $\cV$-module $M$ defines a grading (by eigenvalues) on $M$ into finite-dimensional spaces.

\medskip\noindent{}\textbf{Definition 3.} \textit{A VOA $\cV$ is called} {rational} \textit{if}

\noindent{i)} $\cV\in\mathrm{Irr}(\cV)$ \textit{and  $\cV$ is isomorphic to its
contragredient as $\cV$-modules;}

\noindent{ii)} \textit{writing $\cV=\oplus_{n\in\bbZ}\cV_n$ for the grading by $L_0$ of $\cV$, we have
 $\cV_n=0$ for $n<0$ and $\cV_0$ is 1-dimensional;}
 
\noindent{iii)} \textit{every weak $\cV$-module is completely reducible.}\medskip

 See \cite{LeLi} for the details, which play no role in the following. There is no standard 
 definition of rationality --- we chose Definition 3 to guarantee the existence of modular data.
 All VOAs considered in this paper are rational in this sense.
 
Rational VOAs $\mathcal{V}$ (or if you prefer, rational CFT) realise
modular data as follows. The primaries consist of  the finitely many modules 
$M\in\mathrm{Irr}(\cV)=\Phi$, and the
vacuum 0 is $\cV$ itself. Define their \textit{characters} 
 to be the graded dimensions
\begin{equation}ch_M(\tau)=\mathrm{Tr}_Mq^{L_0-c/24}\,,\label{voachar}\end{equation}
using the above grading by $L_0$, where as always we write $q=e^{2\pi\i\tau}$. These 
$ch_M$ will be holomorphic throughout the upper half-plane $\mathbb{H}=\{\tau\in
\bbC:\mathrm{Im}(\tau)>0\}$ \cite{zhu}. Collect these 
finitely many characters into a column vector $\vec{ch}(\tau)$. Then \cite{zhu} showed
there is a representation $\rho$ of SL$_2(\bbZ)$ such that
\begin{equation}\label{cheq}
\vec{ch}\left(\frac{a\tau+b}{e\tau+d}\right)=\rho\left(\begin{matrix}a&b\cr e&d\end{matrix}\right)
\vec{ch}(\tau)\,,\ \forall \left(\begin{matrix}a&b\cr e&d\end{matrix}\right)\in\mathrm{SL}_2(\bbZ)\,.
\end{equation}
The matrices $S,T$ now defined by (\ref{modrep}) constitute modular data \cite{Hua}. It
may or may not be unitary. Unlike braided subfactors, which are naturally associated to modular invariants, 
VOAs see only one chiral half of the rational CFT, and
so  capture only the notion of modular data.

There is a rational VOA $\cV(L)$ corresponding to any even positive-definite lattice $L$
(the central charge equals the dimension $n$ of $L$), recovering the modular
data of (\ref{moddataL}). The character corresponding to coset $[v]\in L^*/L$ is
$ch_{[v]}(\tau)=\eta(\tau)^{-n}\sum_{x\in[v]}q^{x\cdot x/2}$.
There is a rational VOA $\cV(\mathfrak{g}_k)$ corresponding to
affine algebra $\mathfrak{g}^{(1)}$ at positive integral level $k$, recovering the modular
data of \cite{kac}. The character $ch_\lambda$ corresponding to highest-weight $\lambda$ coincides with
the affine algebra character $\chi_\lambda$, specialised to $\tau\in\mathbb{H}$.
Affine algebra and lattice VOAs overlap for the simply-laced $\mathfrak{g}$
at level 1, corresponding to the root lattices $L$. Finite group modular data $\cD^\omega G$ is 
recovered by taking the orbifold
$\cV^G$ (i.e. subVOA of fixed-points of $G$) of a subgroup $G$ of automorphisms of a
\textit{holomorphic} VOA (i.e. a rational VOA $\cV$ with Irr$(\cV)=\{\cV\}$) --- it has been conjectured that any such $\cV^G$
will itself be rational, but a general proof of this still seems far away.   Examples of holomorphic
VOAs  are the $\cV(L)$ for $L$ self-dual --- e.g. for $c=8$ there
is the $E_8$ root lattice. In this case the orbifold theory is under better control, see e.g.
\cite{KacT}, and indeed this is the main one we'll consider. A sexier example of holomorphic VOA 
though is the $c=24$ \textit{Moonshine module} $V^\natural$.

Incidentally, Theorem 1 suggests that: \textit{a VOA corresponding to
a quantum double should be a conformal subalgebra of a holomorphic VOA}.

\medskip
\noindent\textbf{Definition 4.} \textit{Let $\rho$ be a $d$-dimensional representation of 
SL$_2(\bbZ)$, with diagonal  $T:=\rho({1\,1\atop 0\,1})$. A} {weakly holomorphic
vector-valued modular
 function} \textit{$\vec{ch}(\tau)$ with  multiplier $\rho$ is a holomorphic function
$\mathbb{H}\rightarrow\bbC^d$ satisfying (\ref{cheq}), with $q$-expansion}
\begin{equation}\label{qexp}\vec{ch}(\tau)=q^\lambda \sum_{n=0}^\infty \vec{ch}_nq^n\end{equation}
\textit{for some diagonal matrix $\lambda$ satisfying $T=e^{2\pi\i\lambda}$,
where each Fourier coefficient $\vec{ch}_n\in\bbC^d$ is independent of $\tau$. Let $\cM(\rho)$
denote the space of all such $\vec{ch}$ for $\rho$.}\medskip

Equation (\ref{qexp}) means $\vec{ch}(\tau)$ is meromorphic at the cusp.
\cite{zhu} showed that the character vector of a rational VOA is a vector-valued modular function
in this sense. These character vectors  will help us study and identify the VOA.

The characters of rational VOAs also satisfy other conditions.
For instance, by definition (\ref{voachar}) the Fourier coefficients $\vec{ch}_n$ in (\ref{qexp})
are nonnegative integers. The \textit{vacuum} character 
$ch_{\cV}=ch_0$ begins with $1q^{-c/24}$. In a unitary VOA,  $\lambda$ 
in (\ref{qexp}) can be chosen so that $\la_M>-c/24$ for all $M\ne \cV$ in 
$\mathrm{Irr}(\cV)$.

Recall charge-conjugation $C=S^2$ from Definition 1. A consequence of (\ref{cheq}) is
that $C\vec{ch}=\vec{ch}$. For convenience we assume in the remainder of this subsection
that $C=I$, i.e. we have a representation
of PSL$_2(\bbZ)$. This holds for most modular data considered in this paper
(other than $\cD^\omega\bbZ_{2n+1}$), but at the end of  Section 2.4 we explain
 how to reduce $C\ne I$ to the $C=I$ case.

Given any $d$-dimensional PSL$_2(\bbZ)$-representation $\rho$,  \cite{BG1,BG2} 
describe how to find  all $\vec{ch}(\tau)\in\cM(\rho)$. 
In particular, there is a $d\times d$ matrix $\Xi(\tau)$ with the
property that $\vec{ch}(\tau)\in\cM(\rho)$  iff
$\vec{ch}(\tau)$ is of the form $\vec{X}(\tau)=\Xi(\tau) \vec{P}(J(\tau))$, where $\vec{P}(x)$ is any 
column vector whose entries are polynomials in $x$, and $J(\tau)$ is the Hauptmodul 
$J(\tau)=q^{-1}+196884q+\cdots$.  Thus it suffices to find 
$\Xi(\tau)$ for $\rho$.  
$\Xi(\tau)$ is determined through a differential equation it satisfies. This differential equation 
and the relevant initial condition is determined from two $d\times d$ matrices
of numbers. One is  a diagonal matrix $\Lambda$ satisfying among other things the relation 
$e^{2\pi\i\Lambda}=T$; the more elusive one is called $\chi$. In terms of $
\Lambda$ and $\chi$, the matrix $\Xi(\tau)$  has $q$-expansion
$$\Xi(\tau)=q^\Lambda\sum_{n=-1}^\infty q^n\Xi_n=q^\Lambda(Iq^{-1}+\chi+\sum_{n=1}^\infty q^n\Xi_n)\,,$$
where $\Xi_{-1}$ is the $d\times d$ identity matrix, $\Xi_0$ is 
the matrix $\chi$,  and for $n=1,2,3,\ldots$ the matrix $\Xi_n$ is recursively 
defined by the commutator 
\begin{equation}\label{recurs}[\Lambda,\Xi_n]+(n+1)\Xi_n=\sum_{l=-1}^{n-1}\Xi_l(f_{n-l}
(\Lambda-I)+g_{n-1}(\chi+[\Lambda,\chi]))\,,\end{equation}
where $f_n,g_n$ are defined by
$(J-240)\Delta/E_{10}=\sum f_nq^n$ and 
$\Delta/E_{10}=\sum g_nq^n$ for the discriminant form
$\Delta=\eta^{24}$ and Eisenstein series $E_{10}=E_4E_6$. That is, the 
$(i,j)$-entry of
$\Xi_n$ is the $(i,j)$-entry of \eqref{recurs}, divided by $\Lambda_{ii}-\Lambda_{jj}+n+1$.

The following method suffices to find $\Lambda,\chi$ for the $\rho$ considered in this paper (though invariably there
are more elegant methods). First, decompose $\rho$ into irreps. What is 
special about irreps $\rho$ is that their $\cM(\rho)$  is cyclic (Theorem 4.1 of \cite{BG2}): the space $\cM(\rho)$
 is a module over the ring  $\bbC[J,\nabla_1,\nabla_2,\nabla_3]$ of differentiable operators
 $$\nabla_1=\frac{E_{10}}{\Delta}\mathrm{D}_0\,,\ \nabla_2=\frac{E_8}{\Delta}
 \mathrm{D}_1\mathrm{D}_0\,,\
\nabla_3=\frac{E_6}{\Delta}\mathrm{D}_2\mathrm{D}_1\mathrm{D}_0
\,,\ \mathrm{for}\ \mathrm{D}_k=q\frac{\mathrm{d}}{\mathrm{d}q}-\frac{k}{6}E_2\,,$$
  and if $\rho$ is irreducible, then $\cM(\rho)$
(and hence $\Lambda$ and $\chi$) is generated over that ring by any nonzero $\vec{ch}\in
\cM(\rho)$.  All irreps occurring in this paper are subrepresentations of the
modular data coming from even lattices, and so the desired nonzero modular function can be built
from lattice theta functions. 

Knowing $\Lambda$ and $\chi$ for a given $S,T$ is equivalent to knowing $\Lambda$ and
$\chi$ for $S,\omega T$ for any third root of unity $\omega$, but the explicit equivalence 
\cite{BG2} is not easy.  Write $\Xi^{(k)}(\tau)$ for the matrix $\Xi(\tau)$ for $S,\xi_3^kT$ and assume
$\Xi^{(0)}(\tau)$ (hence $\Lambda^{(0)}$ and $\chi^{(0)}$) is known. Then the columns of 
${\Xi}^{(1)}$ are linear combinations over $\bbC$ of the columns of
$\eta^{-16}E_4^2{\Xi}^{(0)}$ and
$\eta^{-16}\left(E_6\mathrm{D}_0{\Xi}^{(0)}-E_4^2{\Xi}^{(0)}
(\Lambda^{(0)}-1)\right)$, while the columns of $\Xi^{(2)}$ 
are linear combinations over $\bbC$ of $\eta^{-8}E_4{\Xi}^{(0)}$ and
$\eta^{-8}\left(\mathrm{D}_1\mathrm{D}_0{\Xi}^2-E_4{\Xi}^{(0)}(\Lambda^{(0)}-I)(\Lambda^{(0)}-\frac{7}{6}I)
\right)\,.$

The contragredient $(\rho^T)^{-1}$ is handled
similarly \cite{BG2}: the columns of $E_{10}E_{14}\Delta^{-2}(\Xi^T)^{-1}$ are linearly
independent vectors in $\cM((\rho^T)^{-1})$.

\section{Comparing the Haagerup, $Sym(3)$ and SO(13)}

%%%%%%%%%%%%%%%%%%%%%%%%
\subsection{The tube algebra of $Sym(3)$}

Recall the discussion in Section 1.2.
Write $S_3=Sym(3)=\{e,u,u^2,\tau, \tau u, \tau u^2\}$.
The three conjugacy classes are $K_e=\{e\}$, $K_u=\{u,u^2\}$
and $K_\tau=\{\tau, \tau u, \tau u^2\}$. Write $\widehat{S}_3=\{1,\epsilon,\sigma\}$ where $\epsilon$
is $sgn$ and $\sigma$ is 2-dimensional. Then each pair $(K_g,\pi)$, with $\pi\in 
\widehat{{\cal Z}}_g$, gives rise  to a simple component of the tube algebra.

There are two half-braidings for endomorphism $\rho=e$ as
$\rho_{(K_e,{\bf 1})}=\rho_{(K_e,\epsilon)}=e$, and a third
attached to $\rho_{(K_e,\sigma)}=e+e$. We denote these half-braidings by 
$e^{(1)},e^{(2)}$ and $2e$, respectively. 
The conjugacy class $K_u$ provides an endomorphism 
$\rho_{(K_u,\pi)}=u+{u^2}$ for every irrep $\pi$  of 
the centraliser ${\cal Z}_u\simeq \bbZ_3$, and hence
three half-braidings $(u+ u^2)^{(1)}$, $(u+
u^2)^{(2)}$, $(u+ u^2)^{(3)}$.
The conjugacy class $K_\tau$ has 
centraliser ${\cal Z}_\tau\simeq \bbZ_2$, so the endomorphism
$\rho=\tau+{\tau u}+{\tau u^2}$
has two half-braidings $(\tau+ \tau u+ \tau u^2)^{(1)}$,
$(\tau+ \tau u+\tau u^2)^{(2)}$.

\medskip\epsfysize=3in\centerline{ \epsffile{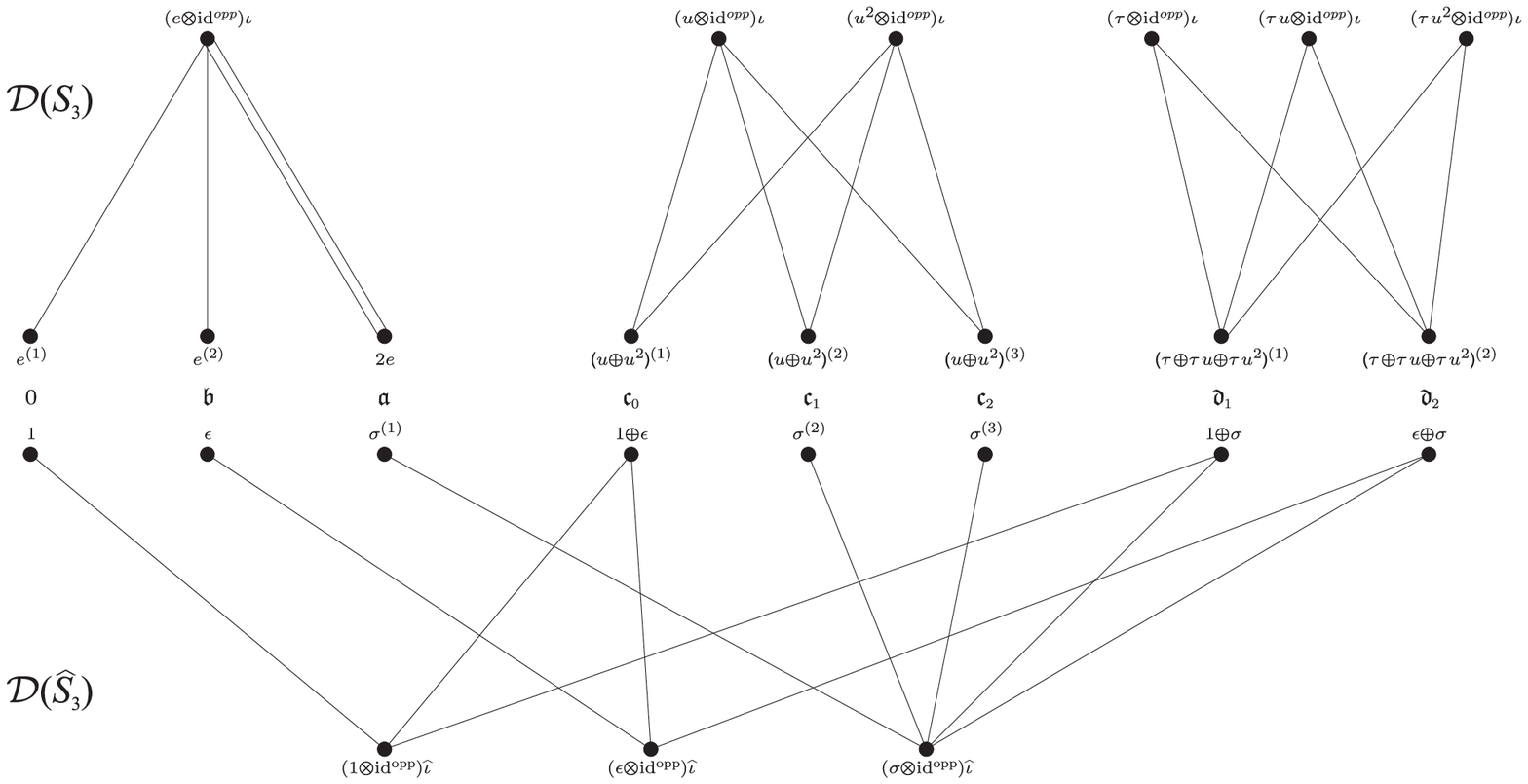}}\medskip

\centerline{Figure 1. Dual principal graphs for doubles of $Sym(3)$ and $\widehat{Sym(3)}$}
\label{s3fig}\medskip

This enables us to match up Figure 1, where the bottom graph
has been drawn in \cite[Fig.\ 1]{iz3} or \cite[Fig.\ 31]{EP}: the
upper graph arises from the Longo-Rehren inclusion of $S_3$, and
the lower describes the Longo-Rehren inclusion of the dual system
$\widehat{S}_3$, see \cite[Remark (i) in Page 154]{iz1} as well as Section 1.2 above for the
general description of the induction-restriction graphs between
$A$-$B$ and $B$-$B$ sectors of the Longo-Rehren inclusion
$A\subset B$ using the structure of the tube algebra.  The middle
vertices in Figure 1 also describe how the half-braidings from
$S_3$ and $\widehat{S}_3$ match up.

%Moreover, the dual canonical sectors $[\theta_{{\rm LR}(G)}]$ and
%$[\theta_{\rm LR(\widehat{G})}]$ of the Longo-Rehren inclusions of
%the systems $G$ and $\widehat{G}$ are, respectively,
%\begin{eqnarray}\label{dualendgroup}
%[\theta_{{\rm LR}(G)}]=\bigoplus_{\pi\in \widehat{G}}{\rm
%dim}(\pi)[\a_e]\,,\qquad [\theta_{\rm LR(\widehat{G})}]={\rm
%Ind}_{id}^{\widehat G}\,. 
%\end{eqnarray}
% By Eq.\
%(\ref{dualendgroup}), 

{The Longo-Rehren dual sectors associated to
$S_3$ and $\widehat{S_3}$ are respectively
\begin{eqnarray}&&[\theta_{\cD{S_3}}]=  [K_e, {\rm
Ind}_{\mathrm{id}}^{\widehat G}]\,=  \sum_{\pi\in \widehat{G}}{\rm
dim}(\pi)[K_e, \pi]\, = [\epsilon^{(1)}]+[\epsilon^{(2)}]+ 2[2\epsilon]\,,\quad\\
&&\quad [\theta_{\rm
\cD\widehat{S}_3}]= \sum_g [K_g, {\rm
Ind}_{K_g}^{ G}\mathrm{id}]\,  = [1]+[1+ \epsilon]+[1+\sigma]\,,\end{eqnarray}
specialised from the discusson of Section 1.3. Observe
that $[\theta_{\cD\widehat{S}_3}]$ has also been computed in
\cite{iz3,EP}. These braided subfactors yield the canonical modular invariants
(recall \eqref{gpmodinv1},\eqref{gpmodinv2}) $Z_{22} = |ch_0 + ch_{\frak{b}} + 2 ch_{\frak{a}}|^2$
and  $Z_{55} = |ch_0 + ch_{\frak{c}_0}+  ch_{\frak{d}_1}|^2$ with full systems
$S_3 \times S_3^{\mathrm{opp}}$, $\widehat{S_3} \times \widehat{S_3}^{\mathrm{opp}}$ 
respectively by \cite[page 357]{EP},} where the names of the primaries of $\cD S_3$ are taken 
from the middle row of { Figure 1} and the names $Z_{22},Z_{55}$ come from the list 
in $\cite{EP}$. 

As explained at the end of Section 1.3, the subsystem $\bbZ_3$ of $\Delta$ corresponds to
an intermediate subfactor $A\subset C\subset N\otimes N^{\mathrm{opp}}$
 where the dual canonical endomorphism of $C\subset N\otimes N^{\mathrm{opp}}$ is
 $\gamma=\sum_{\alpha\in\bbZ_3}\alpha\otimes\alpha^{\mathrm{opp}}$ and the
 double $C$-$C$ system is the quantum double $\cD\bbZ_3\cong K^0_{\bbZ_3}(\bbZ_3)$.
 Then $A\subset C$ is a braided subfactor of index $2$, with canonical
 endomorphism $[\theta]=[0]+[\mathfrak{b}]$ and associated modular invariant
\begin{equation} Z=|ch_0+ch_\mathfrak{b}|^2+2|ch_\mathfrak{a}|^2+2|ch_{\mathfrak{c}_0}|^2+
2|ch_{\mathfrak{c}_1}|^2+2|ch_{\mathfrak{c}_2}|^2\,.\label{Z3S3}\end{equation}

%%%%%%%%%%%%%%
\subsection{The Haagerup tube algebra}
\label{s3haag}

Let $\Delta=\{\hbox{id},u,u^2,\rho,\rho u,\rho u^2\}$ be
the noncommutative $N$-$N$ system (the even
vertices of the top graph of Figure 2)  of the Haagerup subfactor $N\subset M$
of index
 $\delta+1$, where $\delta=(3+\sqrt{13})/2$. The fusions (product of sectors)
are given by:
$$[u]^3=[\hbox{id}]\,,\quad [u][\rho]=[\rho][u]^2\,,\quad
[\rho]^2=[\hbox{id}]+[\rho]+[\rho u]+[\rho u^2]\,.$$
These have statistical dimensions dim$(u)=1$ and dim$(\rho)=\delta$. 
Let $\widehat{\Delta}=\{\hbox{id},a,b,c\}$ be the commutative $M$-$M$
system (the even vertices of the bottom graph of Figure 2), whose fusion rules are:
\begin{align}&&[a]^2=[\hbox{id}]+[a]+[b]+[c]\,,\quad [b]^2=
[\hbox{id}]+[c]\,,\quad [c]^2=[\hbox{id}] +
2[a]+[b]+ 2[c]\,,\nonumber\\
&&[a][b]=[a]+[c]\,,\quad [a][c]=[a]+ [b]+ 2[c]\,,\quad
[b][c]=[a]+[b]+[c]\,.\nonumber\end{align}
Hence we obtain the statistical dimensions  dim($a)=\delta$, 
dim($b)=\delta-1$, dim($c)=\delta+1$.

\medskip\epsfysize=2.5in\centerline{ \epsffile{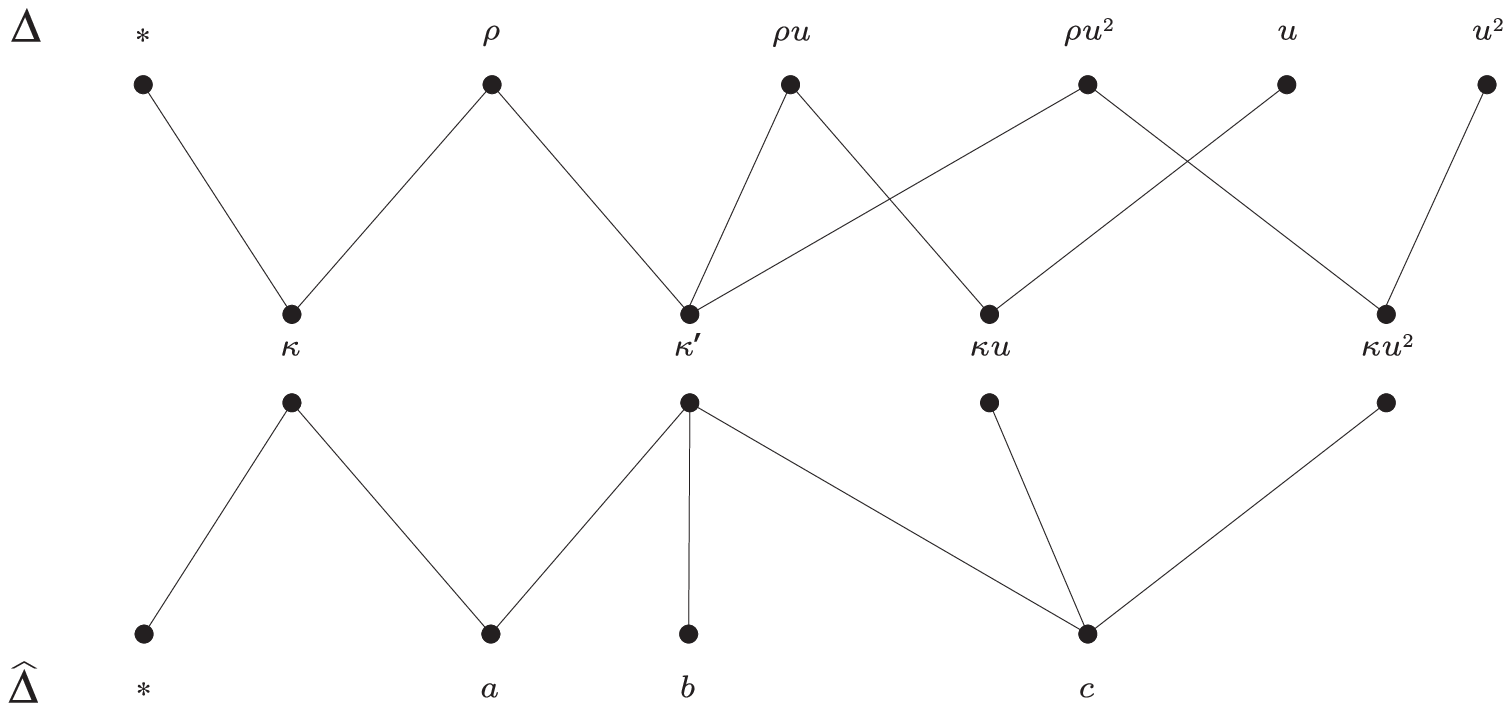}}\medskip

\centerline{Figure 2. Principal graphs of the Haagerup $(5+\sqrt{13})/2$
subfactor} \label{graphshaag2}

\medskip

The $M$-$N$ sectors, $[\kappa u^i]$ and $[\kappa']$ where $\kappa$ is the
inclusion $N\subset M$, are the odd vertices of both graphs in Figure 2.
From this figure we see we can choose 
the endomorphism $c={\kappa}u\overline{\kappa}\cong{\kappa}{u^2}
\overline{\kappa}$. These principal graphs encode multiplication by $\kappa$
--- e.g. $[\kappa][\rho u^i]=[\kappa']+[\kappa u^{i}]$. From this we
obtain the statistical dimensions dim$(\kappa u^i)=\lambda$ and dim$(\kappa')
=\lambda\,(\delta-1)$. The remaining products $M$-$N\times N$-$N\rightarrow
M$-$N$ and $M$-$M\times M$-$N\rightarrow M$-$N$, as well as 
$M$-$N\times N$-$M\rightarrow M$-$M$ and $N$-$M\times M$-$N\rightarrow N$-$N$, were
computed by Bisch \cite{Bisch} and are generalised in Section 4.2 below.

The fusions of $\Delta$ and $\widehat{\Delta}$ should remind one  of $S_3$. In particular, the system $\Delta$ has
$\bbZ_{3}=\{\hbox{id},u,{u^2}\}$ as a subsystem and is some sort of
perturbation of the usual $S_{3}$ multiplication table $[u]^{3}
= [\hbox{id}],\quad [u][\rho] = [\rho][u]^{2},\quad
[\rho]^{2} = [\hbox{id}]$. Similarly, the $\widehat{\Delta}$
system reduces to the character ring $\widehat{S_{3}}$ when
ignoring $[c]$,  where $b$ should be regarded as $S_3$-representation
$\epsilon$  and $a$ as $\sigma$. These similarities were our first indication of a relation between
the Haagerup subfactor and $S_3$. We find in Section 4 that this $\Delta\leftrightarrow
S_3$ relation generalises naturally, whereas $\widehat{\Delta}\leftrightarrow \widehat{S}_3$
does not.

The structure of the tube algebra of the system $\Delta$ has been
studied in \cite[Sect.\ 8]{iz3}. Izumi introduced
the following endomorphisms in $\Sigma(\Delta)$:
\begin{equation} \mu=\rho+ \rho u+\rho u^2\,,\quad \pi_1={\rm id}+\mu\,,\quad
\pi_2=2({\rm id})+\mu\,,\quad \sigma=u+{u^2}+\mu\,, \nonumber\end{equation}
 and proved that $\pi_1$ and $\pi_2$ have only one half-braiding, $\sigma$ has
three ($\sigma^{(1)},\sigma^{(2)},\sigma^{(3)})$, and finally $\mu$ has six  
($\mu^{(1)},\ldots,\mu^{(6)}$) and that these half-braidings exhaust all of them, so that the quantum double 
of $\Delta$ has 12 primaries as in the middle row of Figure 3. The induction-restriction graph of the Longo-Rehren inclusion of the
system $\Delta$ is given in \cite[Fig.\ 5]{iz3}, see also the top of Figure 3. From this we read off
the Longo-Rehren dual sector associated to
$\Delta$, namely $[\theta_{\cD\Delta}] = [0]+[\mathfrak{b}]+ 2[\mathfrak{a}]$,
using the labelling of Figure 3.
 This corresponds to the canonical modular invariant (recall Section 1.3)
$Z_{22} = |\chi_0 + \chi_{\mathfrak{b}} + 2 \chi_{\mathfrak{a}}|^2$, 
with corresponding full system $\Delta \times \Delta^{\mathrm{opp}}$, where the name $Z_{22}$
comes from the list of   $\cite{EP3}$ (see also  Section 3.4  below).  From 
Figure 3 we obtain
the quantum dimensions $1+3\delta,2+3\delta,2+3\delta,3\delta$ for $\mathfrak{b}, \mathfrak{a}, 
\mathfrak{c}_{j}, \mathfrak{d}_l$, respectively.

\medskip\epsfysize=3.5in\centerline{ \epsffile{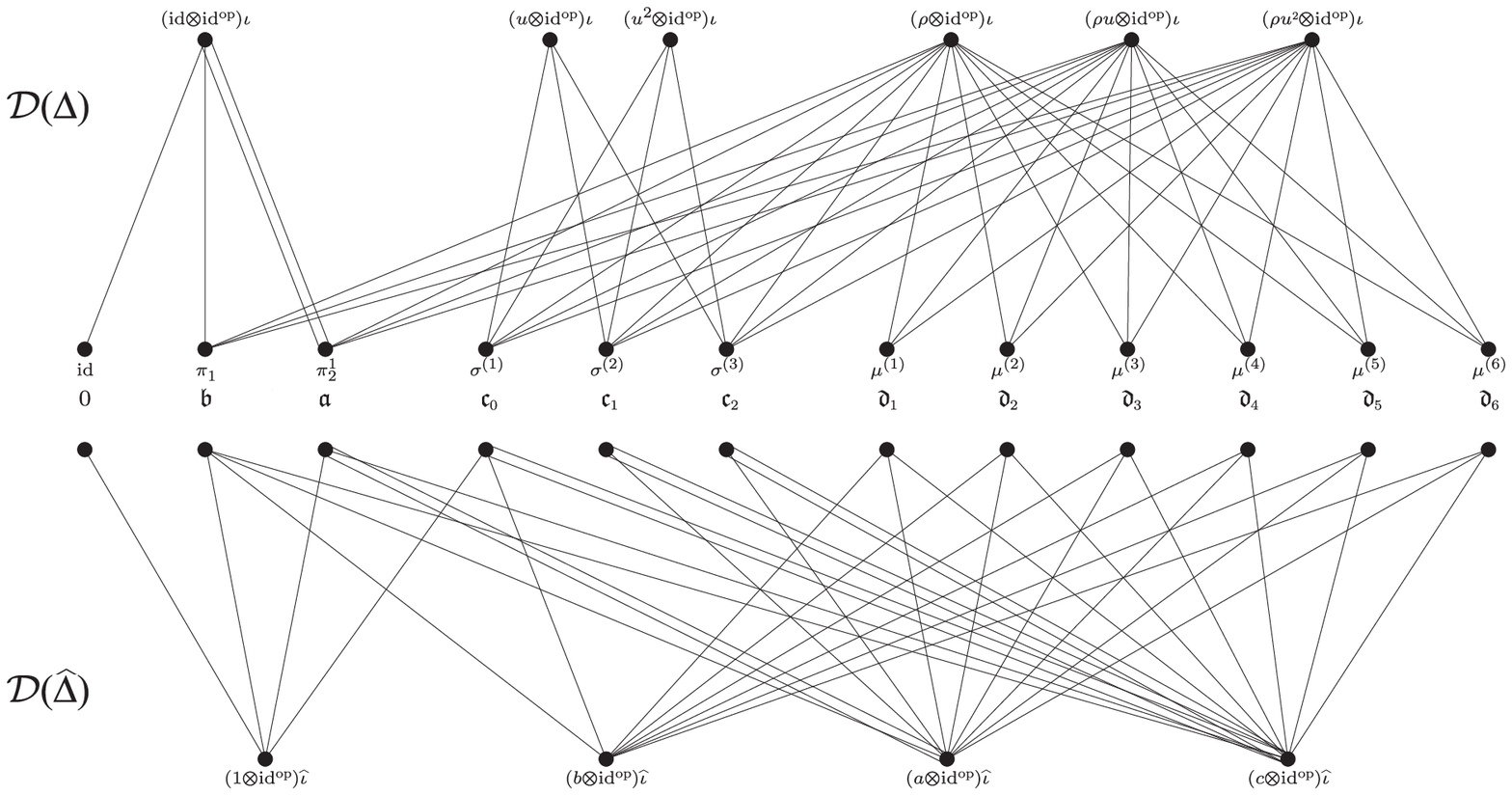}}\medskip

\centerline{Figure 3. Dual principal graphs for doubles of $\Delta$ and $\widehat{\Delta}$
for Haagerup}
\medskip

More subtle is the dual principal graph corresponding to  $\widehat{\Delta}$, on the bottom of
 Figure 3, with dual canonical endomorphism
$ [\theta_{\cD\widehat{\Delta}}]=  [0]+[\mathfrak{b}]+[\mathfrak{a}]
+[\mathfrak{c}_0]$ and modular invariant $Z_{11} = |\chi_0 + 
\chi_{\mathfrak{b}} + \chi_{\mathfrak{a}} +\chi_{\mathfrak{c}_0}|^2$ corresponding
to $\widehat\Delta \times \widehat\Delta^{\mathrm{opp}}$. 
 To our knowledge this appears here for the first time. To derive it, first note that the
 only monomial modular invariants  for $\cD\Hg$ are $Z_{22}=|ch_0+ch_{\mathfrak{b}}+2ch_{\mathfrak{a}}|^2$, $Z_{11}=|ch_0+ch_{\mathfrak{b}}+ch_{\mathfrak{a}}+ch_{\mathfrak{c}_0}|^2$,
 and $Z_{33}=|ch_0+ch_{\mathfrak{b}}+2ch_{\mathfrak{c}_0}|^2$ 
(Proposition 4 below with $\nu=3$). The number of sectors in $\widehat{\Delta}$
 is 4, and the canonical endomorphism $\gamma_{\widehat{\Delta}}=0\otimes 0^{\mathrm{opp}}
 +a\otimes a^{\mathrm{opp}}+b\otimes b^{\mathrm{opp}}+c\otimes c^{\mathrm{opp}}$.
 Since $\langle \gamma_{\widehat{\Delta}},\gamma_{\widehat{\Delta}}\rangle=4$,
 we have $\langle\widehat{\theta},\widehat{\theta}\rangle=4$ by Frobenius reciprocity.
 Thus the canonical modular invariant $Z$ associated to $\widehat{\theta}$ satisfies
 $\sum_i Z_{0,i}^2=4$, forcing it to be $Z_{11}$ and
 fixing $\widehat{\theta}$. Hence $\widehat{\Delta}$ has $1^2+1^2+1^2+1^2=4$ 
 irreducible sectors, call them $e_0,e_1,e_2,e_3$, one of which (say $e_0$) is the identity 
 $\widehat{\alpha}_{0}$. To obtain the dual principal graph, in which $e_i$ is connected
 to the primary $x$ with $\langle \widehat{\alpha}_{e_i},x\rangle$ edges, 
 use $\langle \widehat{\alpha}_x,\widehat{\alpha}_y
 \rangle=\langle\widehat{\theta}x,y\rangle$, which holds for any type I $\widehat{\theta}$
 (see Theorem 1). The fusions for the double Haagerup are explicitly given in \cite{EP3}
(or Section 3.2 below). The 
 edges from $\mathfrak{a}$  are determined from the calculations
  $\langle \widehat{\alpha}_{\mathfrak{a}},\widehat{\alpha}_{\mathfrak{a}}\rangle=6=1^2+1^2+2^2$
 and $\langle\widehat{\alpha}_{\mathfrak{a}},\widehat{\alpha}_0\rangle=1$: 
call $e_1$ the sector not adjacent to $\mathfrak{a}$, and $e_2$ the one
connected to $\mathfrak{a}$ with 2 edges. The additivity of the statistical 
dimension
$\mathrm{dim}\,\widehat{\alpha}_{\mathfrak{a}}=\mathrm{dim}\,\mathfrak{a}$ then identifies
$e_1=b,e_2=a,e_3=c$.
 Likewise, the edges from $\mathfrak{b}$ come from $\langle\widehat{\alpha}_{\mathfrak{b}},
 \widehat{\alpha}_{\mathfrak{b}}\rangle=4$ and $\langle\widehat{\alpha}_{\mathfrak{b}},
 \widehat{\alpha}_{0}\rangle=1$, and those from $\mathfrak{d}_l$ come from
 $\langle \widehat{\alpha}_{\mathfrak{d}_l},\widehat{\alpha}_{\mathfrak{d}_{l}}\rangle=3$
 and $\langle\widehat{\alpha}_{\mathfrak{d}_l},\widehat{\alpha}_0\rangle=0$.
 From $\langle\widehat{\alpha}_{\mathfrak{c}_j},\widehat{\alpha}_{\mathfrak{c}_{j}}\rangle=5+
 \delta_{j,0}$,  $\langle\widehat{\alpha}_{\mathfrak{c}_j},\widehat{\alpha}_{0}\rangle=\delta_{j,0}$, 
 and statistical dimensions, we obtain the final edges.  
  
 We have $A\subset N\otimes N^{\mathrm{opp}}$, the double of $\Delta$, with subsystem
 $\bbZ_3\subset\Delta$. Hence by the Galois theory of Izumi \cite{iz1}, we have as in
 Section 1.3 and last subsection an intermediate subfactor $A\subset C\subset N\otimes N^{\mathrm{opp}}$
 where the dual canonical endomorphism of $C\subset N\otimes N^{\mathrm{opp}}$ is
 $\gamma=\sum_{\alpha\in\bbZ_3}\alpha\otimes\alpha^{\mathrm{opp}}$ and the
 double $C$-$C$ system is the quantum double $\cD\bbZ_3\cong K^0_{\bbZ_3}(\bbZ_3)$.
 Then $A\subset C$ is a braided subfactor of index $2+3\delta$, with canonical
 endomorphism $[\theta]=[0]+[\mathfrak{b}]$ and associated modular invariant
(recall \eqref{Z3S3})
\begin{equation} Z=|ch_0+ch_\mathfrak{b}|^2+2|ch_\mathfrak{a}|^2+2|ch_{\mathfrak{c}_0}|^2+
2|ch_{\mathfrak{c}_1}|^2+2|ch_{\mathfrak{c}_2}|^2\,.\label{Z3Hg}\end{equation}

\subsection{Haagerup modular data $\cD \Hg$}

The modular data $\cD \Hg$ for  the quantum double of the even part of the Haagerup subfactor 
was first computed in \cite{iz3} and simplified somewhat in \cite{EP3}.
It is necessary for all that follows though to significantly
simplify it further. The result is
\begin{subequations}\begin{align}\label{Thaag}
T&={\rm diag}(1,1,1,1,\xi_3,\overline{\xi_3},\xi_{13}^6,
\xi_{13}^{-2} ,\xi_{13}^2,\xi_{13}^{5},\xi_{13}^{-6},\xi_{13}^{-5})\,,&\\ \label{Shaag}
S&=\frac{1}{3}{\scriptsize \left(\begin{array}{cccccccccccc}
x&1-x&1&1&1&1&y&y&y&y&y&y\cr
1-x&x&1&1&1&1&-y&-y&-y&-y&-y&-y\cr  1&1&2&-1&-1&-1& 0&0&0&0&0&0\cr 
1&1&-1&2&-1&-1& 0&0&0&0&0&0\cr 1&1&-1&-1&-1&2& 0&0&0&0&0&0\cr 
1&1&-1&-1&2&-1& 0&0&0&0&0&0\cr y&-y& 0&0&0&0&c(1)&c(2)& c(3)& c(4)& c(5)& c(6)\cr
y&-y& 0&0&0&0&c(2)&c(4)& c(6)& c(5)& c(3)& c(1)\cr y&-y& 0&0&0&0&c(3)&c(6)& c(4)& c(1)&
 c(2)& c(5)\cr y&-y& 0&0&0&0&c(4)&c(5)& c(1)& c(3)& c(6)& c(2)\cr y&-y& 0&0&0&0&c(5)&c(3)&
  c(2)& c(6)& c(1)& c(4)\cr y&-y& 0&0&0&0&c(6)&c(1)& c(5)& c(2)& c(4)& c(3)
\end{array}\right)}\,,&\end{align}\end{subequations}
for $x=(13-3\sqrt{13})/26$, $y=3/\sqrt{13}$ and $c(j)=-2y\,\cos(2\pi j/13)$. 

It is important to note that the last 6 diagonal entries of $T$ are $T_{\mathfrak{d}_l,\mathfrak{d}_{l}}
=\xi_{13}^{6l^2}$  for $1\le l\le 6$, and the bottom-right $6\times 6$ submatrix has entries $S_{\mathfrak{d}_l,
\mathfrak{d}_{l'}}=c(ll')/3$ for $1\le l,l'\le 6$.
Our expression for that $6\times 6$  submatrix is considerably simpler than the corresponding expressions in
\cite{iz3,EP3,HRW}. The proof of the equivalence is easy, e.g. use \eqref{SNT}. 
A direct derivation is given in Theorem 5 below.

$\cD \Hg$ should be compared to the modular data $\cD S_3$ for the (untwisted) double of $S_3$, which
can be obtained from (\ref{moddataG}) (or see (\ref{STDk}) below). As would be anticipated from 
the previous subsection, they are very similar, except that $x$ becomes $1/2$ and
the last 6 rows/columns collapse into 2.

The given matrix $T$ forces central charge $c$ to be a multiple of 24 (since $T_{0,0}=1$), 
but multiplying $T$ by a  third root of unity (and leaving $S$ unchanged) allows
us to consider the Haagerup at any multiple $c$ of 8.
Next section we generalise $\cD \Hg$ in two ways: it can be twisted by a $\bbZ_3$
(this twist is analogous to the $H^3(B_G;\bbT)$-twist of finite group modular data, or the level $k\in
\bbZ_{>0}$ of affine algebra modular data); and it
lies in an infinite sequence corresponding to the odd dihedrals.

\subsection{Clarifying the Haagerup-$Sym(3)$ relation}

To emphasise that the relations between $S_3$ and the Haagerup aren't spurious,
let us consider possible character vectors realising these modular data
(recall Section 1.4). This will also
help identify a  VOA realisation. In this subsection we focus on the simplest possibility,
central charge $c=8$ (the smallest possible), although similar results hold for any other
multiple of 8.

Consider first the modular data $\cD S_3$ for the quantum double of $S_3$, with primaries
$0,\mathfrak{b},\mathfrak{a},\mathfrak{c}_i,\mathfrak{d}_l$ as in Figure 1.  This
8-dimensional  PSL$_2(\bbZ)$-representation decomposes into
3 copies of the 1-dimensional irrep with $T=e^{-2\pi\i/3}$, a 2-dimensional irrep with
kernel the principal congruence subgroup $\Gamma(2)$, and a 3-dimensional irrep with kernel 
containing $\Gamma(3)$. These 1-, 2-, and
3-dimensional irreps are also subrepresentations of the modular data for the lattices $E_8$,
$D_8$ and $A_2\oplus E_6$ respectively, and as explained at the end of Section 1.4 the
corresponding theta functions provide all the information needed to extract $\Lambda$ and
$\chi$ for these irreps. From this we obtain for $\cD S_3$ at $c\equiv_{24} 8$ the
matrices  $\Lambda_{S3}^{\equiv 8}={\rm diag}(2/3,2/3,2/3,2/3,0,1/3,2/3,1/6)$
and 
$$\chi_{S3}^{\equiv 8}={\scriptsize \left(\begin{matrix}39& 47& 81& 81& 8748& 
1215& 128& 5120\cr  47& 39& 81& 81& 8748& 1215& -128& -5120\cr
  81& 81& 167& -81& -8748& -1215& 0& 0\cr 81& 81& 
-81& 167& -8748& -1215& 0& 0\cr  3& 3& -3& -3& -12& 18& 0& 0\cr  27& 27& 
-27&-27& 1458& -152& 0& 0\cr  128& -128& 0& 0& 0& 0& 120& -5120\cr  16& 
-16& 0& 0& 0& 0& -16& 140\end{matrix}\right)}\,$$
(throughout the paper we write $a\equiv_n b$ for $a\equiv b$ (mod $n$)).

Once we know $\Lambda,\chi$, then the simple recursion \eqref{recurs} yields
the full $q$-expansion of $\Xi(\tau)$.
From $\Xi$ we know completely explicitly all possible weakly holomorphic
vector-valued modular functions 
for this SL$_2(\bbZ)$-representation. Specialising to $c=8$,  we find there is a unique 
possible character vector, namely the first column of $\Xi^{\equiv 8}_{S3}$:
\begin{equation}\label{chS3+08}
{\scriptsize\left(\begin{matrix} ch_0(\tau)\cr ch_{\mathfrak{b}}(\tau)\cr ch_{\mathfrak{a}}
(\tau)=ch_{\mathfrak{c}_0}(\tau)\cr ch_{\mathfrak{c}_1}(\tau)\cr
ch_{\mathfrak{c}_2}(\tau)\cr ch_{\mathfrak{d}_1}(\tau)\cr ch_{\mathfrak{d}_2}(\tau)\end{matrix}\right)}={\scriptsize \left(\begin{matrix}
q^{-1/3}(1+39q+699q^2+5761q^3+35593q^4+\cdots)\cr 
q^{2/3}(47+671q+5825q^2+35459q^3+\cdots)\cr
81q^{2/3}(1+17q+143q^2+877q^3+\cdots)\cr 3+243q+2916q^2+21870q^3+\cdots\cr 
27q^{1/3}(1+22q+221q^2+1476q^3+\cdots)\cr 128q^{2/3}(1+16q+136q^2+832q^3+\cdots)\cr 
16q^{1/6}(1+36q+394q^2+2776q^3+\cdots)
\end{matrix}\right)}\,.
\end{equation}

It is easy to realise $\cD S_3$ (and hence the character vector 
(\ref{chS3+08})) at $c=8$. Start with the lattice VOA $\mathcal{V}(E_8)$ corresponding to
the $E_8$ root lattice. Its automorphism group is the compact Lie group 
$E_8(\bbR)$. Then  the VOA realising (\ref{chS3+08}) is the orbifold 
of $\mathcal{V}(E_8)$ by some subgroup $G$ of $E_8(\bbR)$ isomorphic to $S_3$ 
(i.e. is the subVOA of $\mathcal{V}(E_8)$ fixed by $G$). This is indeed possible (see Theorem
4.3 of \cite{KacT}).

Now turn to the double of the Haagerup, also at $c=8$, with primaries labelled as in Figure 3.
Its 12-dimensional PSL$_2(\bbZ)$-representation  decomposes into $1+1+3+7$, 
where the 1- and 3-dimensional irreps are as before, and the 7-dimensional irrep occurs 
as a subrepresentation of the modular data of the 4-dimensional lattice $A_352_1[1,1/4]$,
 in the gluing notation of Conway-Sloane \cite{CS}. From this we quickly obtain 
 for $c\equiv_{24} 8$ its matrices $\Lambda$ and $\chi$:
\begin{align}&\,{\rm
diag}(2/3,2/3,2/3,2/3,0,1/3,5/39,20/39,32/39,2/39,8/39,11/39)\,,\nonumber\\
&\,{\scriptsize \left(\begin{array}{cccccccccccc}6& 80& 81& 81& 8748& 1215
& 3549& 273& 13& 5538& 2275& 1378\cr 
80& 6& 81& 81& 8748& 1215& -3549& -273& -13& -5538& -2275& -1378\cr 
81& 81& 167& -81& -8748& -1215& 0& 0& 0& 0& 0& 0\cr 
81& 81& -81& 167& -8748& -1215& 0& 0& 0& 0& 0& 0\cr 
3& 3& -3& -3& -12& 18& 0& 0& 0& 0& 0& 0\cr 
27& 27& -27& -27& 1458& -152& 0& 0& 0& 0& 0& 0\cr
7& -7&0& 0& 0& 0& -88& -14& -1& 50& 63& 64\cr
42& -42&0&0&0& 0&-1484&92&16&2940& -192&-1041\cr
119&-119&0& 0&0&0&-2142&987&11&-24990&-6035&4641\cr
5& -5& 0& 0& 0& 0&17&13&-3&-2& 35&-14\cr
13&-13&0& 0&0&0&174&-1&-5&294&-147& 51\cr
14&-14&0&0& 0& 0&448&-77& 7&-343& 125& -24\end{array}\right)}\end{align}
%7& -7& 0& 0& 0& 0& -1& -14& 50& 64& -88& 63\cr 
%42& -42& 0& 0& 0& 0& -1484& 92& 16& 2940& -192& -1041\cr 
%119& -119& 0& 0& 0& 0& -2142& 987& 11& -24990& -6035& 4641\cr 
%5& -5& 0& 0& 0& 0& 17& 13& -3& -2& 35& -14\cr 
%13& -13& 0& 0& 0& 0& 174& -1& -5& 294& -147& 51\cr 
%14& -14& 0& 0& 0& 0& 448& -77& 7& -343& 125& -24
As before, this gives us the full $q$-expansion of $\Xi^{\equiv 8}_{hg}$.

There are only two possible character vectors for the Haagerup modular
data at $c=8$,  namely $\gamma=0$ or $\gamma=1$ in
$${\scriptsize\left(\begin{matrix} ch_0(\tau)\cr ch_{\mathfrak{b}}(\tau)\cr ch_{\mathfrak{a}}(\tau)=ch_{\mathfrak{c}_0}(\tau)\cr ch_{\mathfrak{c}_1}(\tau)\cr
ch_{\mathfrak{c}_2}(\tau)\cr ch_{\mathfrak{d}_1}(\tau)\cr ch_{\mathfrak{d}_2}(\tau)\cr ch_{\mathfrak{d}_3}(\tau)\cr ch_{\mathfrak{d}_4}(\tau)\cr ch_{\mathfrak{d}_5}(\tau)\cr ch_{\mathfrak{d}_6}(\tau)
\end{matrix}\right)}={\scriptsize\left(\begin{matrix}q^{2/3}\left(q^{-1}+(6+13\gamma)
+(120+78\gamma)q+(956+351\gamma)q^2+(6010+1235\gamma)q^3+\cdots\right)\cr 
q^{2/3}\left((80-13\gamma)+(1250-78\gamma)q+(10630-351\gamma)q^2+ 
(65042-1235\gamma)q^3+\cdots\right)\cr 
q^{2/3}\left(81+1377q+11583q^2+71037q^3+\cdots\right)\cr 
 3+243q+2916q^2+21870q^3+\cdots\cr 
q^{1/3}\left(27+594q+5967q^2+39852q^3+\cdots\right)\cr q^{5/39}\left((7-\gamma)
+(292-6\gamma)q+(3204-43\gamma)q^2+(23010-146\gamma)q^3+\cdots\right)\cr
q^{20/39}\left((42+16\gamma)+(777+121\gamma)q+(7147+547\gamma)q^2+
(45367+2000\gamma)q^3+\cdots\right)\cr 
q^{32/39}\left(\gamma q^{-1}+(11\gamma+119)+(73\gamma+1623)q+(300\gamma+12996)q^2+
(76429+1063\gamma)q^3+\cdots\right)\cr q^{2/39}\left((5-3\gamma)+(229-50\gamma)q
+(2738-252\gamma)q^2+(19942-1032\gamma)q^3+\cdots\right)\cr
 q^{8/39}\left((13-5\gamma)+(347-37\gamma)q+(3804-212\gamma)q^2+
(26390-794\gamma)q^3+\cdots\right)\cr
 q^{11/39}\left((14+7\gamma)+(441+61\gamma)q+(4445+303\gamma)q^2+
(30329+1167\gamma)q^3+\cdots\right)
\end{matrix}\right)}$$
The proof for this (and for \eqref{chS3+08}) is similar to that sketched in Section 5.2.
That there are sensible character vectors here is strong evidence for the
existence of the corresponding VOA.

We see remarkable similarities between $\cD \Hg$ and 
 $\cD S_3$  when $c=8$. In particular, the corresponding characters 
 $ch_{\mathfrak{a}},ch_{\mathfrak{c}_i}$ are 
 \textit{identical}, as are the sums $ch_0+ch_{\mathfrak{b}}$.
The remaining characters $ch_{\mathfrak{d}_1},\ldots,ch_{\mathfrak{d}_6}$ of $\cD \Hg$ can have
no relation  with characters $ch_{\mathfrak{d}_1},ch_{\mathfrak{d}_2}$ of  $\cD S_3$, 
since the exponents of their $q$-expansions are unrelated.

This suggests that there is a rational VOA 
at $c=8$, call it $\mathcal{V}_8$, which has conformal subalgebras realising  the  $\cD \Hg$
and $\cD S_3$ modular data (\textit{conformal subalgebra} was defined at the beginning of Section 1).
The  $\mathcal{V}_8$ characters will be $ch_0+ch_{\mathfrak{b}}$ (the vacuum character) and
some multiples of $ch_{\mathfrak{a}},ch_{\mathfrak{c}_i}$, for the $ch$'s given in (\ref{chS3+08}). 

Indeed, the knowledge of modular invariants says how this should go. We should look for  modular invariants
of extension type  (i.e. sum of squares)  for both $\cD S_3$ and $\cD \Hg$, which start
with $|ch_0+ch_{\mathfrak{b}}|^2$. Both $\cD S_3$ and $\cD\Hg$  have one, namely
\eqref{Z3S3} and \eqref{Z3Hg}. 
 In the case of $S_3$ we readily identify the corresponding chiral extension: $\mathcal{V}_8$ 
corresponds to the quantum double $\cD \bbZ_3$ of the cyclic group $\bbZ_3$. The modular
invariant says two inequivalent $\mathcal{D}
\bbZ_3$-primaries correspond to each $ch_{\mathfrak{a}}$, $ch_{\mathfrak{c}_i}$. 

So $\mathcal{V}_8$ is the $\bbZ_3$-orbifold (i.e. the $\bbZ_3$-invariant part)
 of the rational
VOA $\mathcal{V}(E_{8})$, for some choice of order-3 element $g$ of 
$E_8(\bbR)$.  The group $E_8(\bbR)$ contains
four inequivalent order 3 elements: orbifolding by any of these would
give a VOA with 9 primaries realising modular data associated to a quantum
double of $\bbZ_3$. Half the time this double is twisted and half the time it
isn't. Take one of the order 3 orbifolds corresponding to the untwisted double
$\cD\bbZ_3$. To recover the $\mathcal{D}S_3$ VOA,
choose an order-2 element $h$ of $E_8(\bbR)$ for which $gh=hg^{-1}$: the $h$-orbifold of
$\mathcal{V}_8$ realises $\mathcal{D}S_3$. Unfortunately the Haagerup VOA won't itself be 
an orbifold of $\mathcal{V}_8$ (see Section 5.1 below), which
 begs Question 3 of Section 1.

To better identify $\mathcal{V}_8$,
we should repeat the analysis for  $\cD \bbZ_3$.
Label its 9 primaries by $(i,j)\in\bbZ_3^2$, where $(0,0)$ is the vacuum.
Note that charge-conjugation $C=S^2$ is no longer the identity: it sends $(i,j)$
to $(-i,-j)$. Hence the characters $ch_{(i,j)}$ and $ch_{(-i,-j)}$ are equal
and we should project to the subrepresentation Span$\{ch_{(0,0)},
(ch_{(0,1)}+ch_{(0,2)})/2,(ch_{(1,0)}+ch_{(2,0)})/2,(ch_{(1,1)}+ch_{(2,2)})/2,
(ch_{(1,2)}+ch_{(2,1)})/2\}$ on which
$\bbZ_2\cong\langle C\rangle$ acts trivially. This decomposes into
$1+1+3$, where the 1- and 3-dimensional irreps are as before. Choosing the
above basis,
we obtain $\Lambda_{Z3}^{\equiv 8}=\mathrm{diag}(2/3,2/3,2/3,0,1/3)$
and $$\chi_{Z3}^{0,\equiv 8}={\scriptsize\left(\begin{matrix} 86 &   162 &  162 &  17496 &  2430\cr
  81 &   167 &   -81 &   -8748  &  -1215\cr 81 &   -81 &   167 &   -8748 &   -1215\cr
  3  &   -3  &   -3  &    -12   &    18\cr 27 &   -27  &  -27  &   1458  &   -152\end{matrix}\right)}\,.$$
It is immediate that the only character vector possible for $\mathcal{V}_8$ is
\begin{equation}
{\scriptsize\left(\begin{matrix}ch_{(0,0)}\cr ch_{(0,1)}=ch_{(0,2)}=ch_{(1,0)}=ch_{(2,0)}\cr 
ch_{(1,1)}=ch_{(2,2)}\cr ch_{(1,2)}=ch_{(2,1)}\end{matrix}
\right)}={\scriptsize\left(\begin{matrix} q^{-1/3}(1+86q+1370q^2+11586q^3+\cdots)\cr
q^{2/3}(81+1377 q + 11583 q^2 +71037q^3+\cdots )  \cr   
   3 + 243 q + 2916 q^2+21870q^3+\cdots  \cr  q^{1/3}(27 + 594 q + 5967 q^2+39852q^3+\cdots)\end{matrix}\right)}\,.\end{equation}    
We recognise this as the character vector for the lattice VOA $\mathcal{V}(A_2 E_6)$,
so this is the VOA $\mathcal{V}_8$ containing both the $S_3$ and Haagerup VOAs. This $
\mathcal{V}_8$ also has, like $\mathcal{V}(E_8)$,  an interpretation as an affine algebra VOA, 
and the containment $\mathcal{V}_8\subset \mathcal{V}(E_8)$ corresponds to the
conformal embedding $A_{2,1}E_{6,1}\subset E_{8,1}$. It can also be realised explicitly
as an orbifold --- see \cite{KacT}. The only
task remaining is {to}  identify the Haagerup VOA as a subalgebra of  $\mathcal{V}(A_2 E_6)$.

\subsection{The Haagerup and SO(13)}

The relationship with $S_3$ concerns the primaries $\mathfrak{a},\mathfrak{c}_i$ and
to a lesser extent $0,\mathfrak{b}$. There also is a striking relationship with
the affine algebra modular data  $B_{6,2}$ concerning the primaries $\mathfrak{d}_l$ (and to a
lesser extent $0,\mathfrak{b}$).  $\cV(B_{6,2})$ has central charge $c=12$, and 
10 primaries we'll denote
$0,\mathfrak{b}=2\Lambda_1$, $\mathfrak{a}_1=\Lambda_6$, $\mathfrak{a}_2=
\Lambda_1+\Lambda_6$, $\mathfrak{d}_1=\Lambda_1,\ldots,\mathfrak{d}_5=\Lambda_5,
\mathfrak{d}_6=2\Lambda_6$. The $T$-matrix is diag$(-1,-1;-\i,
\i; -\xi_{13}^{6l^2})$, while the $S$-matrix is \cite{KaWa}
\begin{equation}
S=\frac{1}{3}{\scriptsize \left(\begin{array}{cccccccccc}
y/2&y/2&3/2&3/2&y&y&y&y&y&y\cr  y/2&y/2&-3/2&-3/2&y&y&y&y&y&y\cr  
3/2&-3/2&3/2&-3/2&0&0&0&0&0&0\cr 3/2&-3/2&-3/2&3/2&0&0&0&0&0&0\cr
y&y&0&0&-c(1)&-c(2)& -c(3)& -c(4)& -c(5)& -c(6)\cr
y&y&0&0&-c(2)&-c(4)& -c(6)& -c(5)& -c(3)& -c(1)\cr y&y& 0&0&-c(3)&-c(6)& -c(4)& -c(1)&
 -c(2)& -c(5)\cr y&y&0&0&-c(4)&-c(5)& -c(1)& -c(3)& -c(6)& -c(2)\cr y&y&0&0&-c(5)&-c(3)&
  -c(2)& -c(6)& -c(1)& -c(4)\cr y&y&0&0&-c(6)&-c(1)& -c(5)& -c(2)& -c(4)& -c(3)
\end{array}\right)}\,,\end{equation}
where $y=3/\sqrt{13}$ and $c(j)=-2y\cos(2\pi j/13)$ as before.
Ignoring the first 4 primaries, the only difference with $\cD\Hg$ are some signs.

This strongly suggests relations between $\cV(B_{6,2})$ and the (still hypothetical)
Haagerup VOA involving the Goddard--Kent--Olive coset construction \cite{GKO}.
In the VOA language, the coset construction was developed in \cite{FrZh}; see also the
lucid treatment in Section 3.11 of \cite{LeLi}). There are several ways this could go:
e.g. the Haagerup VOA at $c=8$ could be a coset of $\cV(B_{6,2})$ by a $c=4$
 subVOA $\cV_4$. In this case the characters of $\cV(B_{6,2})$ (which can be determined
by the Weyl-Kac character formula \cite[Ch.13]{kac})  would be built from the characters of
$\cV_4$ and the Haagerup. 

Although the possible relation between the Haagerup and $\cV(B_{6,2})$, in particular
 those involving the coset construction, seems very intriguing,
for reasons of space we limit the discussion in this paper to relations between the Haagerup
and $S_3$.

\section{Generalising the modular data of the Haagerup}

\subsection{Dihedral groups and orthogonal algebras}

Last section we observed that the Haagerup modular data is closely related to that of the symmetric
group $S_3=Sym(3)$ and affine so$(13)$ at level 2. More generally, we propose next subsection a two-parameter 
generalisation  of the Haagerup, related to the odd dihedral groups and
affine so$(2m+1)$ at level 2. 
 
First, we compute the modular data for quantum doubles of 
the dihedral group $D_\nu=\langle \tau,u\,|\,\tau^2,u^\nu,\tau u=u^{-1}\tau\rangle$ 
where $\nu=2n+1$. Of course $S_3\cong D_3$. The twist group is $H^3(B_{D_\nu};
\bbT)=\bbZ_{2\nu}\cong\bbZ_2\times\bbZ_\nu$; because the Schur multipliers of 
$D_\nu$ and cyclic groups all
vanish, the modular data is cohomologically trivial and is given by (\ref{moddataG}). The 
conjugacy classes of $D_{\nu}$ have representatives $e,u^h,\tau$ 
for $1\le h\le n$, with centralisers $D_{\nu},\langle u\rangle,\langle \tau\rangle$ respectively. There 
are two 1-dimensional
irreps of $D_{\nu}$ (call them $\epsilon_0=1,\epsilon_1$) and $n$ 2-dimensional ones (call them
$\sigma_i$ for $1\le i\le n$); denote the $\nu$ 1-dimensional irreps of $\langle u\rangle\cong\bbZ_\nu$ by
$\pi_j$ for $0\le j< \nu$, and the two 1-dimensional irreps of $\langle \tau\rangle\cong\bbZ_2$
by $\epsilon_0,\epsilon_1$ again. The primaries fall into four classes:\smallskip

\noindent\textit{1.} two primaries: the vacuum $0:=(e,1)$ and $\mathfrak{b}:=(e,\epsilon_1)$;

\noindent\textit{2.}  $n$ primaries, labeled $\mathfrak{a}_i:=(e,\sigma_i)$ for $1\le i\le n$;

\noindent\textit{3.}  $n\nu$ primaries, labeled $\mathfrak{c}_{h,j}:=(u^h,\pi_j)$ for $1\le h\le n, 0\le j< \nu$;

\noindent\textit{4.} two primaries, labeled $\mathfrak{d}_l:=(\tau,\epsilon_l)$, $l=1,2$.

\smallskip Fix a sign $s=\pm $ and integer $\omega\in\bbZ$. Corresponding to these 4 classes we have the modular data
\begin{align}\label{STDk}
T=&\,{\rm diag}(1,1;1,\ldots,1;\exp[2\pi\i\,(\omega h^2+\nu hj)/\nu^2];t,-t)\,,\cr
S=&\,\frac{1}{\nu}\left(\begin{matrix}
\frac{1}{2}_{2\times 2}&1_{2\times n}&1_{2\times n\nu}&\nu B\cr 1_{n\times 2}&2_{n\times n}&D&0_{n\times 2}\cr 1_{n\nu\times 2}&D^t&E&
0_{n\nu\times 2}\cr \nu B^t&0_{2\times n}&0_{2\times n\nu}&s\nu F\end{matrix}\right)\,,\end{align}
where $t=1,\i$ for $s=1,-1$ resp., $k_{a\times b}$ for any number $k$ is the $a\times b$ matrix with constant entry $k$, 
$$B=\frac{1}{2}\left(\begin{matrix}1&1\cr -1&-1\end{matrix}\right)\ {\rm and}\ F=\frac{1}{2}\left(\begin{matrix}1&-1\cr -1&1\end{matrix}\right)\,,$$
$D_{i,(h,j)}=2\cos(2\pi ih/\nu)$, and
$E_{(h,j),(h',j')}=2\cos(2\pi (2\omega hh'+\nu hj'+\nu h'j)/\nu^2)$. We'll denote this modular data by
$\cD^{s,\omega}D_\nu$. The untwisted double is $\cD^{+,0}D_\nu$. 
Note that  $\cD^{s,\omega+k\nu}D_\nu$ is equivalent for any integer $k$, so the
twist does indeed live in $\bbZ_{2\nu}$. 
We've separated the order 2 from the order $\nu$ twists, for later convenience.
The quantum-dimensions ($=$statistical dimensions) $S_{*0}/S_{00}$ of $\mathfrak{b},\mathfrak{a}_i,
\mathfrak{c}_{h,j},\mathfrak{d}_l$ are 1, 2, 2, $\nu$, respectively, and the global dimension
$1/S_{0,0}$ is $2\nu$.  

The $S,T$ entries live in $\bbQ[\xi_{4\nu^2}]$. Choose any Galois
automorphism $\sigma\in \mathrm{Gal}(\bbQ[\xi_{4\nu^2}]/\bbQ)$ fixing $\i$ and define $\ell\in
\bbZ_{\nu^2}$ by $\sigma\xi_{\nu^2}=\xi_{\nu^2}^\ell$. Then applying $\sigma$ to
$\cD^{s,\omega}D_\nu$, as explained in Section 1.1, yields modular data equivalent to $\cD^{s,\ell \omega}
D_\nu$. In particular the fusion ring of $\cD^{s,\omega}D_\nu$ depends only on $\nu$, $s$
and gcd$(\omega,\nu)$. In fact, if $\ell\in\bbZ_{\nu}^\times$, then
$\cD^{s,\ell^2 \omega}D_\nu$ and $\cD^{s,\omega}D_\nu$ are equivalent.

The second ingredient going into our generalisation of the Haagerup modular data is the
 affine algebra modular data 
$B_{m,2}$, which has central charge $c=2m$. In particular \cite{KaWa}, $B_{m,2}$
 has $m+4$ primaries,
labelled $0,\mathfrak{b}=2\Lambda_1$, $\mathfrak{a}_1=\Lambda_m$, $\mathfrak{a}_2=
\Lambda_1+\Lambda_m$, $\mathfrak{d}_l=\Lambda_l$ for 
$1\le l<m$ and $\mathfrak{d}_m=2\Lambda_m$. Write $\mu=2m+1$. The $T$-matrix is 
$\xi_{12}^{-m}\mathrm{diag}(1,1;\xi_8^m,-\xi_8^m; \xi_{\mu}^{ml^2})$, while the $S$-matrix is
\begin{equation}S=\left(\begin{matrix}\frac{x}{2}_{2\times 2}&B&x_{2\times m}\cr  
B^T&F&0_{2\times m}\cr  x_{m\times 2}&0_{m\times 2}&H\end{matrix}\right)\,,\end{equation}
where $x=1/\sqrt{\mu}$, $B$ and $F$ are as above,  and $H_{l,l'}=2x\cos(2\pi ll'/\mu)$.
The quantum-dimensions of $\mathfrak{b},\mathfrak{a}_i,\mathfrak{d}_l$  are 
$1,\sqrt{\mu},2$, respectively, with global dimension $2\sqrt{\mu}$.

There is no obvious twist of the $B_{m,2}$.
Curiously \cite{CGR},  affine so$(\nu^2)$ at level 2 coincides with $\cD^\omega D_\nu$ for 
a specifically chosen twist $\omega$ (this is one of the few overlaps of affine algebra and finite
group modular data \cite{CGR}), but that fact seems to have no role
in our story (in particular we are not interested in the values $m=(\nu^2-1)/2$).

\subsection{Generalising the  Haagerup modular data}

Sections 2.3, 2.5 and 3.1 suggest a generalisation of  $\cD\Hg$.
Write $\nu=2n+1$  for $n\ge 0$,   and choose $\omega\in\bbZ$  
 as before. Write $m=2(n^2+n+1)$, $\mu=2m+1=\nu^2+4$, and $\delta=(\nu+\sqrt{\mu})/2$.
Note that $\delta$ satisfies $\delta^2=\nu\delta+1$ and lies in the range
$\nu<\delta<\nu+\frac{1}{2}$. Again the primaries fall into four classes:\smallskip

\noindent\textit{1.} two primaries, denoted 0 and $\mathfrak{b}$;

\noindent\textit{2.} $n$ primaries, denoted $\mathfrak{a}_i$ for  $1\le i\le n$;

\noindent\textit{3.} $n\nu$ primaries, denoted $\mathfrak{c}_{h,j}$ for $1\le h\le n, 0\le j< \nu$;

\noindent\textit{4.} $m$ primaries, denoted $\mathfrak{d}_l$ for $1\le l\le m$.

\smallskip 
Breaking $S$ and $T$ into 16 blocks, as in the previous subsection, we get
\begin{align}
T=&\,{\rm diag}(1,1;1,\ldots,1;\exp[2\pi\i(\omega h^2+\nu hi)/\nu^2];\exp[2\pi\i\, l^2m/\mu])\,,\cr
S=&\,{1\over \nu} \left(\begin{matrix}A&1_{2\times n}&1_{2\times n\nu}&B'\cr 
1_{n\times 2}&2_{n\times n}&D&0_{n\times m}\cr 1_{n\nu\times 2}&D^t&E&
0_{n\nu\times m}\cr B'{}^t&0_{m\times n}&0_{m\times n\nu}&-\nu H\end{matrix}\right)\,,\label{genhS}\end{align}
where $k_{a\times b}$, $D$, $E$ and $H$ are as in Section 3.1 (so
$-\nu H_{l,l'}=-2y\cos(2\pi ll'/\mu)$), and 
$$A=\frac{1}{2}\left(\begin{matrix}1-y&1+y\cr 1+y&1-y\end{matrix}\right)\ {\rm and}\ 
B'=y\left(\begin{matrix}1&1&\cdots&1\cr -1&-1&\cdots&-1\end{matrix}\right)$$
for $y={\nu\over\sqrt{\mu}}$. Denote this modular
data by $\cD^{\omega}\Hg_{2n+1}$; we  call $\cD^{0}
\Hg_{\nu}$ \textit{Haagerup-Izumi modular data}. To our knowledge this modular 
data is new except when $n\le 1$ and 
$\nu$ divides $\omega$; in particular twisting the Haagerup data was
unanticipated in the literature. We discuss further generalisations next 
subsection.

The quantum dimensions (or statistical dimensions) $S_{*,0}/S_{0,0}$ for
$\mathfrak{b},\mathfrak{a}_i,\mathfrak{c}_{h,j},\mathfrak{d}_l$ are respectively
$1+\nu\delta$, $2+\nu\delta$, $2+\nu\delta$, and $\nu\delta$. The global
dimension is $1/S_{0,0}=\nu\,(\nu\delta+2)$. Note that the submatrix
$-\nu H$ has exactly half of the $\mu-1$ rows and columns of $S$.

The fusions can be computed directly using Verlinde's formula \eqref{verl}
(or from Proposition 2 below). As
mentioned in Section 1.1,
it suffices to consider all $i,j,k$ different from 0. For $\cD^\omega\Hg_\nu$, 
all of those fusion coefficients are 1 except for
\begin{align} &N_{\mathfrak{b},\mathfrak{a}_i,\mathfrak{a}_i}=2\,,\ N_{\mathfrak{b},\mathfrak{c}_{h,i},\mathfrak{c}_{h,i}}=2\,,\ N_{\mathfrak{b},\mathfrak{d}_l,\mathfrak{d}_l}=0\,,\ 
N_{\mathfrak{a}_i,\mathfrak{a}_j,\mathfrak{a}_k}=2\ \mathrm{if}\ k\equiv_\nu si+s'j\,, 
&\nonumber\\
&N_{\mathfrak{c}_{h,i},\mathfrak{c}_{h',j},\mathfrak{c}_{h'',k}}=2\ \mathrm{if}\ h''\equiv_\nu sh+s'h'\, ,
\ k\equiv_\nu si+s'j+2\omega\,(sh+s'h'-h)/\nu\,,&\nonumber\\
&N_{\mathfrak{a}_i,\mathfrak{c}_{h,j},\mathfrak{c}_{h,k}}=2\ \mathrm{if}\ 
si\equiv_\nu j-k\,,\ 
N_{\mathfrak{d}_l,\mathfrak{d}_{l'},\mathfrak{d}_{l''}}=0\ \mathrm{if}\ l''\equiv_\mu sl+s'l'\,, 
&\nonumber\end{align}
where $s,s'\in\{\pm 1\}$ are arbitrary signs.

As in Section 3.1, $\cD^\omega \Hg_{2n+1}$ depends up to equivalence on the
value of the twist $\omega$ mod $\nu$; the fusion ring depends up to equivalence 
only on $n$ and gcd$(\nu,\omega)$;
and $\cD^\omega \Hg_{2n+1}\cong \cD^{\ell^2 \omega}\Hg_{2n+1}$ for any $\ell$ coprime to $\nu$.

The computation of fusions reduces to the identity 
\begin{equation}8\sum_{d=1}^k\cos\left(\frac{2\pi ad}{2k+1}\right)
\cos\left(\frac{2\pi bd}{2k+1}\right)\cos\left(\frac{2\pi cd}{2k+1}\right)
=-4+(2k+1)s\nonumber\end{equation}
proved using $2\cos(x)=e^{\i x}+e^{-\i x}$,
where $0\le s\le 4$ is the number of pairs $(s_a,s_b)$ of signs $\pm 1$ such
that $c\equiv_{2k+1}s_aa+s_bb$. Modularity reduces to the Gauss sum
$$\sum_{a=0}^{2k}\exp[2\pi\i\,ca^2/(2k+1)]=\sqrt{2k+1}\left(\frac{c}{2k+1} \right)
\left\{\begin{matrix} 1&{\rm for}\ k\ {\rm even}\cr \i&{\rm for}\ k\ {\rm odd}
\end{matrix}\right.\,.$$
Since $\mu\equiv_8 5$, $\mu$ can never be a 
perfect square, and the Jacobi symbol $\left(\frac{m}{\mu} \right)$ will always equal $-1$.
See Proposition 2 below for a more elegant argument.

Of course $\cD^0\Hg_3$ recovers the original Haagerup double $\cD\Hg$, given in Section 2.3.
The similarity with the modular data  $\cD^{+,\omega} D_{\nu}$  and $B_{m,2}$ 
is evident
 (ignoring the first 4 primaries of $B_{m,2}$, the only difference with 
 $\cD^\omega\Hg_\nu$ is that the bottom-right corner of both $T$ and
$S$ are off by some 6th root of 1). In particular, class 1 runs through the 
2-dimensional irreps of $D_\nu$, while $h$ and $j$ in class 3 parametrise the 
size-2 conjugacy classes in $D_\nu$ and 
the irreps of $\bbZ_\nu$, respectively; class 4 runs through the fundamental 
nonspinorial
weights of so$(\mu)$.  

The mysterious `13' in $\cD \Hg$ is thus $4+3^2$, where the 3 here references the normal
subgroup $\bbZ_3$  of $S_3$.
It is tempting to guess that the 4 is $2^2$, where the `2' is the involution of $S_3$. This
suggests the further generalisation of $\cD^\omega \Hg_\nu$
given in the Wildly Optimistic Guess of Section 1.

\subsection{Further generalisations}

The way in which
$\cD^\omega\Hg_\nu$ is built from $\cD^{+,\omega}D_\nu$ and 
$B_{m,2}$ leads to the notion of \textit{grafting} modular data. 
A simple instance is provided by the following proposition, 
but it can be massively generalised. For instance it would be interesting to
extend it to vacuum blocks larger than $2\times 2$.

For simplicity restrict here to \textit{unitary} modular data (recall Section
1.1) with real matrix $S$.
Call modular data $(\Phi,0,S,T)$ \textit{$\bbZ_2$-laminated} if it can be written in block form
$\Phi=(0,\mathfrak{b};\mathfrak{a}_1,\ldots,\mathfrak{a}_m;\mathfrak{d}_1,\ldots,\mathfrak{d}_n)$ 
such that $S_{0,\mathfrak{b}}=S_{0,0}$, $N_{\mathfrak{b},\mathfrak{a}_i}^{\mathfrak{a}_i}\ne 0$,
$T_{\mathfrak{b},\mathfrak{b}}=T_{0,0}$ and $S_{\mathfrak{b},\mathfrak{d}_j}<0$. 
It is elementary to verify that modular data is $\bbZ_2$-laminated iff it can be written in the
form
\begin{equation}
S=\left(\begin{matrix}x&x&\vec{a}&\vec{d}\cr x&x&\vec{a}&-\vec{d}\cr
\vec{a}{}^T&\vec{a}{}^T&A&0_{m\times n}\cr \vec{d}{}^T&-\vec{d}{}^T&0_{m\times n}&D\end{matrix}
\right)\,,\qquad T=\mathrm{diag}(r,r;\vec{s};\vec{t})\,,\end{equation}
for some numbers $x,r\in\bbR$, row vectors $\vec{a},\vec{s}\in\bbR^m$, $\vec{d},\vec{t}\in\bbR^n$ and matrices $A$, $D$.
The condition $S_{\mathfrak{b},0}=S_{0,0}$ implies the matrix $N_{\mathfrak{b}}$ with entries
$N_{\mathfrak{b},i}^j$ for all $i,j\in\Phi$ is a permutation matrix (see e.g.\ \cite{Gan2}); in conformal
field theory, such a primary $\mathfrak{b}$ is called a \textit{simple-current}.
The modular data  $\cD^{\pm,\omega}D_\nu$ and $B_{r,2}$ are both $\bbZ_2$-laminated,
for any odd $\nu$ and any twist $(\pm,\omega)$ and any rank $r$, as is
$A_{1,4}$.

\medskip\noindent\textbf{Proposition 2.} \textit{Consider $\bbZ_2$-laminated modular data
$(\Phi,0,S,T)$ and $(\Phi',0',S',T')$.} \smallskip

\noindent{(a)} \textit{There are integers $M>0$, $L\ge 0$, $N>0$ such that $2x=1/\sqrt{M}$, $N$ and $4^L$ divide $4M$,
$2^{1-L}\sqrt{M}\vec{a}\in\bbZ^m$, $2\sqrt{M/N}\vec{d}\in\bbZ^n$, and 
the gcd over all components $2^{1-L}\sqrt{M}a_i$ is 1 (likewise for the components
$2\sqrt{M/N}d_k$).}\smallskip

\noindent{(b)} \textit{Suppose $x=x'$ (that is, the global dimensions agree) and
$r=r'$ (i.e. the central charges $c,c'$ agree mod 24). Then the following
defines modular data:}
$$\widehat{S}=\left(\begin{matrix} x&x&\vec{a}&\vec{d}{}'\\
x&x&\vec{a}&-\vec{d}{}'\\ \vec{a}{}^T&\vec{a}{}^T&A&0\\ \vec{d}{}^{\prime T}&
-\vec{d}{}^{\prime T}&0&D'\end{matrix}\right)\,,\ \widehat{T}=\mathrm{diag}(r,r;
\vec{s};\vec{t}\,')\,.$$

\noindent{(c)} \textit{Suppose $x>x'$ and $r=-r'\omega$ for some 3rd root 
$\omega$ of 1 (that is, $c-c'\equiv_{8}4$). Define
$x_\pm=x\pm x'$. Then
$$\widehat{S}=\left(\begin{matrix} x_-&x_+&\vec{a}&\vec{a}{}'\\
x_+&x_-&\vec{a}&-\vec{a}{}'\\ \vec{a}{}^T&\vec{a}{}^{T}&A&0\\ \vec{a}{}^{\prime T}&
-\vec{a}{}^{\prime T}&0&-A'\end{matrix}\right)\,,\ \widehat{T}=\mathrm{diag}(r,r;
\vec{s};-\omega\vec{t}\,')\,$$
satisfies all conditions of modular data except possibly that the fusion
coefficients are nonnegative integers. Its fusion coefficients are all integers
iff $4|(M'-M)$, and both $L>0,L'>0$. 
$\widehat{S},\widehat{T}$ define modular data if in addition
$N'_{\mathfrak{a}'_k,\mathfrak{a}'_{k'},\mathfrak{a}'_{k''}}\le 8^{L'}/(M'-M)$.}

\medskip\noindent\textit{Proof.} Recall \eqref{galois}. Clearly, any $\sigma$
permutes $\{0,\mathfrak{b}\}$, so the number $x^{-2}$ will be fixed by any
$\sigma$. But $x^{-2}$ is an algebraic integer, being the sum of squares of
certain eigenvalues $S_{i,0}/S_{0,0}$ of integer matrices. Therefore $x^{-2}$
is indeed an integer. Computing 
$N_{\mathfrak{a}_i,\mathfrak{d}_k,\mathfrak{d}_k}/
N_{\mathfrak{a}_j,\mathfrak{d}_k,\mathfrak{d}_k}=a_i/a_j$ and
$N_{\mathfrak{a}_i,\mathfrak{d}_k,\mathfrak{d}_{k'}}/
N_{\mathfrak{a}_i,\mathfrak{d}_k,\mathfrak{d}_k}=d_{k'}/d_k$, we see that there are
$\beta,\gamma>0$ for which both vectors $\beta\vec{a}$ and $\gamma\vec{d}$
are integral. The gcd condition follows by choosing these 
$\beta,\gamma$ as small as  possible.
Because $\vec{d}\cdot\vec{d}=1/2$, we know $\gamma^2\in\bbQ$. Because 
$N_{\mathfrak{a}_i,\mathfrak{d}_k,\mathfrak{d}_k}=2a_id_k^2/x\in\bbZ$ and each
quantum-dimension $d_k/x$ is an algebraic integer, we know $a_i/x\in\bbZ$ and
hence $(\beta x)^{-1}\in\bbZ$. From $2x^2+\vec{a}\cdot\vec{a}=1/2$ we now see
that 4 divides $x^{-2}$, and $(\beta x)^{-1}$ is a power $2^L$ of 2. This gives us (a).

The proof of (b) and (c) is now straightforward. Most interesting are the
fusion calculations in (c). Note that $4x^2x'{}^2/(x^2-x'{}^2)=1/(M'-M)$, 
$a_i/x\in 2^{L}\bbZ$ and $a_k'/x'\in 2^{L'}\bbZ$. We obtain
\begin{align}&&\widehat{N}_{\mathfrak{b},\mathfrak{b},\mathfrak{b}}=\frac{4}{M'-M}\,,\
\widehat{N}_{\mathfrak{b},\mathfrak{b},\mathfrak{a}_i}=\frac{2a_i}{x\,(M'-M)}\,,\ 
\widehat{N}_{\mathfrak{b},\mathfrak{a}_i,\mathfrak{a}_{i'}}=\frac{a_ia_{i'}}{x^2\,(M'-M)}+\delta_{i,i'}\,,\nonumber\\
&&\widehat{N}_{\mathfrak{b},\mathfrak{b},\mathfrak{a}'_k}=\frac{2a'_k}{x'\,(M'-M)}\,,\
\widehat{N}_{\mathfrak{b},\mathfrak{a}_i,\mathfrak{a}'_k}=\frac{a_ia'_k}{xx'\,(M'-M)}\,,\ 
\widehat{N}_{\mathfrak{b},\mathfrak{a}'_k,\mathfrak{a}'_{k'}}=\frac{a'_ka'_{k'}}{x'{}^2\,(M'-M)}-\delta_{k,k'}\,,\nonumber\\
&&\widehat{N}_{\mathfrak{a}_i,\mathfrak{a}_{i'},\mathfrak{a}_{i''}}=\frac{a_ia_{i'}a_{i''}}{2x^3\,(M'-M)}+{N}_{\mathfrak{a}_i,\mathfrak{a}_{i'},\mathfrak{a}_{i''}}\,,\
\widehat{N}_{\mathfrak{a}_i,\mathfrak{a}_{i'},\mathfrak{a}'_k}=\frac{a_ia_{i'}a'_k}{2x^2x'\,(M'-M)}\,,\nonumber\\  
&&\widehat{N}_{\mathfrak{a}_i,\mathfrak{a}'_k,\mathfrak{a}'_{k'}}=\frac{a_ia'_ka'_{k'}}{2xx'{}^2\,(M'-M)}\,,\
\widehat{N}_{\mathfrak{a}'_k,\mathfrak{a}'_{k'},\mathfrak{a}'_{k''}}=\frac{a'_ka'_{k'}a'_{k''}}{2x'{}^3\,(M'-M)}-N'_{\mathfrak{a}'_k,\mathfrak{a}'_{k'},\mathfrak{a}'_{k''}}\,.
\nonumber\end{align}
QED\medskip

We call $\widehat{S},\widehat{T}$ the \textit{graft} of $S',T'$ onto $S,T$. 
Of course $\cD^\omega\Hg_\nu$ is the graft of $B_{m,2}$  onto
$\cD^{\pm,\omega}D_\nu$; here $M'=\mu$, $M=\nu^2$ (so $M'-M=4$), $\omega=\xi_3^{n^2+n-2}$, and 
$L=L'=2$. An easy example of (b) is grafting  $\cD^{-,0}D_\nu$ onto
$\cD^{+,\omega}D_\nu$ to form $\cD^{-,\omega}D_\nu$.

\smallskip
An easy generalisation of $\cD^\omega\Hg_\nu$ is to replace 
$\bbZ_\nu$ with any abelian group $K$ of odd order $\nu$. Take a trivial twist $\omega=0$, for 
simplicity. The entries of $S$ and $T$ involving  classes 1 and 4 are identical to
those in (\ref{genhS}). The role of $D_\nu$ 
is played by the semi-direct product $D_K=K\sdprod \bbZ_2 $, so the class 2 primaries
are labelled by the 2-dimensional irreps $i$ of $D_K$, and the class 3 ones by pairs
$(h,j)$ where $h$ runs through representatives of the cardinality-2 conjugacy classes of $D_K$ and $j$ runs 
through the irreps of $K$. The $T$-entries for class 2 are again 1's, while those for class 3 are
the evaluation $j(h)$. The submatrices $D$ and $E$ are again cosines, namely $D_{i,(h,j)}
=i(h)$ and $E_{(h,j),(h',j')}=2\,\mathrm{Re}(j'(h)\,j(h'))$.  Although the modular data thus generalises
naturally to arbitrary odd order abelian groups, we see in Section 4.1 that the subfactor
realisation for noncyclic $K$ is more subtle.

It is tempting to look for additional twists of $\cD\Hg_\nu$. After all,
we've suggested that the generalised Haagerup  is a twin of the dihedral $D_\nu$, and the latter can be twisted by $\bbZ_{2\nu}$ and not merely $\bbZ_\nu$.
The $\bbZ_\nu$ twist of both $\cD D_\nu$ and $\cD \Hg_\nu$ affects the class 3 primaries.
The independent $\bbZ_2$ twist of $\cD D_\nu$ affects the
class 4 primaries, so this suggests one should look for additional twists (by $\bbZ_2$
or perhaps $\bbZ_\mu$) of $\cD \Hg_\nu$, which sees the $m$ primaries of class 4
and leaves the other classes untouched.

We haven't gotten this to work. For instance, one
 approach to probing  some of the missing twists would be a
Galois automorphism acting on $S$ and $T$ entry-by-entry. 
The entries of $\cD^\omega \Hg_\nu$ lie in 
$\bbQ[\xi_{\nu^2\mu}]$, with Galois group  $\bbZ^\times_{\nu^2}
\times\bbZ^\times_\mu$. The class 3 twists $\omega\in\bbZ_\nu$ are stable under $\bbZ^\times_{\nu^2}$, 
as explained last subsection, so we should consider the effect of $\bbZ^\times_\mu$
on $\cD^\omega \Hg_\nu$. 

For $\ell\in\bbZ^\times_\mu$, $\sigma_\ell y=\pm y$ for some sign. If $\sigma_\ell y=-y$, then
$\sigma_\ell$ sends $\cD^\omega \Hg_\nu$ to \textit{nonunitary} modular data
with exactly the same fusion coefficients: the vacuum
primary is still the first one, but the positive column of $S$ is the second one. This means this
modular data won't have a subfactor realisation (it could perhaps have a planar algebra
interpretation, if the requirement of positive-definiteness there is dropped), though it can still have a VOA one.
The notion of twist should presumably preserve unitarity, so we won't regard this $\sigma_\ell
\cD^\omega\Hg_\nu$ as a twist. Nevertheless it may be interesting to search for 
its VOA realisations.

For this reason we should focus on those
 $\ell\in\bbZ^\times_\mu$ for which $\sigma_\ell y=+y$. Now,
it is typical for small $n$ that $\mu$ is a prime power (this is true for all $\nu<19$ except $\nu=9$).
In this case, any such $\sigma_\ell$ will send $\cD^\omega\Hg_\nu$ to equivalent modular
data. The search for new twists using Galois seems to fail. 

However when $\mu$ is not a prime power, this idea bears fruit. Given the prime decomposition
$\mu=\prod_{i=1}^k p_i^{m_i}$, where each $m_i>0$, then $\sigma_\ell\in \mathrm{Gal}(\bbQ[\xi_{\nu^2
\mu}]/\bbQ)\cong \bbZ_{\nu^2\mu}^\times$ fixes $\delta$ iff $\prod_{i=1}^k \left(\frac{\ell}{p_i}
\right)=+1$. We can also require here that $\ell\equiv_{\nu^2}1$ (the 
$\nu^2$-part of $\ell$ merely shuffles
the twists $\omega$). The automorphism $\sigma_\ell$ maps $\cD^\omega \Hg_\nu$ to
equivalent modular data iff $\left(\frac{\ell}{p_i}\right)=+1$ $\forall i$. Given any set $P$
(possibly empty) containing an even number of distinct prime divisors of $\mu$, pick any
$\ell_P\in\bbZ_{\nu^2\mu}^\times$ with $\ell_P\equiv_{\nu^2}1$ and $\left(\frac{\ell}{p_i}\right)=-1$
iff $p_i\in P$. By $\cD^{\omega,P}\Hg_\nu$ we mean the modular data $\sigma_{\ell_P}S,
\sigma_{\ell_P}T$; up to equivalence it is well-defined, and is inequivalent to $\cD^\omega
\Hg_\nu$ even though it has the same fusions. Thus as long as $k>1$,
i.e. $\mu$ is not a prime power, this construction yields new modular data.
In Section 4.1 we find subfactor realisations for this new Galois twist, at least when $\nu=9$ and $\nu=19$.

We learned in Section 2.4 that $\cD \Hg$ also sees $\cD \bbZ_3$, in some ways more
directly than it does $\cD S_3$. Since $\cD \bbZ_\nu$ has only a $\bbZ_\nu$ worth of
twists, it is certainly not inconceivable that the correct generic
twist group for $\cD \Hg_\nu$ is indeed $\bbZ_\nu$.

\subsection{Miscellanea involving modular invariants}

\cite[Sect.\ 6]{EP} studied the quantum double $S_3$ modular data in detail,
and similar techniques \cite{EP4} can be used for their twists, as noted in 
\cite{E2}. Using \cite[Theorem 4.3]{EP3},\cite{EP4} (with the correction to \cite{EP3} where
$Z_{11}$ was erroneously recorded as being sufferable, when we have now seen in
Section 2.2 of this paper that it is indeed a dual canonical modular invariant),
  we can summarize these results in the
following table (\textit{nimless} means there is no compatible nimrep, while \textit{nimble} means there is; \textit{sufferable} means there is a subfactor realisation):

$$\begin{array}{|c||c|c|c|c|} \hline
&\rm{nimless} & \rm{nimble, insufferable} & \rm{sufferable} & \rm{total}
\\ \hline
\cD^{+,0}S_3& 14&6&28&48\\ \hline
\cD^{-,1}S_3& 1&0&4&5\\ \hline
\cD^{+,1}S_3& 1&0&8&9\\ \hline
\cD^{-,0}S_3& 9&7&12&28\\ \hline
\cD^0\Hg_3 & ?& ?& 8\le ?\le 13&28\\ \hline
\end{array}$$

\centerline{\textbf{Table 1.} Modular invariants for the doubles of $Sym(3)$ and Haagerup}\medskip

The commutant of $S,T$ of course forms an algebra under matrix multiplication,
but the sufferable modular invariants themselves have a fusion structure, in
the sense that the product of two of them is always a linear combination over
nonnegative integers of other sufferable modular invariants. 
The modular invariants for $\cD^{s,\omega}D_\nu$ and $\cD^{s,\ell\omega}D_\nu$, 
including this fusion structure, are naturally isomorphic for any $\ell\in
\bbZ^\times_\nu$, using $\sigma_\ell\cD^{s,\omega}D_\nu\cong\cD^{s,\ell\omega}D_\nu$;
the same applies to $\cD^\omega\Hg_\nu$ and $\cD^{\ell\omega}\Hg_\nu$.

\cite{EP4,E2} remarked that there is an (injective) homomorphism from the 28 
$\cD \Hg$  modular invariants to the 28 sufferable  $\cD^{+,0}S_3$ ones.
We can define a natural injective homomorphism $\phi_{(S_3,s\omega)}^{(S_3,+0)}$ from the sufferable
twist $s,\omega$ quantum double $S_3$  modular invariants into the untwisted 
$S_3$ ones. We can actually embed the sufferable $-1$ into $+1$ modular 
invariants, the sufferable $-1$ into
the $-0$ modular invariants and the sufferable $-0$ into $+0$ modular invariants. 

This can be generalised and understood as follows. As is clear from Sections 3.1 and 3.2,
 primaries for both $\cD^{s,\omega}D_\nu$ and $\cD^\omega\Hg_\nu$ fall naturally into
 4 classes, and so their  modular invariants decompose 
 naturally into $4\times 4$ blocks as did their $S,T$ matrices. Write $S_{[i],[j]},T_{[i],[i]},M_{[i],[j]}$ for these blocks.

\medskip\noindent\textbf{Theorem 3.} \textit{Write $k$ for the number of distinct prime 
divisors of $\mu$. Define  a map
$$^\omega\phi_{\Hg \nu}^{D \nu} M=\left(\begin{array}{cccc}M_{[1],[1]}&M_{[1],[2]}&M_{[1],[3]}&0_{2\times 2}\\
M_{[2],[1]}&M_{[2],[2]}&M_{[2],[3]}&0_{n\times 2}\\M_{[3],[1]}&M_{[3],[2]}&M_{[3],[3]}&0_{n\nu\times 2}\\
0_{2\times 2}&0_{2\times n}&0_{2\times n\nu}&(a-b)I_{2\times 2}\end{array}\right)$$
where $M_{[1],[1]}=\left(\begin{array}{cc}a&b\\ c&d\end{array}\right)$. Then $^\omega\phi_{\Hg\nu}^{D\nu}$
is a bijection between all modular invariants of $\cD^\omega \Hg_\nu$ with $b\ne 0$, and
all modular invariants of $\cD^{-,\omega}D_\nu$ with $b\ne 0$. $^\omega\phi_{\Hg\nu}^{D\nu}$
is a $2^{k -1}$-to-1 surjection from all modular invariants of $\cD^\omega \Hg_\nu$ with $b= 0$, onto
all modular invariants of $\cD^{-,\omega}D_\nu$ with $b= 0$. Moreover, any modular invariant
of $\cD^{-,\omega}D_\nu$ is a modular invariant of $\cD^{+,\omega}D_\nu$ (but not conversely). 
}\medskip

\noindent\textit{Proof.} Let $M$ and $M'$ be modular invariants for $\cD^\omega\Hg_\nu$ and
$\cD^{-,\omega}D_\nu$, respectively. Let $S,T$ (resp. $S',T'$) denote the modular data
for $\cD^\omega\Hg_\nu$ (resp. $\cD^{-,\omega}D_\nu$). From $M'T'=T'M'$ we see that
$M'_{[i],[4]}=0=M'_{[4],[i]}$ for $i=1,2,3$. We also see that $M'_{[4],[4]}$ is diagonal.
Recall \eqref{Mselrule}: we see from the 0-row of $S$ that $s_\ell$ in classes [1],[2],[3]
is identically +1 for all $\ell\in\bbZ_{\nu^2\mu}^\times$, while  $s_\ell(\mathfrak{d}_l)
=\left(\frac{\ell}{\mu}\right)$ equals $-1$ when e.g. $\ell\equiv_\mu m$. Thus also
$M_{[i],[4]}=0=M_{[4],[i]}$ for $i=1,2,3$.

Evaluating $MS=SM$ at $([1],[4])$ and $([4],[1])$, we get two possibilities (identical conclusions
hold for $M'$):\smallskip

\noindent{(i)} $M_{[1],[1]}=I$ and $M_{[4],[4]}$ is a permutation matrix: $M_{\mathfrak{d}_l,\mathfrak{d}_{l'}}
=\delta_{\pi l,l'}$ for some permutation $\pi$ of $1\le l\le m$;

\smallskip\noindent{(ii)} $M_{[1],[1]}=1_{2\times 2}$ and $M_{[4],[4]}=0$.\smallskip

Looking at the remaining equations $MS=SM$, $MT=TM$, $M'S'=S'M'$ and $M'T'=T'M'$,
we see that (in case (i)) the equations involving $M_{[4],[4]}$ decouple from the others and reduce to
\begin{equation}\pi(l)\,\pi(l')\equiv_\mu\pm ll'\,,\ \pi(l)^2\equiv_\mu l^2\,.\label{m44S}\end{equation}
On the other hand, the equations for the unknown entries of $M_{[i],[j]}$, $i,j\ne 4$,
are identical to the corresponding equations  for $M'_{[i],[j]}$. For example, $(MS)_{[2],[4]}
=(SM)_{[2],[4]}$ and $(M'S')_{[2],[4]}=(S'M')_{[2],[4]}$ both reduce to $M_{\mathfrak{a}_i,0}
=M_{\mathfrak{a}_i,\mathfrak{b}}$ and $M'_{\mathfrak{a}_i,0}=M'_{\mathfrak{a}_i,\mathfrak{b}}$.
Therefore in case (ii), where $M_{[4],[4]}=0=M'_{[4],[4]}$, this means $^\omega\phi_{\Hg\nu}^{D\nu}$ is indeed a bijection.

In case (i), we know $M'_{[4],[4]}=I$, so we need to identify the permutation $\pi$.
Write $\pi 1=\ell$. Then \eqref{m44S} requires $\ell^2\equiv_\mu 1$, hence 
 $\pi l\equiv_\mu \pm l \ell$. So to any $1\le \ell\le m$ satisfying $\ell^2\equiv_\mu 1$,
define $\pi_{\ell} l$ to be the unique number $1\le \pi_{\ell}l\le m$ obeying $\pi_{\ell}l\equiv_\mu\pm
l \ell$. Then $\pi_{\ell}$ satisfies \eqref{m44S} and thus defines the $M_{[4],[4]}$-block
of a case (i) $\cD^\omega\Hg_\nu$  modular invariant. These $2^{k-1}$ 
$\ell$ parametrise the kernel of $^\omega\phi_{\Hg\nu}^{D\nu}$ in case (i). 

$\cD^{+,\omega}D_\nu$ is handled identically.\ QED\medskip

Note that $^\omega\phi_{\Hg\nu}^{D\nu}$ is linear and preserves matrix multiplication. 
In particular, when $\mu$ is a prime-power (as it is for the original Haagerup $\cD \Hg$, and
all $\nu<19$ except $\nu=9$), there will be an algebra isomorphism between
the span of modular invariants for $\cD^\omega \Hg_\nu$ and for $\cD^{-,\omega}D_\nu$, and an (injective but nonsurjective)
algebra homomorphism of these into $\cD^{+,\omega}D_\nu$.

The relation between the 28 sufferable $\cD^{+,0}S_3$ modular invariants and
the 28 $\cD^{0}\Hg_3$ modular invariants would seem to be a coincidence:
e.g. at $\nu=1$, $\cD^{+,0}D_1=\cD^{+,0}\bbZ_2$ has exactly
6 modular invariants, all sufferable, while $\cD^{0}\Hg_1$ has exactly 2.

The dihedral group is only half the story; the other half is the affine algebra data $B_{m,2}$.
Its modular invariants were classified in \cite{Gorth}:
when $\mu$ is not a perfect square (the situation here), the complete list is
$B(d_0,\ell_0)$ and $B(d_1,\ell_1;d_2,\ell_2)=B(d_2,\ell_2;d_1,\ell_1)$ where $d_i|
\mu$, $\mu|d_i^2$, and $1\le \ell_i\le m$ obeys $\ell_i^2\equiv_{d_i^2
/\mu}1$. for matrices $B(d,\ell)$ and $B(d_1,\ell_1;d_2,\ell_2)$ 
defined in \cite{Gorth}. Those modular invariants which are permutation 
matrices are the $B(\mu,\ell)$ for $1\le \ell\le m$ satisfying
$\ell^2\equiv_\mu 1$: their nonzero entries are $B(\mu,\ell)_{z,z}=1$ for $z
\in\{0,\mathfrak{b},\mathfrak{a}_i\}$,
and $B(\mu,\ell)_{\mathfrak{d}_l,\mathfrak{d}_{l'}}=1$ for $l'\equiv_\mu\pm 
\ell l$. We see from this and Theorem 3 that the modular invariants of
$\cD^\omega\Hg_\nu$ come from those of $\cD^{-,\omega}\Hg_\nu$ and
$B_{m,2}$ in a very direct sense.

Recall \cite{ost,Ev} that the sufferable modular
invariants of the double of a finite group $G$ are parametrised by pairs $(H,\psi$) where
$H$ is a subgroup of $G\times G$ and $\psi\in H^2(B_H;\bbC^\times)$ is called
\textit{discrete torsion}. This suggests:

\medskip\noindent\textbf{Question 4.} \textit{Find some analogue for $\cD^0\Hg_\nu$ of
this $(H,\psi)$ parametrisation of
sufferable modular invariants.}\medskip

This $(H,\psi)$ parametrisation belongs most naturally to the $\cD\Delta$ formulation of the double of the group $G$.
As mentioned in Section 1.2, $\cD\widehat{\Delta}$ is also  the double of $G$.
Is there also a parametrisation of the sufferable modular invariants of the double of $G$
which is more natural in the $\cD\widehat{\Delta}$ language?

Of special importance (recall Section 1.3) are the \textit{monomial} modular invariants
of $\cD^0\Hg_\nu$. Three of these are obvious:
\begin{equation}
|ch_0+ch_{\mathfrak{b}}+2\sum ch_{\mathfrak{a}_i}|^2\,,\
|ch_0+ch_{\mathfrak{b}}+\sum ch_{\mathfrak{a}_i}+\sum ch_{\mathfrak{c}_{h,0}}|^2\,,\
|ch_0+ch_{\mathfrak{b}}+2\sum ch_{\mathfrak{c}_{h,0}}|^2\,.\label{monom}\end{equation}

\noindent\textbf{Proposition 4.} \textit{When $\nu=p$ or $pq$ for (not necessarily distinct)  primes
$p,q$, the only monomial modular invariants are the three in \eqref{monom}.}

\medskip\noindent\textit{Proof.}  We need to find all eigenvectors $u$ of $S,T$ with eigenvalue 
1, with $u_0=1$ and all other $u_x\in\{0,1,2,\ldots\}$, so $Z=|\sum_x u_xch_x|^2$. We know
$u_{\mathfrak{b}}=1$ since $M_{\mathfrak{b},\mathfrak{b}}=1$ for any modular invariant
(clear from the proof of Theorem 3), hence $u_{\mathfrak{d}_l}=0$ for all $l$. We need to
determine $u_i:=u_{\mathfrak{a}_i}$ and $u_{h,j}:=u_{\mathfrak{c}_{h,j}}$.

$Tu=u$ implies $u_{h,j}=0$ unless $\nu$ divides $hj$. Therefore if $u_{h,j}\ne 0$ we have:

\noindent{(i)} for $\nu$ prime, $j$ must be 0;

\noindent{(ii)} for $\nu=p^2$, only for $j=0$ or $p|h$ and $p|j$;

\noindent{(iii)} for $\nu=pq$ ($p\ne q$), only for $j=0$, or $p|h$ and $q|j$, or $q|h$ and $p|j$.

$Su=u$ implies, for all $1\le h,h''\le n$ and $0\le i<\nu$, 
\begin{align} &&2+2\sum u_{j}+2\sum u_{h',j}\cos(2\pi
ih'/\nu)=\nu u_i \,,\label{sub}\\ &&2+2\sum\cos(2\pi jh/\nu)+2\sum\cos(2\pi(
hj+h'i)/\nu)u_{h',j}=\nu u_{h,i}\,,\label{succ}\\ &&\sum\cos(2\pi hj/\nu)\,u_{h'',j}=\sum\cos(2\pi h''j/\nu)\,u_{h,j}\,,\label{suc}\end{align}
where we sum over $h',j$ (\eqref{suc} was obtained by hitting  \eqref{succ} with cos$(2\pi h''i/\nu)$ and summing over $i$).
In analysing these equations, it is useful to recall (for $k$ odd)
\begin{equation}
\sum \xi_k^\ell=2\sum \cos(2\pi l/k)=\mu(k)\,,\end{equation}
the sums taken over the $1\le \ell\le k$ and $1\le l< k/2$ which are coprime to $k$. 
The M\"obius function $\mu(k)$ 
 equals $(-1)^n$ if $k$ is a product of $n\ge 0$ distinct primes, and 
0 otherwise.
The Galois symmetry $\sigma_\ell$ for $\ell\in\bbZ^\times_\nu$ (recall
\eqref{galsym}) maps 
$\mathfrak{a}_i$ to $\mathfrak{a}_{\pm\ell i}$ and $\mathfrak{c}_{h,j}$ to $\mathfrak{c}_{\pm\ell h,
\ell j}$. Therefore $u_i=u_{\mathrm{gcd}(i,\nu)}$ and $u_{h,j}=u_{\mathrm{gcd}(h,\nu),jh/
\mathrm{gcd}(h,\nu)}$, so $u_{h,-j}=u_{h,j}$.

For $\nu$ prime, this means $u_i=u_1$ and $u_{h,j}=u_{1,0}\delta_{j,0}$. Plugging this into
\eqref{sub} gives $u_1+u_{1,0}=2$, which correspond to the three
 solutions \eqref{monom}.

For $\nu=p^2$, we read off from \eqref{suc} with $h=p,h'=1$ that $u_{p,jp}=u_{p,0}-u_{1,0}$
for all $j\ne 0$. Then comparing \eqref{sub} for $i=1,p$ gives $2=u_p+u_{1,0}=u_1+u_{p,0}$.
We can now force $u_p=u_1$ by \eqref{succ} at $h=p,i=1$, and again we recover \eqref{monom}.

Finally, turn to $\nu=pq$.  \eqref{suc} with $h=1$ and $h'=p,q$ gives respectively
$u_{p,qj}=u_{p,0}-u_{1,0}$, $u_{q,pj}=u_{q,0}-u_{1,0}$ for $j\ne 0$. Then \eqref{sub} at
$i=1,p,q$ gives $u_q+u_{q,0}=u_p+u_{p,0}=2$ and $u_{1,0}=2+u_1-u_p-u_q$. Now
\eqref{succ} at $h=i=p$ and $h=i=q$ give $u_1=u_p=u_q$ and $u_{1,0}=u_{p,0}=u_{q,0}$
and we're done. QED\medskip

The $\nu=3$ modular invariants \eqref{Z3S3},\eqref{Z3Hg} generalise,
for any $\cD^{\pm,\omega}D_\nu$, $\cD^\omega\Hg_\nu$, to
\begin{equation}\label{modinvDnuZnu}
|ch_0+ch_{\mathfrak{b}}|^2+2\sum|ch_{\mathfrak{a}_i}|^2+2\sum|ch_{\mathfrak{c}_{h,j}}|^2\,.
\end{equation}
This corresponds to the $\bbZ_\nu$ subsystem in
both $\cD D_\nu$ and $\cD^0\Hg_\nu$, and to the VOA realising $\cD^\omega\bbZ_\nu$
and containing both  $\cD^{+,\omega}D_\nu$ and $\cD^\omega\Hg_\nu$ VOAs.

\section{Subfactors for Haagerup--Izumi modular data}

\subsection{Izumi's subfactors and their modular data}

In Section 5 we address the question of realising this generalised Haagerup modular data
$\cD^\omega\Hg_\nu$ by
VOAs. In this section we address their subfactor realisation. 

In Section 7 of \cite{iz3}, Izumi suggests associating subfactors to any odd abelian group $K$ of 
order ${\nu}$,
using endomorphisms in the Cuntz algebra $\mathcal{O}_{\nu+1}$. (Warning: Izumi's $n$ is our 
$\nu=2n+1$.)  Let $\cA(K)$ be the set of all  $\nu\times\nu$ complex matrices $A=(A_{g,h})$, $g,h\in K$, satisfying
\begin{eqnarray}
&&A_{g,h}=\overline{A_{h,g}}\,,\ A_{g,h}=A_{-h,g-h}=A_{h-g,-g}\,,\
A_{g,0}=\delta_{g,0}-1/(\delta-1)\,,\label{7.1}\\ &&\sum_{m\in K}A_{g+m,h}\overline{A_{m,h}}=\delta_{g,0}
-\delta_{h,0}/\delta\,,\label{7.4}\\ &&\sum_{m\in K}A_{m,g+h}\overline{A_{m+k,g}A_{m+l,h}}=
A_{g+l,k}A_{h+k,l}-\delta_{g,0}\delta_{h,0}/\delta\,,\label{7.5}\end{eqnarray}
for all $g,h,k,l\in K$ (recall $\delta=(\nu+\sqrt{\nu^2+4})/2$). These equations imply \cite{iz3}
\begin{eqnarray}
&&|A_{g,h}|=\sqrt{\delta}\,,\\
&& A_{g+k,h}A_{h,k}=A_{-k,g}A_{g,h-k}\,,\label{7.7}\end{eqnarray}
for any $g,h,k\in K$ with $g\ne 0,\,h\ne 0,\,g\ne h$.
Izumi shows that to any matrix $A\in \cA(K)$, 
 there corresponds a (nonbraided) subfactor of index $\delta+1$ with principal
graph the $\nu$-star (see the top graphs of Figures 2 and 4). We call these \textit{subfactors of
Izumi type $K$.} Next subsection we
determine much of the data of these subfactors. 
His Theorem 8.4 identifies the modular data $S$, $T$ of the
even part of the quantum double of this subfactor in terms of $m(\nu^2+1)$ variables
 $\omega_l,C^l_{g,h}$, $g,h\in K$, $1\le l\le m$, satisfying the $m(\nu^2+1)$ equations:
 \begin{equation}
 \sum_{g\in K}C_{0,g}^l=\omega_l-\overline{\omega_l}/\delta\,,\
 \omega_lC^l_{g,h}-\sum_{k\in K} A_{g+k,2h}C^l_{h,k}= \delta_{h,0}\overline{\omega_l}/\delta\,. \label{8.10}\end{equation}
Some  solutions to \eqref{8.10}, occurring when $\omega^\nu_l=1$,
  are redundant (i.e. correspond to
$\mathfrak{b},\mathfrak{a}_i$ or $\mathfrak{c}_{h,j}$) and should be 
ignored. In particular, 
precisely $\sigma_1(\nu)-1$ (i.e. the sum of divisors $d>1$ of $\nu$) 
solutions to \eqref{8.10} with $\omega_l= 1$ are redundant, and 
 for any other root of unity $\omega_l$ of order $d$ dividing $\nu$, 
the number of redundant solutions is $\sum d'\phi(\nu/d')/2$ where $\phi$ is
Euler's totient and the sum is over all divisors $d'<d$ of $d$. The 
entries of the $S$ and $T$
matrices equal those of $\cD^0\Hg_\nu$ except possibly for those in the bottom-right
$m\times m$ block, which are given by $T_{\mathfrak{d}_l,\mathfrak{d}_l}
=e^{-2\pi\i c/24}\omega_l$ and 
\begin{equation}S_{\mathfrak{d}_l,\mathfrak{d}_{l'}}=\frac{1}{\nu\delta+2}\left(\omega_l\omega_{l'}
+\delta\sum_{g,h\in K}\overline{C^l_{g,g+h}C^{l'}_{-g,h}}\right)\,.\label{8.12}\end{equation}
When some $\omega_l$ has more than one solution in \eqref{8.10}, \eqref{8.12} 
is ambiguous as it isn't obvious which solutions to \eqref{8.10} to use.

Izumi shows in his Example 7.1 that $\cA(\bbZ_3)$ contains exactly 
two matrices, both corresponding to
the Haagerup subfactor. In Appendix C he exactly 
solves \eqref{8.10} for $A\in \cA(\bbZ_3)$, and in this 
way obtains (a complicated expression for) $\cD\Hg$.

In his Example 7.2, Izumi shows $\cA(\bbZ_5)$ contains exactly 4 matrices, 
again corresponding to a single (new) subfactor.
Izumi does not solve \eqref{8.10} for them.
This analysis of \eqref{7.1}-\eqref{7.5} for $\nu=3$ and $\nu=5$ is as far as Izumi went (this is what 
we mean when we call it a \textit{hypothetical} family of subfactors).

Nevertheless we can push a little further this analysis. The group Aut$(K)$ acts naturally
on the matrix $A$ by shuffling the entries: $\alpha\in\mathrm{Aut}(K)$ sends $A$ to 
$A^\alpha$ defined by $A^\alpha_{g,h}=A_{\alpha g,\alpha h}$. For example, when 
$K=\bbZ_\nu$,
Aut$(K)=\bbZ_\nu^\times$ acts by multiplication. In this way, given any $A\in \cA(K)$, Aut$(K)$
embeds into Aut$(\bbQ[A])$ and fixes $\delta$, where $\bbQ[A]$ is the field generated over
$\bbQ$ by all entries $A_{g,h}$. Call $A,A'\in \cA(K)$ \textit{equivalent} if 
$A'=A^{\alpha}$ for some $\alpha\in\mathrm{Aut}(K)$. Equivalent matrices give rise to
equivalent subfactors and equivalent modular data.
For later convenience, for $x\in\bbQ[A]$, write $\Tr\, x$ for the orbit sum $\sum_{y\in \langle x\rangle}y$ over the $K$-orbit $\langle x\rangle$.

It is very difficult to solve \eqref{8.10} directly (see Appendix C of \cite{iz3} for the $\nu=3$
argument, which is already quite involved). Our strategy for finding the modular data is
to reduce it to computer calculations. A key observation is that if $S,T$ are the modular data
of $\cD^0\Hg_\nu$, and if $S',T$ is a second modular data, with
\begin{equation}\label{Sapprox}
\left|S'_{\mathfrak{d}_l,\mathfrak{d}_{l'}}-S_{\mathfrak{d}_l,\mathfrak{d}_{l'}}\right|<\frac{1}{18\sqrt{\mu}}\end{equation}
for all $1\le l,l'\le m$, and all other entries of $S'$ and $S$ are equal, then \eqref{SNT} implies
$S'=S$ everywhere.
Maple 11 (x86 64 Linux) was used in these calculations. It is difficult to rigourously analyse 
the resulting error, but various consistency checks (e.g. unitarity of the
numerically obtained $S'$ matrix, or the values of $|S'_{\mathfrak{d}_l,\mathfrak{d}_{l'}}-
S_{\mathfrak{d}_l,\mathfrak{d}_{l'}}|$) all indicate it is on the order of $10^{-7}$ or better,
far more precision than is actually needed. More details are given below.
Computer calculations were \textit{not} used in the determination of $\cA(K)$
in Theorem 5.

\medskip\noindent\textbf{Theorem 5.} (a) \textit{There is a unique
$A\in \cA(1)$. Its modular data is $\cD^0
\Hg_1$.}

\smallskip\noindent(b) \textit{There are precisely
2 matrices in $\cA(\bbZ_3)$. They are equivalent and realise $\cD^0\Hg_3$.}

\smallskip\noindent{(c)} \textit{There are precisely 4 matrices  in $\cA(\bbZ_5)$. 
They are equivalent and realise $\cD^0\Hg_5$.}

\smallskip
\noindent{(d)} \textit{There are precisely 6 matrices in $\cA(\bbZ_7)$. 
They are equivalent and realise $\cD^0\Hg_7$.}

\smallskip\noindent{(e)} \textit{There are precisely 12 matrices in $\cA(\bbZ_9)$. They form
2 equivalence classes and define inequivalent subfactors; one realises the $T$ matrix of $\cD^0\Hg_9$, and the other
realises the $T$ of $\cD^{0,\{5,17\}}\Hg_9$ (defined in Section 3.3).}

\smallskip\noindent{(f)} \textit{$\cA(\bbZ_3\times\bbZ_3)$ is empty.}

\medskip\noindent\textit{Proof.} Let's begin by identifying $\cA(K)$ for these groups (done
in \cite{iz3} for $K=\bbZ_3,\bbZ_5$). It is trivial that $\cA(1)=\left\{\left(\frac{\delta-2}{\delta-1}\right)
\right\}$, and for $\bbZ_3$ and  $\bbZ_5$ \eqref{7.1} fixes everything in terms of complex numbers $a$ and $a,b$,
 respectively:
\begin{equation}A(\bbZ_3) ={\scriptsize \frac{1}{\delta-1}\left(\begin{array}{ccc}\delta-2&
-1&-1\\-1&-1&a\\-1&\overline{a}&-1\end{array}\right)}\ ,\ 
A(\bbZ_5)={\scriptsize \frac{1}{\delta-1}\left(\begin{array}{ccccc}\delta-2&-1&-1&-1&-1\\
-1&-1&a&b&a\\-1&\overline{a}&-1&b&b\\-1&\overline{b}&\overline{b}&-1&a\\-1&\overline{a}&
\overline{b}&\overline{a}&-1\end{array}\right)}\,.\nonumber\end{equation}
 For $\nu=3$, \eqref{7.5} forces  $a^2-a+\delta=0$. The two roots of this polynomial are
 complex conjugates, and indeed $-1\in\bbZ_3^\times$ here acts as $a\leftrightarrow\overline{a}$.
Similarly, the generator $2\in\bbZ^\times_5$ acts on $A(\bbZ_5)$ by sending $a\mapsto b\mapsto
\overline{a}\mapsto\overline{b}\mapsto a$. Then $a$ has minimal polynomial 
$x^4-\delta x^3+(3\delta-2)x^2-\delta^2x+\delta^2=0$ and the other
3 roots of this polynomial are $\overline{a},b,\overline{b}$ --- indeed, the Galois group over
$\bbQ[\delta]$ of this polynomial is $\bbZ_5^\times$ with this action.
Given $a$, this and the equation 
$1-a-\overline{a}+a\overline{b}+\overline{a}b=0$ then determine $b$.

\smallskip
\noindent\textit{Computation of $\cA(\bbZ_7)$} is handled similarly. \eqref{7.1} and \eqref{7.7} give us
\begin{equation}
A(\bbZ_7)={\scriptsize \frac{1}{\delta-1}\left(\begin{array}{ccccccc}\delta-2&-1&-1&-1&-1&-1&-1\\
-1&-1&a&b&c&b&a\\-1&\overline{a}&-1&b&d&d&b\\-1&\overline{b}&\overline{b}&-1&c&d&c\\-1&
\overline{c}&\overline{d}&\overline{c}&-1&b&b\\ -1&\overline{b}&\overline{d}&\overline{d}
&\overline{b}&-1&a\\ -1&\overline{a}&\overline{b}&\overline{c}&\overline{b}&
\overline{a}&-1\end{array}\right)}\,.
\end{equation}
The generator $3\in\bbZ^\times_7$ sends $a\mapsto c\mapsto d\mapsto \overline{a}\mapsto
\overline{c}\mapsto\overline{d}\mapsto a$ and $b\mapsto\overline{b}\mapsto b$.
\eqref{7.7} gives $ad=bc$, while \eqref{7.4} with $(g,h)=(1,1),(3,1),(2,3)$ give
$$(b+\overline{b}+2)(c+\overline{c})=-1+\Tr\, a
=(b+\overline{b}+2)(a+\overline{a})=(b+\overline{b}+2)(d+\overline{d})\,,$$
where the last two equalities come from the $\bbZ^\times_7$ symmetry.
If $b+\overline{b}\ne -2$, then $a=d=\overline{c}$ and $\delta b=a^3$ and we get two
incompatible equations for $s:=a+\overline{a}$, namely
 $s^4+s^2(\delta-4)-\delta s=\delta-2$
and $s^3+(\delta-3)s-\delta s^2=5\delta+1$.

Thus $b+\overline{b}=-2$ (hence $b=-1\pm\i\sqrt{\delta-1}$) and $\Tr\, a
=1$. Let $t=a+d+\overline{c}$. Then $2\mathrm{Re}(t)=1$ and (from \eqref{7.4})
$2\mathrm{Re}(t\overline{b})=-3\delta-5$ gives $2t=1+(b+1)(3-\delta)$ and we obtain
$a^3-ta^2+b\overline{t}a-\delta b=0$. The other two roots of this polynomial are $d$ and $\overline{c}$
(which to call $d$ follows from $2\,\mathrm{Re}\,(b+c-a\overline{b})=\delta$),
and its Galois group over $\bbQ[\delta,b]$ is $\langle 3^2\rangle<\bbZ^\times_7$.
\smallskip

\noindent\textit{Computation of $\cA(\bbZ_9)$}. We get as before
\begin{equation}
A(\bbZ_9)={\scriptsize \frac{1}{\delta-1}\left(\begin{array}{ccccccccc}
\delta-2&-1&-1&-1&-1&-1&-1&-1&-1\\
-1&-1&a&b&c&d&c&b&a\\ -1&\overline{a}&-1&b&e&f&f&e&b\\ 
-1&\overline{b}&\overline{b}&-1&c&f&g&f&c\\-1&\overline{c}&\overline{e}&\overline{c}&-1&d&f&f&d
\\-1&\overline{d}&\overline{f}&\overline{f}&\overline{d}&-1&c&e&c\\
-1&\overline{c}&\overline{f}&\overline{g}&\overline{f}&\overline{c}&-1&b&b\\
-1&\overline{b}&\overline{e}&\overline{f}&\overline{f}&\overline{ e}&\overline{b}&-1&a\\
-1&\overline{a}&\overline{b}&\overline{c}&\overline{d}&\overline{c}&\overline{b}&\overline{a}&-1
\end{array}\right)}\,.\end{equation}
Moreover, \eqref{7.7} yields $gb=cf$ and $a\overline{c}=b\overline{e}=d\overline{f}$.
The generator $2 \in\bbZ^\times_9$ sends $a\mapsto e\mapsto d\mapsto 
\overline{a}\mapsto\overline{e}\mapsto\overline{d}\mapsto a$, $b\mapsto f\mapsto\overline{c}
\mapsto\overline{b}\mapsto\overline{f}\mapsto c\mapsto b$, and $g\mapsto\overline{g}\mapsto g$.

Up to the relations $gb=cf$, $a\overline{c}=b\overline{e}=d\overline{f}$, there are precisely
17 $\bbZ_9^\times$ orbits in $b,\overline{b},c,\overline{c},f,\overline{f}$ of degree $\le 3$.
Using \eqref{7.4} with $g=3$ and \eqref{7.5} when 3 divides both $g$ and $h$, we can compute
these 17 orbit sums $\Tr$ as linear expressions in $s:=\Tr\, b=2\,\mathrm{Re}(b+c+f)$.
For example, $\Tr\, g= 1-s$, $\Tr\, (bg)=3(s-\delta)$, $\Tr\, (b\overline{g})=2\delta+1+s\,(2-\delta)$.
But $\Tr\, b \,\Tr\, g=\Tr\, (bg)+\Tr\, (b\overline{g})$, so $s^2+(4-\delta)s+1-\delta=0$. Choose one
of these 2 solutions for $s$ (they both yield $A\in\cA(\bbZ_9)$).

From $g+\overline{g}=1-s$ we get $g=(1-s\pm\sqrt{(\delta-6)s-3\delta})/2$. Again select one
of these $g$ (they are complex conjugates). Let $t:=b+\overline{c}+\overline{f}$. Then 
$t+\overline{t}=s$ and $tg+\overline{tg}=3(s-\delta)$, which fixes $t=s/2\pm(3\delta-6-\delta s-s)
\sqrt{(\delta-6)s-3\delta}/36$ (use the same sign as in $g$). Then $b\overline{c}+b\overline{f}
+\overline{c}\overline{f}=\overline{tg}$ and $b\overline{cf}=\delta\overline{g}$, so $b,\overline{c},
\overline{f}$ are the 3 roots of $x^3-tx^2+\overline{gt}x-\delta\overline{g}=0$. Let $b$ be any
of these roots; which to call $\overline{c}$ and $\overline{f}$ is fixed by the relation
$\Tr\, (b^2\overline{c})=6\delta+3+(8-\delta)s$.

Now \eqref{7.4} at $(g,h)=(1,1)$ and (4,2) gives $a^2+va/u+\delta\overline{u}/u=0$
for $(u,v)=(\overline{b}+\overline{c}^2/\overline{f}-1,1+\overline{b}c+b\overline{c})$
and $(f+\overline{c}-\overline{bc}/\delta,\delta-f-\overline{f})$ respectively;
subtracting these expresses $a$ rationally
in terms of $b,c,f$.
%$a^2+va+\delta\overline{u}/u$ where $u=f+\overline{c}-
%\overline{bc}/\delta$ and $v=(\delta-f-\overline{f})/u$. Take either solution. 
Then $d=af/c$
and $e=bc/a$. These values of $a,\ldots,g$ will solve \eqref{7.1}-\eqref{7.5}. There are
2 choices for $s$, 2 for $g$, and 3 for $b$. The $\bbZ^\times_9$ symmetry accounts
for this $g$ and $b$ ambiguity, leaving 2 inequivalent solutions.

\smallskip\noindent\textit{Computation of  $\cA(\bbZ_3^2)$}. 
Write the $(1,0)$-row of the matrix $(\delta-1)A(\bbZ_3^2)$
as $(-1,-1,t,u,v,w,x,y,z)$. \eqref{7.7} with $g=h=(1,0)$ and $k=(0,1),(1,1)$ gives
$y\overline{u}=x\overline{v}=z\overline{w}=:\delta a$ for some complex number $a$ with $|a|=1$. 
From \eqref{7.7} with $g=(1,1),h=(1,0),k=(1,2)$ yields $\delta t=xyz$, and applying the
$\bbZ_3^2$-automorphism $(i,j)\mapsto (i,-j)$ to this we get $\delta t=uvw$.
Comparing these gives $a^3=1$.

Putting $g=(0,1),h=(1,0)$ into \eqref{7.4} yields
\begin{equation}\label{z3z3}
1-\delta=(\overline{a}u+a\overline{u}-1)(\overline{a}v+a\overline{v}-1)\,.\end{equation}
Hitting this with the $\bbZ_3^2$-automorphism $(i,j)\mapsto(i+j,j)$ gives
the same equation with $w$ replacing $u$. Thus $\overline{a}u+a\overline{u}-1=\overline{a}v+
a\overline{v}-1$ and \eqref{z3z3} says there is a real number whose square is negative.
This impossibility establishes part (f).

\smallskip\noindent\textit{Determining the modular data.} Consider first $K=\bbZ_3$ and fix
$A\in\cA(\bbZ_3)$.  Let $\omega_l,C_{g,h}^l$ be some solution to \eqref{8.10}. Then we learn
from the above proof that any automorphism $\sigma\in\mathrm{Gal}(\overline{\bbQ}/\bbQ[\delta])$ 
corresponds to an $\ell\in \bbZ_3^\times$ in such a way that $\sigma A_{g,h}=A_{\ell g,
\ell h}$.  Certainly $\sigma$ sends $\omega_l$ to some other root of unity $w_{l^\sigma}$
of the same order $d$, and moreover $\sigma\overline{\omega_l}=\overline{\omega_{l^\sigma}}$
(since $\sigma$ commutes with complex conjugate in cyclotomic fields).
Hence $\sigma$ sends the solution $\omega_l,C_{g,h}^l$ of \eqref{8.10} to
another solution $\omega_{l^\sigma},C^{l^\sigma}_{g,h}$ given by $\sigma(C^l_{g,h})=
C^{l^\sigma}_{\ell g,\ell h}$. Incidentally, \eqref{8.12} then implies the Galois
actions 
\begin{equation}\sigma(S_{\mathfrak{d}_l,\mathfrak{d}_{l'}})=S_{\mathfrak{d}_{l^\sigma},
\mathfrak{d}_{l'{}^{\sigma}}}\,,\ \sigma(T_{\mathfrak{d}_l,\mathfrak{d}_{l}})=
T_{\mathfrak{d}_{l^\sigma},\mathfrak{d}_{l{}^{\sigma}}}\,\nonumber\end{equation}
(the second equation is simply the statement $\sigma\omega_l=\omega_{l^\sigma}$; to see
the first, compute instead $\sigma\overline{S}$).
This is stronger than the usual Galois action \eqref{galois} in modular data,
though of course it is compatible with $\cD^0\Hg_\nu$. 

Since $\delta\in\bbQ[\xi_\mu]$, when $\mu$ is coprime to the order $d$ of $\omega_l$
 any automorphism in $\mathrm{Gal}(\bbQ[\xi_d]/\bbQ)$ will lift to a $\sigma$ 
fixing $\delta$ (when gcd$(d,\mu)>1$, only half will).
So if \eqref{8.10} can be solved for a $d$th root of unity $w_l$, for $d$
coprime to $\mu$, then it can be solved for \textit{all} $d$th roots of 1.
That is, we would have at least $\phi(d)$ solutions. As long as $d$ doesn't
divide $\nu$, none of these are redundant. But Izumi tells us there are 
precisely $m$ nonredundant (independent) solutions to \eqref{8.10}
(corresponding to the $m$ half-braidings, giving rise to the $m$ primaries
of type $\mathfrak{d}_l$).

Hence any solution to \eqref{8.10} involves $\omega_l$ of order $d$
where $\phi(d)\le m$ (for $d$ coprime to $\mu$) or $\phi(d)\le 2m$ (otherwise).
There is thus a finite number of possibilities for $\omega_l$. A little
thought reduces the number further --- e.g. when gcd$(d,\mu)=1$ it suffices
to consider only $\omega_l=\xi_d$. We find that (for $\nu=3$)  only
12 possibilities for $\omega_l$ need be considered, namely $\pm 1, \pm\xi_3,\i,\pm\xi_7,\pm \xi_9,
\pm\xi_{13},-\xi_{13}^6$. To rule out these
possibilities,  we compute
determinants. In particular, for $\omega_l\ne 1,\xi_3$, we show that the linear system \eqref{8.10} 
is inconsistent by showing  the $(\nu^2+1)\times(\nu^2+1)$ augmented matrix (formed from 
the coefficients and constant vector) has nonzero determinant. To show $\omega_l=1$
(respectively $\omega_l=\xi_3$) has at most 3 (resp. 1) independent solutions, we 
deleted 3 (resp. 1) equations and 2 (resp. 0) variables from \eqref{8.10} so that the resulting
reduced coefficient matrix has nonzero determinant.
We evaluated these determinants numerically using Maple. If some
determinant was fairly close to 0, we found some Galois associate of $\omega_l$
with large determinant (i.e. of order 1 or higher).

In this way we show that $\omega_l$ must be one of the values $\xi_{\mu}^{ml^2}$. We can now
determine the modular data, rigourously identifying it with
$\cD^0\Hg_\nu$, as follows.  Numerically
compute the solution to \eqref{8.10} and plug it into \eqref{8.12} to obtain 
an estimate $S'$. We find that $S'$ matches $S$ to within $2\times 10^{-9}$, well within
the 0.015 required by \eqref{Sapprox}.

This argument for $K=\bbZ_5$ and $\bbZ_7$ is identical. For them exactly 14 and 22, respectively,
values of $\omega_l$ need to be considered. The error $|S'-S|$ is less than $2\times 10^{-9}$
and $5\times 10^{-8}$ respectively, again well within \eqref{Sapprox}.

The modular data calculation for $K=\bbZ_9$ is more complicated for two reasons, both
related to the compositeness of $\mu=5\cdot 17$. Fix $A\in\cA(\bbZ_9)$ (the result for the
other $A$ is automatic by Galois considerations). First, an automorphism $\sigma$ will
permute the entries of $A$ iff $\sigma$ fixes both $\delta$ and $s=\Tr\,b$ (incidentally,
Gal$(\bbQ[s]/\bbQ)\cong\bbZ_2^2$, so $s$ lies in a cyclotomic
extension of $\bbQ$). Hence potentially only 1/4 of Gal$(\bbQ[\xi_d]/\bbQ)$ 
lifts to such $\sigma$. This means more possibilities for $\omega_l$ to be
eliminated (218 possibilities is overkill), and we find that the $T$ matrix for $A$ exactly matches
that of either $\cD^0\Hg_9$ or $\cD^{0,\{5,17\}}\Hg_9$. More significantly, because of equalities 
like $\omega_1=\omega_{16}$, there won't always be unique solutions to \eqref{8.10},
and \eqref{8.12} is ambiguous. This means we cannot identify $S$ without more work. 
QED\medskip

The key difference between $K=\bbZ_9$ and $\bbZ_3^2$ is the size of their automorphism
groups (6 vrs 48). Nonuniqueness for $K=\bbZ_9$ arises because $\mu=85$ is 
composite: more generally, where $k>1$ is the number of distinct prime divisors of $\mu$,
 there will be at least $\left({k\atop 2}\right)$ inequivalent $A$ and indeed
inequivalent subfactors, realising the inequivalent modular data of type $\cD^{0,P}\Hg_\nu$
defined in Section 3.3.

This nonuniqueness of subfactors with identical principal graphs is not a 
surprise: for a simple example,  the principal graph fails to uniquely
identify the subfactor even in index $<4$ {(two subfactors realise each of  the graphs $E_6,E_8$ 
although one is the opposite of the other) whilst in index $4$ the affine graphs $A^{(1)}_{2n-5}$ 
and its orbifold $D^{(1)}_n$ both have  $n-2$ inequivalent subfactors  \cite{IK}.}

Exact expressions for the matrices $A\in \cA(K)$ appear in the proof of Theorem 5. 
Numerical estimates for these $A$ may also be of interest. A convenient way to express
this, for $K$ cyclic, is in terms of $j_2,j_3,\ldots,j_{n+1}\in\bbR$: for $0<g<h<\nu$ we have
$$A_{g,h}=\frac{\sqrt{\delta}}{\delta-1}\exp[\i(j_h-j_g-j_{h-g})]\,,$$
where $j_1=0$ and $j_{n+1+i}=j_{n+1}+j_n-j_{n-i}$ for $1\le i<n$
(see Lemma 7.3 of \cite{iz3}).
\begin{align}
j_2^{(3)}\approx&\, 1.292076\,;\nonumber\\
(j_2^{(5)},j_3^{(5)})\approx&\,(0.1846862,1.5984702)\,;\nonumber\\
(j_2^{(7)},j_3^{(7)},j_4^{(7)})\approx &\,(2.471228,0.51685555,0.2137724)\,;\nonumber\\
(j_2^{(9)},\ldots,j_5^{(9)})\approx&\,
(2.396976693,2.079251103,-0.2079168419,-2.508673987)\,;\nonumber\\
(j_2^{(9)\prime},\ldots,j_5^{(9)\prime})\approx&\,
(-2.364737070,1.031057162,1.569692175,0.3383837765)\,.\nonumber\end{align}
$j^{(9)}$ corresponds to $\cD^0\Hg_9$ and $j^{(9)\prime}$ to $\cD^{0,\{5,17\}}\Hg_9$.

It should be possible to identify $\cA(\bbZ_\nu)$ for the next couple
$\nu$'s, using the method of $\nu=9$ described in the proof. To probe the answer, we
used Maple to obtain numerical solutions to \eqref{7.1}-\eqref{7.5}. 
We don't have any proof that these correspond to actual solutions (although the numerics
are convincing), and we certainly have no proof that there are no other solutions
for these $\nu$. 
The inequivalent numerical solutions we obtained for $K=\bbZ_{\nu}$, $11\le \nu\le 19$, are: 
\begin{align}
j^{(11)}\approx &\,(0.9996507,2.7258434,-0.5714203,-1.7797340,
1.2675985)\,,\nonumber\\
j^{(11)\prime}\approx&\,(-2.6444397,-1.7629598,-2.6444440,
2.7572657,0.1128260)\,;\nonumber\\
j^{(13)}\approx &\,( -3.1050384,0.5993399,-0.111708,
-0.969766,1.336848,1.00483129)\,;\nonumber\\
j^{(15)}\approx&\,(-1.0777623,-.7748018,-2.171863,-1.6068402,-.257508,2.092502,.72289565)\,;\nonumber\\
j^{(17)}\approx&\,(-1.466074,.291489,3.130735,-2.693185,1.398153,-.611938,-1.667078,-1.754821)\,;\nonumber\\
j^{(19)}\approx&\,(-2.677465,1.088972,-.899442,.015448,-1.240928,-.493394,1.839879,\nonumber\\
&\qquad -1.525884,-2.084374)\,;\nonumber\\
j^{(19)\prime}\approx&\,( .896858,-.882585,-2.369855,-1.873294,-1.711620,-.119360,2.972018,\nonumber\\
&\qquad-2.460652,.041334)\,. \nonumber\end{align}

For $\nu=13,15,17$ respectively, $\mu=173,229,293$ is prime. In these cases, choosing
$\omega_l=\xi_\mu^{ml^2}$ yields a unique solution to \eqref{8.10} (according to Maple)
and plugging into \eqref{8.12} yields a very close approximation (of order $10^{-6}$ or so)
$S'$ to the $S$ matrix of $\cD^0\Hg_\nu$. With high confidence we can expect $j^{(\nu)}$
for $\nu=13,15,17$ to describe a subfactor of Izumi type $\bbZ_\nu$ realising the modular
data $\cD^0\Hg_\nu$.

For $\nu=11,19$ respectively, $\mu=5^3,5\cdot 73$ are composite and the situation is
more subtle. We find (up to the numerical accuracy of Maple) that for these $j^{(\nu)}$,
$\omega_l=\xi_\mu^{ml^2}$ in \eqref{8.10} yields the correct multiplicity of solutions, and thus
the corresponding $T$ matrix agrees with that of $\cD^0\Hg_\nu$. Likewise, the solution
$j^{(19)\prime}$ corresponds to $\cD^{0,\{5,73\}}\Hg_{19}$. As with $\nu=9$, the
existence of  higher multiplicities means we cannot determine $S$ unambiguously from
\eqref{8.12} for these 3 $A$s. We expect $j^{(11)},j^{(19)},j^{(19)\prime}$ to correspond
to actual solutions of \eqref{7.1}-\eqref{7.5}, hence subfactors of Izumi type $\bbZ_{11},\bbZ_{19},
\bbZ_{19}$ respectively, with modular data $\cD^0\Hg_{11},\cD^0\Hg_{19},\cD^{0,\{5,73\}}\Hg_{19}$
resp.

The only unexpected solution here is $j^{(11)\prime}$. Following the method described in
the proof of Theorem 5, we obtain the corresponding $T$ matrix: for each $1\le l\le m$,
\begin{equation}
T_{\mathfrak{d}_l,\mathfrak{d}_l}=\left\{\begin{array}{cc}\xi_{\mu}^{ml^2}&\ \mathrm{if}\ 5|l\cr
\xi_{\mu}^{5ml^2}&\ \mathrm{otherwise}\end{array}\right.\,.\end{equation}
As before, $S$ is not uniquely determined by \eqref{8.12}. We do not have a guess for
what this $S$ matrix is, though its entries $S_{\mathfrak{d}_l,\mathfrak{d}_{l'}}$
should lie in $\bbQ[\cos(2\pi/25)]$. Just as extra solutions like $j^{(9)\prime}
,j^{(19)\prime}$ occur whenever \textit{distinct} primes divide $\mu$, we would expect
that solutions like $j^{(11)\prime}$ occur whenever a nontrivial prime power 
divides $\mu$.

On the basis of these observations, it is tempting to make the following conjecture.

\medskip\noindent\textbf{Conjecture 1.} \textit{The modular data $\cD^{0,P}\Hg_\nu$ for any  $\nu=2n+1$ and any set $P$ (possibly empty) consisting of
an even number of distinct prime divisors of $\mu$, is
realised by the even part of the quantum double of a subfactor of Izumi type
 $\bbZ_\nu$. All $A\in \cA(K)$ should correspond to modular
data obeying $T^{2\mu\nu}=I$. }\medskip

We expect $T_{\mathfrak{d}_l,\mathfrak{d}_l}^{2\mu}=1$ for reasons clear from the proof of
Theorem 5. Existence is accumulating suggesting that (up to equivalence) there is a unique
matrix $A\in\cA(\bbZ_\nu)$, and a unique subfactor of Izumi type $\bbZ_\nu$, iff
$\mu=\nu^2+4$ is prime. It would be interesting to see if  $\cA(K)$ can be
nonempty when $K$ is not cyclic.

\medskip\noindent\textbf{Question 5.} \textit{Interpret the $\bbZ_\nu$ twist
in $\cD^\omega\Hg_\nu$  in a von Neumann algebra
formalism.}

\medskip We would expect this to be in terms of twisted systems of endomorphisms, or
twisted fusion categories, as what happens for twisted doubles of
finite groups.

\subsection{Principal graphs, $\alpha$-induction etc for Izumi's subfactors}

Consider  a subfactor $N\subset M$ of Izumi type $K$, an abelian group of
odd order $\nu$, defined by some matrix $A\in\cA(K)$,  with $N$-$N$ and $M$-$M$ systems $\Delta$ and $\widehat{\Delta}$
respectively. We argued last subsection that the modular data of
 the double $\cD\Delta$ and $\cD\widehat{\Delta}$ is often (always?)
the untwisted Haagerup--Izumi modular data $\cD^0\Hg_\nu$ or some 
generalisation thereof. This subsection describes the data associated to this subfactor, generalising 
Section 2.2.  
{For the special case of the Haagerup subfactor, i.e. $K=\bbZ_3$, much of this
was computed in \cite{Bisch}, though without appreciating the underlying
$\bbZ_3$ structure which organises everything.}

Write $\kappa$ for the inclusion $N\subset M$. Write $K'$ for a set of representatives of
the equivalence classes $(K\setminus 0)/\pm$; then $g\in K'$
labels the conjugacy classes $|g|:=\{g,-g\}$ in the dihedral group
$K\sdprod\bbZ_2$ (where $\bbZ_2$ acts by inverse). 
Recall $\delta=(\nu+\sqrt{\mu})/2$,
$n=(\nu-1)/2$, $m=(\nu^2+3)/2$, $\mu=\nu^2+4=2m+1$ and write $\lambda=\sqrt{1+\delta}$.

\medskip\epsfysize=2.5in\centerline{ \epsffile{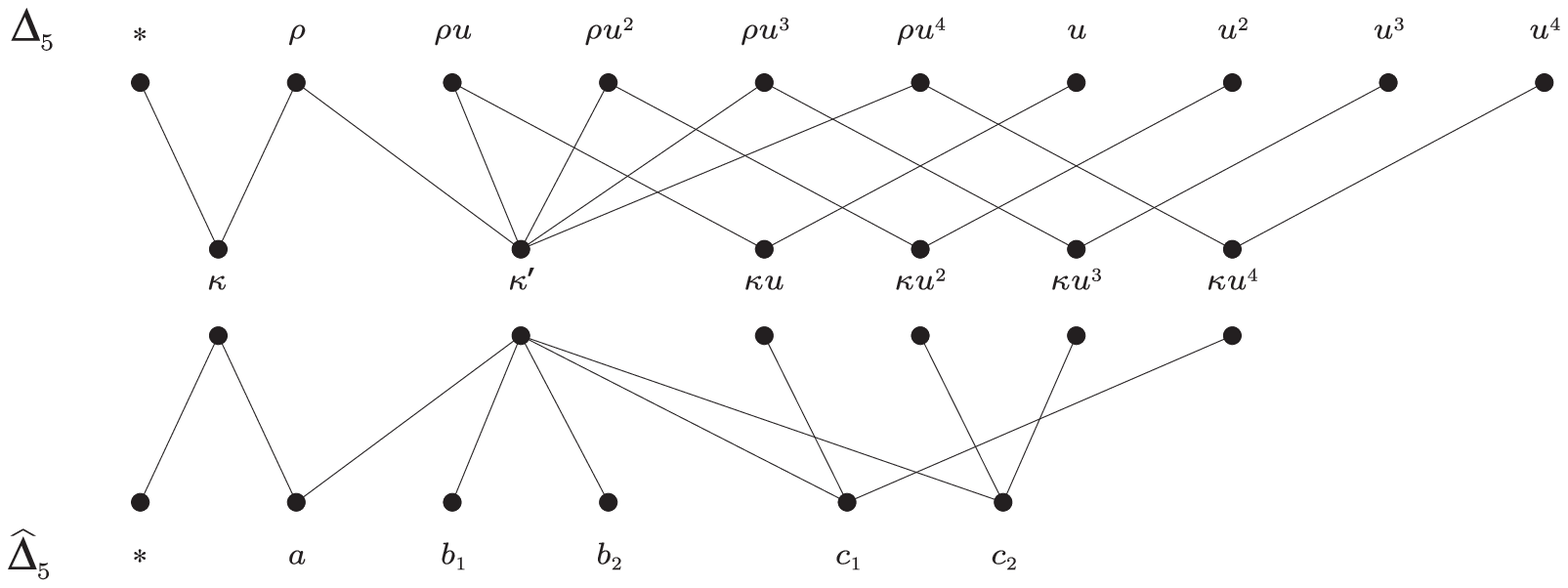}}\medskip

\centerline{Figure 4. Principal graphs for Izumi's $\bbZ_5$ subfactor}
\medskip

\medskip\noindent\textbf{Theorem 6.} \textit{Let $N\subset M$ be a subfactor of Izumi type
 $K$ as above. It has index $\delta+1$. The $N$-$N$ system 
is $\Delta=\{\mathrm{id},u,\ldots,{u_{\nu-1}},\rho,\rho u,\ldots,\rho{u_{\nu-1}}\}$, and its 
sectors  have fusions generated by
\begin{equation}\label{fusgen}[u_g][u_h]=[u_{g+h}]\,,\ [u_g][\rho]=[\rho][u_{-g}]\,,\ [\rho]^2=
[\mathrm{id}]+\sum[\rho u_g]\,.\end{equation}
The $M$-$N$ system consists of endomorphisms $\kappa u_g$ ($g\in K$) and 
$\kappa'$.
The principal graph consists of $\nu$ segments, each of length 3, sharing a common central 
vertex. The product $M$-$N\times N$-$N\rightarrow M$-$N$ is
\begin{eqnarray}&&[{\kappa u_g}][\rho u_h]=[\kappa']+[\kappa u_{-g+h}]\,,\ 
[\kappa'][{\rho u_g}]=\sum[\kappa u_h]+(\nu-1)[\kappa']\,,\nonumber\\
&&[\kappa u_g][u_h]=[\kappa u_{g+h}]\,,\ [\kappa'][{u_g}]=[\kappa']\,.\end{eqnarray}
The product $N$-$M\times M$-$N\rightarrow N$-$N$ is:}
\begin{eqnarray}
&&[\overline{\kappa u_g}][\kappa u_h]=[u_{h-g}]+[\rho u_{h+g}]\,,\
[\overline{\kappa'}][\kappa u_g]=\sum[\rho u_h]\,,\nonumber\\
&&[\overline{\kappa'}][\kappa']=\sum[u_g]+(\nu-1)\sum[\rho u_g]\,.\end{eqnarray}
\textit{The primaries of the double $\cD\Delta$ are $0,\mathfrak{b},\mathfrak{a}_i,\mathfrak{c}_{i,j},
\mathfrak{d}_l$ for $i\in K'$, $j\in K$, $1\le l\le m$. 
The canonical endomorphism  is 
$\theta_{\cD\Delta}=1+\mathfrak{b}+2\sum_i\mathfrak{a}_i$ with modular invariant 
$Z_\Delta=|ch_0+ch_{\mathfrak{b}}+2\sum_ich_{\mathfrak{a}_i}|^2$. The alpha-inductions are
$\alpha^+_x=\alpha x\otimes 1$, $\alpha^-_x=1\otimes (\alpha x)^{\mathrm{opp}}$ 
where, for all $i\in K',j\in K,1\le l\le m$,
\begin{align}&&[{\alpha}\mathfrak{a}_i]=
2[\mathrm{id}]+\sum[\rho u_j]\,,\ [{\alpha}\mathfrak{b}]=
[\mathrm{id}]+\sum[\rho u_j]\,,\nonumber\\ 
&&[{\alpha}\mathfrak{c}_{i,j}]=[u_i]+[u_{-i}]+\sum[\rho u_{j'}]\,,\
[{\alpha}\mathfrak{d}_l]=\sum[\rho u_j]\,.\end{align}
}\medskip

\noindent\textit{Proof.} In Lemma 7.1 and Theorem 7.2 of \cite{iz3}, we learn that $\Delta$
contains $u_g$ and $\rho u_g$ for all $g\in K$, and the fusions
$[\overline{\kappa}\kappa]=[\mathrm{id}]+[\rho]$, $[u_g][u_h]=[u_{g+h}]$ and $[\rho]^2= [\mathrm{id}]+\sum
[\rho u_g]$.  By definition $\Delta$ consists of all irreducible sectors
in any $(\overline{\kappa}\kappa)^k$, and from this we quickly obtain that $\Delta=\{u_g,\rho u_g
\}$. Because $\langle \kappa\rho,\kappa\rho\rangle=\langle\overline{\kappa}\kappa\rho,\rho
\rangle=\langle(\rho+1)\rho,\rho\rangle=2$ and $\langle\kappa\rho,\kappa u_g\rangle=
\langle(\rho+1)\rho,u_g\rangle=\delta_{g,0}$, we see that $[\kappa\rho]-[\kappa]:=[\kappa']$
is an irreducible $M$-$N$ sector distinct from any $[\kappa u_g]$. We obtain the fusion
$[\kappa' u_g]=[\kappa']$ from the calculation $\langle\kappa' u_g,\kappa'\rangle=
\langle\kappa(\rho-1)u_g,\kappa(\rho-1)\rangle=\langle(\rho^2-1)u_g,\rho-1\rangle
=1$. $[\overline{\kappa}][\kappa']$ has dimension $\lambda^2(\delta-1)=\nu\delta$ so will
equal the sum of $\nu$ sectors $[\rho u_g]$, perhaps with multiplicities, but $K$-invariance
now forces $[\overline{\kappa}][\kappa']=\sum[\rho u_g]$. As above we see that the
only sectors in $\kappa(\overline{\kappa}\kappa)^k$ are $[\rho u_g],[\kappa']$, so we
have exhausted all $M$-$N$ sectors, and obtain the given principal graph.
All remaining products are now easy to compute. For example, $\overline{\kappa'}\kappa'$
has dimension $\lambda^2(\delta-1)^2=\nu+(\nu-1)\nu\delta$, so by $K$-invariance this must
be $\sum[u_g]+(\nu-1)\sum[\rho u_g]$.

At the end of Section 4 of \cite{iz1} is explained how to recover the induction-restriction
(i.e. dual principal) graph between $A$-$B$ and $A$-$A$ sectors, from the structure of
Tube$(\Delta)$. In particular, we obtain the dual principal graph directly from the half-braidings
as described in Proposition 8.2 of \cite{iz3} --- for instance $\theta$ is read off from the 
half-braidings containing the identity. From the graph we read off the given alpha-inductions.
QED\medskip

From this the dual principal
graph for the double $\cD\Delta$ is the obvious analogue of the top graph of Figure 3.
Note the resemblance between the fusions for $\Delta$ and the group algebra $\bbC D_K$
of $D_K=K\sdprod\bbZ_2$, where $\bbZ_2$ acts on $K$ by inverse. We also see 
that the sectors $[u_g]$, $[\rho u_g]$, $[\kappa u_g]$, and $[\kappa']$ have 
statistical dimensions 
1,  $\delta$, $\lambda$ and $(\delta-1)\lambda$ respectively. It is far from obvious (though
of course an immediate  consequence of Theorem 6) that the map $\alpha$ defined there
is a ring homomorphism.

\medskip\noindent\textbf{Theorem 7.} \textit{Assume the hypotheses and notation
of Theorem 6. In addition, assume the fusions of the double system $\cD\Delta=
\{\mathrm{id},\mathfrak{b},\mathfrak{a}_i,\mathfrak{c}_{i,j},\mathfrak{d}_l\}$
agree with those of $\cD^0\Hg_\nu$, and that we have the canonical endomorphism
$\theta_{\cD\widehat{\Delta}}=1+\mathfrak{b}+\sum\mathfrak{a}_i
+\sum\mathfrak{c}_{h,0}$. Then we can identify 
the $M$-$M$ system $\widehat{\Delta}$ with $\{\mathrm{id},a,b_i,c_i\}_{i\in K'}$
where $c_i={\kappa}
u_i\overline{\kappa}\cong{\kappa}u_{-i}\overline{\kappa}$. 
Its sectors obey the fusions  
\begin{align}&&[a]^2=[\mathrm{id}]+[a]+\sum[b_i]+\sum[c_i]\,,\quad 
 [a][b_i]=[a]+\sum[b_k]-[b_i]+\sum[c_k]\,,\nonumber\\ &&[a][c_i]=[a]+ \sum[b_k]+ 
[c_i]+\sum[c_k]\,,\quad
[b_i][c_j]=[a]+\sum[b_k]+\sum[c_k]\,,\nonumber\\
&&[b_i][b_j]=
\delta_{i,j}[\mathrm{id}]+(1-\delta_{i,j})[a]+\sum[b_k]-[b_{i+j}]-[b_{i-j}]+\sum[c_k]\,,\nonumber\\ &&
[c_{i}][c_{j}]=\delta_{i,j}[\mathrm{id}] +
(1+\delta_{i,j})[a]+\sum[b_k]+ [c_{i+j}]+[c_{i-j}]+\sum[c_k]\,,\nonumber
\end{align}
where we set $[c_0]=0=[b_0]$ and write $[c_{-i}]=[c_i]$, $[b_{-i}]=[b_i]$.
The dual principal graph is as in Figure 4: e.g. there is a (odd) 
central vertex, to which is attached $[a]$, the $[b_i]$'s and
the $[c_i]$'s. The product $M$-$M\times M$-$N\rightarrow M$-$N$ is
\begin{align} 
&[a][\kappa u_g]=[\kappa u_g]+[\kappa']\,,\ [a][\kappa']=\sum[\kappa u_g]+
(\nu-1)[\kappa']\,,\ [b_{i}][\kappa u_g]=[\kappa']\,,\nonumber&\\
&[b_{i}][\kappa']=(\nu-2)[\kappa']+\sum[\kappa u_g]\,,\ 
[c_{i}][\kappa u_g]=[u_{g+i}]+[\kappa u_{g-i}]+[\kappa']\,,\nonumber&\\
&[c_{i}][\kappa']=\sum[\kappa u_g]+\nu[\kappa']\,,\end{align}
while the product $M$-$N\times N$-$M\rightarrow M$-$M$ is
\begin{align}
&[{\kappa u_g}][\overline{\kappa u_h}]=([\mathrm{id}]+[a])\delta_{g,h}+[c_{g-h}](1-\delta_{g,h})\,,\
[\kappa u_g][\overline{\kappa'}]=[a]+\sum[b_{i}]+\sum[c_{i}]\,,&\nonumber\\
&[{\kappa'}][\overline{\kappa'}]=[\mathrm{id}]+(\nu-1)[a]+(\nu-2)\sum[b_{i}]+
\nu\sum[c_{i}]\,.&\end{align}
The dual canonical modular invariant is $\widehat{Z}=|ch_0+ch_\mathfrak{b}+\sum ch_{\mathfrak{a}_i}
+\sum ch_{\mathfrak{c}_{h,0}}|^2$. Writing $\widehat{\alpha}_x^+=\widehat{\alpha}
x\otimes 1$ and $\widehat{\alpha}_x^-=1\otimes(\widehat{\alpha}x)^{\mathrm{opp}}$,
the alpha-inductions are}
\begin{align} &[\widehat{\alpha}\mathfrak{a}_i]=[\mathrm{id}] +2[a]+
\sum[b_h]-[b_i]+\sum[c_h]\,,\ [\widehat{\alpha}\mathfrak{b}]=[\mathrm{id}]+
[a]+\sum
[b_h]+\sum[c_h]\,,\nonumber&\\ &[\widehat{\alpha}\mathfrak{c}_{i,j}]=[\mathrm{id}]\delta_{j,0}+
([a]-[b_i])(1-\delta_{j,0})+\sum [b_h]+[c_i]+\sum[c_h]\,,\nonumber&\\
&[\widehat{\alpha}\mathfrak{d}_l]=[a]+\sum[b_i]+\sum[c_i]\,.&\end{align}

\noindent\textit{Proof.} From $\widehat{\theta}$ we obtain immediately $\widehat{Z}$.
The cardinality of $\widehat{\Delta}$ is determined from $\langle\widehat{\gamma},
\widehat{\gamma}\rangle=\langle\widehat{\theta},\widehat{\theta}\rangle=\sum \widehat{Z}_{i,0}^2
=\nu+1$. Fix any $l$ and define $[\mathrm{id}]=[\widehat{\alpha}_0]$ (the
identity sector of $\widehat{\Delta}$), $A=[\widehat{\alpha}(\mathfrak{d}_l-\sum \mathfrak{c}_{h,0}
+\sum 
\mathfrak{a}_i)]$, $B_i=A+\nu\,[\widehat{\alpha}(\mathfrak{b}-\mathfrak{a}_i)]$, 
$C_h=A+\nu\,[\widehat{\alpha}(\widehat{c}_{h,0}-\widehat{b})]$. We compute
(using  $ \langle\widehat{\alpha}{x},\widehat{\alpha}y\rangle=\langle\widehat{\theta},xy\rangle$
and the fusions of the double computed in Section 3.2) that 
$\langle \widehat{\alpha}(\mathfrak{a}_i-\mathfrak{b}),\widehat{\alpha}(\mathfrak{a}_{i'}-\mathfrak{b})\rangle=1+\delta_{i,i'}$,
$\langle \mathfrak{c}_h-\mathfrak{b},\mathfrak{c}_{h'}-\mathfrak{b}\rangle=1+\delta_{h,h'}$,
$\langle \mathfrak{a}_i-\mathfrak{b},\mathfrak{c}_{h}-\mathfrak{b}\rangle=-1$,
$\langle \mathfrak{d}_l,\mathfrak{a}_i-\mathfrak{b}\rangle=\langle\mathfrak{d}_l,
\mathfrak{c}_{h}-\mathfrak{b}\rangle=0$, $\langle \mathfrak{d}_l,\mathfrak{d}_{l}\rangle=\nu$.
Hence the $A$, $B_i$, $C_h$ are mutually orthogonal,
each with norm $\nu$, so they form a basis (over $\bbQ$) for the vector space spanned
by the sectors of $\widehat{\Delta}$. 

Let $[e_k]_{0\le k\le \nu}$ be the irreducible sectors of $\widehat{\Delta}$; choose 
$[e_0]=[\mathrm{id}]$. Then $[A],[B_i],[C_i]$
are linear combinations (over $\bbZ$) of these $[e_k]$: define $\vec{A},\vec{L}^{(1)},\ldots,
\vec{L}^{(2n)}\in\bbZ^{\nu+1}$ by $[A]=\sum \vec{A}_k[e_k]$, $[B_i]=[A]+\nu\sum \vec{L}_k{}^{(i)}
[e_k]$, $[C_i]=[A]+\nu\sum\vec{L}_k{}^{(i+n)}[e_k]$. Then $\vec{A}\cdot\vec{A}=\nu^2$,
$\vec{A}\cdot\vec{L}{}^{(k)}=-\nu$, $\vec{L}{}^{(k)}\cdot\vec{L}{}^{(l)}=1+\delta_{k,l}$.
From $\langle[\mathrm{id}],\widehat{\alpha}_{\mathfrak{a}_{i}}-\widehat{\alpha}_{\mathfrak{b}}
\rangle=\langle[\mathrm{id}],\widehat{\alpha}_{\mathfrak{c}_{h,0}}-
\widehat{\alpha}_{\mathfrak{b}}\rangle=\langle[\mathrm{id}],\widehat{\alpha}_{\mathfrak{d}_{l}}
\rangle=0$ we have $\vec{A}_0=\vec{L}{}^{(k)}_0=0$. 
Thus (reordering the $[e_k]$ if necessary) 
$\vec{A}=(0,s_1\,(\nu+t),\ldots,s_\nu\,(\nu+t),st)$ and
$\vec{L}^{(k)}$ has components $\vec{L}^{(k)}_i=s\delta_{i,\nu}+s_k\delta_{i,k}$,
where $s,s_k\in\{\pm 1\}$  and $t\in\bbZ$. This integer $t$ 
satisfies $\nu^2=t^2+\nu\,(\nu+t)^2$, so $t=-\nu$.

This means that $[a]:=\frac{1}{\nu}[A]$, and hence $[b_i]:=\frac{1}{\nu}[B_i]$ and $[c_i]:=
\frac{1}{\nu}[C_i]$, all being norm 1, are in $\widehat{\Delta}$, possibly up to signs. 
That the signs are indeed all +1's follows by writing say $[\widehat{\alpha}\mathfrak{d}_1]$ 
in terms of $[a],[b_i],[c_i]$: the coefficients should be nonnegative.

The fusions in $\widehat{\Delta}$ now reduce to the fusions in the double. The dual
principal graph and the two remaining products follow by arguments as in
Theorem 6. QED

\medskip
This guess for $\widehat{\theta}$ is the natural one, matching what happens for the
Haagerup subfactor (recall Section 2.2) and what we know about monomial
modular invariants in Proposition 4, which appear to be severely constrained.
The only example we know violating the hypothesis on the fusions of $\cD\Delta$
is the solution $j^{(11)\prime}$ of Section 4.1. 
The noncommutativity of $\Delta$ and commutativity of $\widehat{\Delta}$ arise because
of the presence respectively absence of higher multiplicities in $\theta$ resp. $\widehat{\theta}$. 
The generalisation of Figure 3  can be read off from the alpha-inductions
listed in the Theorem.
These imply the statistical dimensions  dim($a)=\delta$, dim($b_i)=\delta-1$ and dim($c_j)=\delta+1$.

The dihedral group ${K\sdprod \bbZ_2}$ has one nontrivial 1-dimensional irrep and $n$ 2-dimensional
irreps. But for the Haagerup considered in Section 2.2, $a$ corresponded to the  2-dimensional
and $b$ the 1-dimensional irreps, so $\widehat{K\sdprod\bbZ_2}$ and 
$\widehat{\Delta}$ here are
unrelated for $\nu>3$, even after projecting away the $c_i$'s.

%A simple consistency check is  dimension: e.g. the quantum-dimension of
%$\mathfrak{a}_i$ should equal the sum of the statistical dimensions of 1, $2a$, and $nc_j$.
%We have verified these are all satisfied. Alpha-induction is an algebra homomorphism; a special
%case of this is
%\begin{equation} \langle \alpha {x},\alpha  {y}\rangle=
%\langle {\theta},{xy}\rangle\,,\ \langle\widehat{\alpha}{x}
%,\widehat{\alpha}y\rangle=\langle\widehat{\theta},xy\rangle\,, \end{equation}
%valid for all $\cD\Delta$ primaries $x,y$ (the fusions of $\cD^\omega\Hg_\nu$ are
%given in Section 3.2).

%Then the multiplicity of $[c_k]$ in $[c_i][c_j]$ equals
%\begin{align}&\langle\overline{\kappa}u^i\kappa\,\overline{\kappa}u^j\kappa,\overline{\kappa}u^k
%\kappa\rangle=\langle u^i(1+\rho)u^j,(1+\rho)u^k(1+\rho)\rangle\nonumber&\\ 
%&=\langle u^{i+j}+\rho u^{-i+j}\rho,
%u^k+\rho u^k+\rho u^{-k}+u^{-k}+\sum\rho u^{-h}\rangle=1+\delta_{i+j,k}
%+\delta_{i-j,k}+\delta_{i-j,-k}+\delta_{i-j,-k}\,.\nonumber&\end{align}

\section{VOAs for Haagerup-Izumi modular data}

 \subsection{The Haagerup-dihedral diamond}

In Section 2.4 we saw $\cD \Hg$
as a mutation of $\cD S_3$ and contained in a conformal subalgebra of $\cD \bbZ_3$;
now this becomes $\cD^\omega \Hg_\nu$ is a mutation of $\cD^{+,\omega} D_\nu$ and is contained
in a conformal subalgebra of $\cD^\omega \bbZ_\nu$. The branching rules for this embedding
can be read off from the modular invariant \eqref{modinvDnuZnu}.
In some ways (e.g. modular data) the Haagerup-Izumi
closely resembles its sibling the dihedral group, and in other ways (e.g.  VOA realisations) it more closely
relates to its parent the cyclic group.

The trivial case $\cD^0\Hg_1$ is already interesting.
%$T={\rm diag}(1,1,\xi_5^2,\xi_5^3)$ and 
%$$S=\frac{1}{2\sqrt{5}}\left(\begin{matrix}
%\sqrt{5}-1&\sqrt{5}+1&2&2\cr \sqrt{5}+1&\sqrt{5}-1&-2&-2\cr 2&-2&1-\sqrt{5}&1+\sqrt{5}\cr
%2&-2&1+\sqrt{5}&1-\sqrt{5}\end{matrix}\right)$$
In $c=8$, its VOA can be realised as the tensor of the affine algebra VOA $\cV(G_{2,1})$ (which has
$c=14/5$), with the affine algebra VOA $\cV(F_{4,1})$ (which has $c=26/5$). (Hence
there is a VOA realisation for $c$ any positive multiple of 8, by tensoring
by copies of $\cV(E_8)$.) Call this $\cV(\Hg_1^0)$.
The VOA corresponding to $\cD\bbZ_1$ in $c=8$ of course must be $\cV(E_8)=
\cV(E_{8,1})$, and the containment corresponds to the conformal embedding
$G_{2,1}F_{4,1}\subset E_{8,1}$. The VOA corresponding to $\cV D_1$ is
the lattice VOA $\cV(A_1 E_7)$. Although the latter is an orbifold of $\cV(E_8)$,
$\cV(G_{2,1}F_{4,1})$ cannot be, since its modular data is not that of
the (possibly twisted) quantum double of a finite group. This is an important lesson:
we cannot expect the $\cD^\omega \Hg_\nu$ VOAs to be orbifolds,
by some subgroup of automorphisms, of $\cD^\omega\bbZ_\nu$ VOAs.
Hence Question 3.

Another important observation from this trivial case $\nu=1$ is that the intersection of
the holomorphic orbifold $\cV D_1$ with the $\nu=1$ Haagerup-Izumi VOA $\cV(\Hg_1^0)$
is itself a rational VOA. In particular, there are conformal embeddings $\cV(B_{3,1}
G_{2,1})\subset\cV( E_{7,1})$ and $\cV(A_{1,1}B_{3,1})\subset\cV( F_{4,1})$ \cite{SchWa}, and so
both $\cV D_1=\cV(A_1 E_7)$ and $\cV(\Hg_1^0)=\cV(G_{2,1}F_{4,1})$ contain  copies
of the rational VOA $\cV(A_{1,1} G_{2,1}B_{3,1})$. Thus in this baby example at
least, the dihedral and  Haagerup-Izumi VOAs are commensurable in this sense. 
It is tempting to guess though that something similar happens for higher $\nu$. This is what
we mean by the Haagerup-dihedral diamond: at the top is a VOA realising $\cD^\omega\bbZ_\nu$,
which contains as conformal subalgebras $\cD^{+,\omega}D_\nu$ and $\cD^\omega\Hg$, and these
contain a common conformal subalgebra (the analogue of $\cV(A_{1,1} G_{2,1}B_{3,1})$).

\medskip\noindent\textbf{Conjecture 2.} \textit{There is a rational VOA realisation of each
$\cD^\omega \Hg_\nu$ for any sufficiently large central charge $c$ ($c$
a multiple of 8).}\medskip

As $\nu$ increases, we will not be able to realise all $\cD^\omega \Hg_\nu$ --- or even
the $\cD^\omega \bbZ_\nu$ --- as conformal subalgebras of 
$\cV(E_{8})$. In other words, $c>8$ will sometimes be necessary.
This conjecture is open even for $\nu=3$. To make it more plausible, and also to aid in the
construction of the VOAs, in the following sections we supply all the information needed
to identify the possible character vectors. We saw in Section 2.4 that the Haagerup VOA
would see the holomorphic orbifold by $\bbZ_3$ more directly than that by $S_3$, so for
this reason we also supply the relevant data for $\cD^\omega\bbZ_3$.

We should repeat here the suggestion of Section 2.5 that a second approach to constructing 
VOAs realising $\cD^\omega\Hg_\nu$ is to use the affine algebra $\cV(B_{m,2})$ in the
coset construction. For reasons of space, we haven't explored this in this paper.

\subsection{Character vectors for the untwisted  Haagerup double}

In this section we determine the matrices $\Lambda$ and $\chi$ for the untwisted Haagerup
modular data $\cD^0\Hg_3$, at all allowed values of central charge $c$ (i.e. when $8|c$).
As explained in Section 1.4, from this all possible character vectors can be obtained. We 
illustrate this with $c=16$ and $c=24$ (for $c=8$ see Section 2.4).

\subsubsection{Untwisted Haagerup at central charge $c\equiv_{24}16$}

We found $\Lambda,\chi$ for $\cD^0\Hg_3$ at $c\equiv_{24}8$ in Section 2.4. 
Changing the central charge by a multiple of 8 leaves $S$ unchanged but multiplies $T$ by
some 3rd root $\omega$ of 1. 
At the end of Section 1.4 we explain how this affects $\Lambda,\chi$.
Taking the primaries in the usual order $0,\mathfrak{b};\mathfrak{a};\mathfrak{c}_0,
\mathfrak{c}_1,\mathfrak{c}_2;\mathfrak{d}_1,\ldots,\mathfrak{d}_6$, we obtain
for $c\equiv_{24}16$
\begin{align}&&\Lambda={\rm diag}
(1/3,1/3;1/3;1/3,2/3,0;31/39,7/39,19/39,28/39,-5/39,-2/39)\,,\nonumber\\
&&\chi ={\scriptsize \left(\begin{array}{cccccccccccc}
17& 155& 162& 162& 27& 729& 13& 286&65& 13&  1001& 728\cr 
155& 17& 162& 162& 27& 729& -13& -286&-65& -13&  -1001& -728\cr 
162& 162& 334& -162& -27& -729& 0& 0& 0& 0& 0& 0\cr 
162& 162& -162& 334& -27& -729& 0& 0& 0& 0& 0& 0\cr 
1215& 1215& -1215& -1215& -76& 17496& 0& 0& 0& 0& 0& 0\cr 
9& 9& -9& -9& 6& -12& 0& 0& 0& 0& 0& 0\cr 
1925& -1925& 0& 0& 0& 0& -21& -6250&-125& 50&  51331& 43175\cr 
45& -45& 0& 0& 0& 0& -9& 58&43& 12&  -45& -297\cr 
374& -374& 0& 0& 0& 0& 0& 1650& 106& -51& -4250& 3927\cr 
1288& -1288& 0& 0& 0& 0& 23& 5313&-759& 4&  39721&-15432\cr 
1& -1& 0& 0& 0& 0& 1& 0&-1& 1&  3& 1\cr
3& -3& 0& 0& 0& 0& 3& -9& 4& -2& 6& 4\end{array}\right)}\,.\nonumber\end{align}
 From \eqref{recurs} the matrix $\Xi(\tau)$ is recursively computed. Any weakly
holomorphic vector-valued modular function equals $\Xi\vec{P}(J)$ for some polynomial
vector $\vec{P}$.

Specialise now to a VOA character vector at $c=16$. The corresponding polynomials $P_i(x)$ 
are heavily constrained. For one thing, the vacuum character 
$ch_0(\tau)$ starts like 
$1q^{-2/3}+\cdots$ and all other characters $ch_i(\tau)$ start like $a_iq^{h_i-c/24}+
\cdots$ where $h_i>0$. This forces $\vec{P}=(1,0,0,0,e,0,g,h,0,j,0,0)^T$ for \textit{constants}
$e,g,h,j$.
So we've uniquely determined the hypothetical $c=16$ Haagerup character vector,
once 4 numbers (namely $e,g,h,j$) are
specified. Its first few terms are
\begin{equation}
{\scriptsize\left(\begin{matrix}ch_0(\tau)\cr ch_{\mathfrak{b}}(\tau)\cr ch_{\mathfrak{a}}(\tau)=
ch_{\mathfrak{c}_0}(\tau)\cr
ch_{\mathfrak{c}_1}(\tau)\cr ch_{\mathfrak{c}_2}(\tau)\cr ch_{\mathfrak{d}_1}(\tau)\cr ch_{\mathfrak{d}_2}(\tau)\cr ch_{\mathfrak{d}_3}(\tau)\cr
ch_{\mathfrak{d}_4}(\tau)\cr ch_{\mathfrak{d}_5}(\tau)\cr ch_{\mathfrak{d}_6}(\tau)\end{matrix}\right)}=
{\scriptsize\left(\begin{matrix}q^{-2/3}(1+(17+13g+27e+65h+13j)q+(2013+594e+4173h+260j
+91g)q^2+\cdots)\cr%\qquad\qquad\qquad+(59090+5967e+79742h+1989j+468g)q^2+\cdots)\cr 
q^{1/3}((155-13g+27e-65h-13j)+(21245+594e-4173h-260j-91g)q+\cdots)\cr
%\qquad\qquad\qquad+(645870+5967e-79742h-1989j-468g)q^2+\cdots)\cr 
q^{1/3}\left((-27e+162)+(23247-594e)q+(705024-5967e)q^2+\cdots\right)\cr 
%q^{1/3}\left((-27e+162)+(23247-594e)q+(705024-5967e)q^2+\cdots\right)\cr 
q^{-1/3}\left(e+(-76e+1215)q+(-1384e+79704)q^2+(-11580e+1886166)q^3+
\cdots\right)\cr 
(6e+9)+(5832+486e)q+(247131+5832e)q^2+\cdots\cr
q^{-8/39}(g+(-21g+1925-125h+50j)q+(-239g+103631-4376h+506j)q^2+\cdots)\cr
%\qquad\qquad\qquad+(-1230g+2257134-61606h+3310j)q^2+\cdots)\cr
q^{7/39}\left((43h+45+12j-9g)+(4792h+258j+10400-71g)q+%(364932+102652h+2167j-388g)q^2+
\cdots\right)\cr 
 q^{-20/39}\left(h+(106h+374-51j)q+(4980h-634j+34606-43g)q^2+
%(84284h-4688j+930900-172g)q^2+
\cdots\right)\cr 
q^{-11/39}(j+(4j+1288+23g-759h)q+(-84j+141g+79583-25135h)q^2+\cdots\cr
%\qquad\qquad\qquad+(-752j+668g+1816206-375572h)q^2+\cdots)\cr 
q^{-5/39}\left((-h+1+j+g)+(-856h+130j+2704+31g)q+
%(134516-25364h+1472j+228g)q^2+
\cdots\right)\cr q^{-2/39}\left((4h+3-2j+3g)+(1582h-95j+3898+70g)q+
%(173656+41630h-1036j+439g)q^2+
\cdots\right)\end{matrix}\right)}\,.\nonumber\end{equation}

Since the Fourier coefficients of a VOA character vector are nonnegative
integers, we know each $e,g,h,j\in\bbZ_{\ge 0}$.
In addition the $q^{1/3}$-coefficient of $ch_{\mathfrak{a}}$ forces $e\le 6$ and the
$q^{1/3}$-coefficient of $ch_{\mathfrak{b}}$ then forces $g+j+5h\le 24$, so there
are only finitely many possibilities. Other coefficients yield further bounds,
e.g. $h\le 2$, $j\le 11$, $g\le 16$.

There are several other constraints. One is that VOA characters are linear 
combinations over $\bbZ_{\ge 0}$
of the Virasoro characters  at that $c$ for various values of \textit{conformal weight}
$h$ (because the VOA carries a Virasoro representation). 
All these characters are well-known: $ch_{(16,0)}=q^{-2/3}(1+
q^2+q^3+2q^4+\cdots)$ and $ch_{(24,h)}=q^{h-2/3}(1+q+2q^2+3q^3+5q^4+\cdots)$ for 
$h>0$. This isn't so useful here, but it is effective for the $c=24$ case considered next.

More important,
Dong-Mason \cite{DM} prove that the conformal weight 1 part of a rational VOA
$\cV$ is a reductive Lie algebra $\oplus_i\mathfrak{g}_i$ of central charge 
$\sum_ik_i\mathrm{dim}\,\mathfrak{g}_i/(k_i+h^\vee_i)\le c$ and dimension equal to the coefficient of $q^{1-c/24}$ in
$ch_0$. This gives a finite set of possible values for that coefficient, and knowing this
reductive Lie algebra helps in constructing $\cV$. Each $\cV$-module is also a module for the
corresponding affine Lie algebra, so the $\cV$-characters are linear combinations of 
the affine characters. 

This is as far as we'll carry the analysis here.

The matrices $\Lambda,\chi$ for $\cD^0\bbZ_3$, with primaries in the order 
$(00),(01)=(02),(10)=(20),(11)=(22),(12)=(21)$ are
$\Lambda=\textrm{diag}(1/3,1/3,1/3,2/3,0)$ and 
\begin{equation}\chi= {\scriptsize\left(\begin{array}{ccccc} 172 & 324&324& 54&1458\\
162& 334 & -162 & -27 & -729 \\
162 & -162 & 334& -27& -729\\
1215& -1215 &-1215 & -76&17496\\
 9 & -9 & -9 & 6 & -12\end{array}\right)}\,.\nonumber\end{equation}
At $c=16$ there is only one free parameter, $d$: the possible character vectors are
\begin{equation}{\scriptsize\left(\begin{matrix}ch_{(0,0)}\cr ch_{(0,1)}=ch_{(0,2)}
=ch_{(1,0)}=ch_{(2,0)}\cr ch_{(1,1)}=ch_{(2,2)}\cr 
ch_{(1,2)}=ch_{(2,1)}\end{matrix}\right)}=
{\scriptsize\left(\begin{matrix}  q^{-2/3}(1+ (172 + 54 d)q + (23258 + 1188 d) q^2 +\cdots)\\
q^{1/3}((162 - 27 d) + (23247 - 594 d) q + (705024 - 5967 d) q^2+\cdots)\\
q^{-1/3}(d    + (1215 - 76 d)q + (79704 - 1384 d) q^2 + \cdots)\\
 (9 + 6 d) + (5832 + 486 d) q + (247131 + 5832 d) q^2+\cdots
\end{matrix}\right)}\,.\nonumber\end{equation}
This implies $d\in\bbZ_{\ge 0}$ and $d\le 6$. This $d$ equals the parameter 
$e$ of $\cD^0\Hg$.

\subsubsection{Untwisted Haagerup at central charge $c\equiv_{24}0$}

Next turn to the Haagerup at $24|c$. Finding its $\Lambda,\chi$ is now routine:
\begin{align}&&\Lambda = {\rm diag}
(0,0,1,1,1/3,2/3,6/13,11/13,2/13,5/13,7/13,8/13)\,,\nonumber\\
&&\chi= {\scriptsize\left(\begin{array}{cccccccccccc} 
-12& 0& 1/2& 1/2& 18& 6& 6& 1& 10& 5& 4& 4\cr
0& -12& 1/2& 1/2& 18& 6& -6& -1& -10& -5& -4& -4\cr  
65610& 65610& 0& 0& -5832& -243& 0& 0& 0& 0& 0& 0\cr
65610& 65610& 0& 0& -5832& -243& 0& 0& 0& 0& 0& 0\cr
729& 729& 0& 0& -152& 54& 0& 0& 0& 0& 0& 0\cr
8748& 8748& 0& 0& 2430& -76& 0& 0& 0& 0& 0& 0
\cr 1716& -1716& 0& 0& 0& 0& -252& -5& -176& 176& 134& 100
\cr 22451& -22451& 0& 0& 0& 0& -980& 8& 16556& 2464& -96& -385
\cr 104& -104& 0& 0& 0& 0& 0& 5& 44& -56& -15& 20
\cr 910& -910& 0& 0& 0& 0& 77& 8& -847& -32& 120& -34
\cr 3003& -3003& 0& 0& 0& 0& 330& 0& -1540& 704& -216& 55
\cr 5200& -5200& 0& 0& 0& 0& 616& -14& 3388& -595& 132& -20
\end{array}\right)}\,.\nonumber\end{align}
Its possible character vectors are found as before to be
$${\scriptsize\left(\begin{matrix}ch_0(\tau)\cr ch_{\mathfrak{b}}(\tau)\cr ch_{\mathfrak{a}}(\tau)\cr ch_{\mathfrak{c}_0}(\tau)\cr
ch_{\mathfrak{c}_1}(\tau)\cr ch_{\mathfrak{c}_2}(\tau)\cr ch_{\mathfrak{d}_1}(\tau)\cr ch_{\mathfrak{d}_2}(\tau)\cr ch_{\mathfrak{d}_3}(\tau)\cr
ch_{\mathfrak{d}_4}(\tau)\cr ch_{\mathfrak{d}_5}(\tau)\cr ch_{\mathfrak{d}_6}(\tau)\end{matrix}\right)}=
{\scriptsize\left(\begin{matrix}q^{-1}+(-12+10i+18e+{c+d\over 2}+6f+4k+5j+h+4l+6k)+
\cdots\cr (-10i+18e+{c+d\over 2}+6f-4k-5j-h-4l-6g)+\cdots\cr c+(65610-5832e-243f)q+
(7164612-247131e-2916f)q^2+\cdots\cr d+(65610-5832e-243f)q
+(7164612-247131e-2916f)q^2+\cdots
\cr eq^{-2/3}+(729-152e+54f)q^{1/3}+(370332-23236e+1188f)q^{4/3}+\cdots\cr 
fq^{-1/3}+(8748-76f+2430e)q^{2/3}+(1743039-1384f+159408e)q^{5/3}+\cdots
\cr g q^{-7/13}+(1716-252g -176i+134k+100l-5h+176j)q^{6/13}+\cdots
\cr hq^{-2/13}+(22451+8h+2464j+16556i-385l-96k-980g)q^{11/13}+\cdots
\cr iq^{-11/13}+(104+44i-15k+5h-56j+20l)q^{2/13}+\cdots
\cr jq^{-8/13}+(910-32j-34l-847i+77g+8h+120k)q^{5/13}+\cdots
\cr kq^{-6/13}+(3003+704j-1540i+55l-216k+330g)q^{7/13}+\cdots
\cr lq^{-5/13}+(5200-20l+3388i+616g-14h+132k-595j)q^{8/13}+\cdots\end{matrix}\right)}$$
for certain nonnegative integers $c,d,e,f,g,h,i,j,k,l$. By the usual arguments we find a
finite number of possibilities. 
Likewise the matrices for 
$\cD^0\bbZ_3$ are $\Lambda=\textrm{diag}(0,1,1,1/3,2/3)$ and 
$$\chi= {\scriptsize\left(\begin{array}{ccccc}
-12&1&1&36&12\\ 65610&0&0&-5832&-243\\ 65610&0&0&-5832&-243\\ 729&0&0&-152&
54\\ 8748&0&0&2430&-76\end{array}\right)}$$
and we find a character vector depending on 4 bounded nonnegative integers. 

The most intriguing possibility is if these Haagerup and $\bbZ_3$ VOAs are 
both subVOAs of the
Moonshine module. It is tempting to guess that this is the most natural
home of the Haagerup. Then the character $J$ of the Moonshine module will
equal one of $ch_0+ch_{\mathfrak{b}}+2ch_{\mathfrak{a}}$, $ch_0+ch_{\mathfrak{b}}+2ch_{\mathfrak{c}}$, or $ch_0+ch_{\mathfrak{b}}+ch_{\mathfrak{a}}+
ch_{\mathfrak{c}_0}$, and we find that $(c,d,f)$ is either (12,0,0), (0,12,0)
or (0,0,1), and that $(g,h,j,k,l)$ is either (1,2,0,0,0), (0,1,1,0,0), (0,2,0,1,0)
or (0,1,0,0,1), and all other parameters vanish.

\subsection{Character vectors of the  twisted  Haagerup double}

The modular group representation for the $\omega\ne 0$ twist of the Haagerup decomposes into
a sum of 3 irreps: the trivial one, a 4-dimensional one
with kernel $\Gamma(9)$, and a 7-dimensional one with kernel $\Gamma(13)$.
Only the 4-dimensional irrep is new: it is handled by e.g. the root lattice $A_8$.
Actually, 6 different 4-dimensional irreps arise here, varying with the twist and central
charge mod 24, but they are all obtained from the one in $A_8$ by some combination of
3rd roots of unity and taking the contragredient. The effects on $\Lambda,\chi$ of both of
these is discussed at the end of Section 1.4.

\subsubsection{The 1-twisted Haagerup at central charge $c\equiv_{24}8$}

The matrices $\Lambda,\chi$ for  $\cD^1\bbZ_3$ and $\cD^1\Hg_3$ are respectively
$\textrm{diag}(2/3,2/3,7/9,1/9,4/9)$,
\begin{align}&{\scriptsize\left(\begin{array}{ccccc}
80 & 168 & 54& 9504& 1078\\
84 &  164&  -27 &  -4752& -539\\
126 &    -126 &     8 &   -17248 &  1375\\
9&    -9&   -8 &  96 & 10\\
36 &      -36 &   20 &  392&  -340\end{array}\right)}\,,\nonumber\\
&\mathrm{diag}(-1/3,2/3,2/3,7/9,1/9,4/9,5/39,20/39,32/39,41/39,8/39,11/39)\,,\nonumber\\
&{\scriptsize\left(\begin{array}{cccccccccccc}
-42/5&1&0&0&0&0&-17/5&-13/5&3/5&1/5&-7&14/5\\
209456/5&80&84&81&14256&1620&-19054/5&-2366/5&166/5&77/5&-2814&-5812/5\\
228912/5&168&164&-81&-14256&-1620&-1428/5&-1092/5&252/5&84/5&-588&1176/5\\
176556/5&84&-42&8&-17248&1375&-714/5&-546/5&126/5&42/5&-294&588/5\\
114/5&6&-3&-8&96&10&-51/5&-39/5&9/5&3/5&-21&42/5\\
11391/5&24&-12&20&392&-340&-204/5&-156/5&36/5&12/5&-84&168/5\\
546/5&0&0&0&0&0&-559/5&-161/5&16/5&7/5&14&418/5\\
52611/5&0&0&0&0&0&-8134/5&-86/5&206/5&42/5&-486&-4617/5\\
600457/5&0&0&0&0&0&-12733/5&3388/5&412/5&119/5&-6868&24871/5\\
2851862/5&0&0&0&0&0&100947/5&6118/5&437/5&229/5&31027&-49514/5\\
1729/5&0&0&0&0&0&649/5&-174/5&14/5&13/5&-238&437/5\\
5252/5&0&0&0&0&0&2002/5&-567/5&77/5&14/5&27&76/5\end{array}\right)}\,.\nonumber\end{align}
At $c=8$, $\cD^1\bbZ_3$ is realised by the lattice VOA $\cV(A_8)$. $\cD^{+,1}S_3$ is realised
by the affine algebra VOA $\cV(B_{4,2})$. Both of these are conformal embeddings.

\subsubsection{The 1-twisted Haagerup at central charge $c\equiv_{24}16$}

The matrices $\Lambda,\chi$ for $\cD^1\bbZ_3$ and $\cD^1\Hg_3$ are respectively
$\textrm{diag}(1/3,1/3,4/9,7/9,1/9)$,
\begin{align}&\chi= {\scriptsize\left(\begin{array}{ccccc}
166 & 330& 198& 18& 900\\
165 &   331 &  -99&  -9&  -450\\
351&  -351 &  56 &  -25 & 1625\\
2079& -2079&-1694& 53 & 3146\\
27&  -27 &  28&  1&  -102\end{array}\right)}\,,\nonumber\\
&\mathrm{diag}(-2/3,1/3,1/3,4/9,7/9,1/9,31/39,7/39,19/39,28/39,34/39,-2/39)\,,\nonumber\\
&{\scriptsize\left(\begin{array}{cccccccccccc}3&1&0&0&0&0&-1&0&1&-1&1&-1\\
6732&166&165&297&27&1350&-165&-286&87&-165&152&-880\\
7359&330&331&-297&-27&-1350&-165&0&165&-165&165&-165\\
8199&234&-117&56&-25&1625&-117&0&117&-117&117&-117\\
173727&1386&-693&-1694&53&3146&-693&0&693&-693&693&-693\\
126&18&-9&28&1&-102&-9&0&9&-9&9&-9\\
498225&0&0&0&0&0&-1946&-6250&1800&-1875&1925&41250\\
858&0&0&0&0&0&-54&58&88&-33&45&-342\\
31603&0&0&0&0&0&-374&1650&480&-425&374&3553\\
263120&0&0&0&0&0&-1265&5313&529&-1284&1288&-16720\\
920556&0&0&0&0&0&-2673&-2025&1848&-2574&2704&22572\\
13&0&0&0&0&0&0&-9&7&-5&3&1\end{array}\right)}\,.\nonumber\end{align}

\subsubsection{The 1-twisted Haagerup at central charge $c\equiv_{24}0$}

The matrices $\Lambda,\chi$ for $\cD^1\bbZ_3$ and $\cD^1\Hg_3$ are respectively
$\textrm{diag}(1,0,1/9,4/9,7/9)$, 
\begin{align}&{\scriptsize\left(\begin{array}{ccccc}
0&131274&185328&13770&324\\ 
1&6&-135&-54&-27\\ 0&-27&24&-35&10\\ 
0&-594&-2002&290&11\\ 0&-5967&14144&340&-64\end{array}\right)}\,,\nonumber\\
&\mathrm{diag}(0,1,0,1/9,4/9,7/9,6/13,11/13,2/13,5/13,7/13,8/13)\,,\nonumber\\
&{\scriptsize\left(\begin{array}{cccccccccccc} 
-12&1&0&0&0&0&20&2&10&8&12&8\\ 60033&0&65637&92664&6885&162&-18997&-16&-3200&-350&-1617&-792\\
 -18&2&6&-135&-54&-27&20&2&10&8&12&8\\
27&0&-27&24&-35&10&0&0&0&0&0&0\\
594&0&-594&-2002&290&11&0&0&0&0&0&0\\ 5967&0&-5967&14144&340&-64&0&0&0&0&0&0\\ 
1716& 0& 0& 0 &0& 0 &-252& -5& -176& 176& 134& 100\\
22451& 0& 0& 0& 0& 0 &-980& 8 &16556&2464& -96& -385\\
104& 0& 0&0& 0& 0& 0&5& 44&-56&-15& 20\\
910& 0& 0& 0& 0&0& 77& 8&-847&-32&120&-34\\
3003& 0& 0& 0& 0&0& 330& 0&-1540& 704& -216&55\\
5200& 0& 0& 0& 0&0&616&-14&3388&-595&132&-20\end{array}\right)}\,.\nonumber\end{align}

\subsubsection{The 2-twisted Haagerup at central charge $c\equiv_{24}8$}

The matrices $\Lambda,\chi$ for $\cD^2\bbZ_3$ and $\cD^2\Hg_3$ are respectively
$\mathrm{diag}(-1/3,2/3,8/9,2/9,5/9)$, 
\begin{align}&{\scriptsize\left(\begin{array}{ccccc}2&1&1/3&-1/3&-2/3\\ 46683&248&26&-836&-133\\ 
706401&0&156&2584&475\\ 2187&0&14&172&-77\\ 56862&0&65&-2431&-74\end{array}\right)}\,,\nonumber\\
&\mathrm{diag}(-1/3,-1/3,2/3,8/9,2/9,5/9,5/39,20/39,32/39,41/39,8/39,11/39)\,,\nonumber\\
&{\scriptsize\left(\begin{array}
{cccccccccccc}-16/5&26/5&1/2&1/6&-1/6&-1/3&-17/10&-13/10&3/10&1/10&-7/2&7/5\\ 
26/5&-16/5&1/2&1/6&-1/6&-1/3&17/10&13/10&-3/10&-1/10&7/2&-7/5\\ 
46683&46683&248&26&-836&-133&0&0&0&0&0&0\\ 
706401&706401&0&156&2584&475&0&0&0&0&0&0\\
2187&2187&0&14&172&-77&0&0&0&0&0&0\\
56862&56862&0&65&-2431&-74&0&0&0&0&0&0\\
546/5&-546/5&0&0&0&0&-559/5&-161/5&16/5&7/5&14&418/5\\
 52611/5&-52611/5&0&0&0&0&-8134/5&-86/5&206/5&42/5&-486&-4617/5\\
 600457/5&-600457/5&0&0&0&0&-12733/5&3388/5&412/5&119/5&-6868&24871/5\\
 2851862/5&-2851862/5&0&0&0&0&100947/5&6118/5&437/5&229/5&31027&-49514/5\\
 1729/5&-1729/5&0&0&0&0&649/5&-174/5&14/5&13/5&-238&437/5\\
 5252/5&-5252/5&0&0&0&0&2002/5&-567/5&77/5&14/5&27&76/5\end{array}\right)}\,.\nonumber\end{align}

\subsubsection{The 2-twisted Haagerup at central charge $c\equiv_{24}16$}

The matrices $\Lambda,\chi$ for $\cD^2\bbZ_3$ and $\cD^2\Hg_3$ are respectively
$\mathrm{diag}(1/3,1/3,5/9,8/9,2/9)$,
\begin{align}&{\scriptsize\left(\begin{array}{ccccc}160& 336 & 32/3&2/3&196/3\\ 
  168& 328 &  -16/3  &   -1/3&  -98/3\\
5832 &    -5832& -272 &  2&1925\\
32076&    -32076&  196 &  12&-15092\\
729&  -729&  40 &  -4&28\end{array}\right)}\,,\nonumber\\
&\mathrm{diag}(-2/3,1/3,1/3,5/9,8/9,2/9,31/39,7/39,19/39,28/39,34/39,-2/39)\,,\nonumber\\ 
&{\scriptsize\left(\begin{array}{cccccccccccc}3&1&0&0&0&0&-1&0&1&-1&1&-1\\
6717&160&168&16&1&98&-162&-286&84&-162&149&-877\\
7374&336&328&-16&-1&-98&-168&0&168&-168&168&-168\\
221859&3888&-1944&-272&2&1925&-1944&0&1944&-1944&1944&-1944\\
3788856&21384&-10692&196&12&-15092&-10692&0&10692&-10692&10692&-10692\\
5589&486&-243&40&-4&28&-243&0&243&-243&243&-243\\
498225&0&0&0&0&0&-1946&-6250&1800&-1875&1925&41250\\
858&0&0&0&0&0&-54&58&88&-33&45&-342\\
31603&0&0&0&0&0&-374&1650&480&-425&374&3553\\
263120&0&0&0&0&0&-1265&5313&529&-1284&1288&-16720\\
920556&0&0&0&0&0&-2673&-2025&1848&-2574&2704&22572\\
13&0&0&0&0&0&0&-9&7&-5&3&1\end{array}\right)}\,.\nonumber\end{align}

\subsubsection{The 2-twisted Haagerup at central charge $c\equiv_{24}0$}

The matrices $\Lambda,\chi$ for $\cD^2\bbZ_3$ and $\cD^2\Hg_3$ are respectively
$\textrm{diag}(1,0,8/9,5/9,2/9)$,
\begin{align}&{\scriptsize\left(\begin{matrix}0&131274&6&528&9282\\ 
1&-3&-1&-7&-8\\ 0&-104247&3&-1001&12376\\ 
0&-12393&-7&232&476\\ 0&-729&5&22&-224\end{matrix}\right)}\,,\nonumber\\
&\mathrm{diag}(0,0,1,2/9,5/9,8/9,6/13,11/13,2/13,5/13,7/13,8/13)\,,\nonumber\\
&{\scriptsize\left(\begin{array}{cccccccccccc}-15/2&9/2&1/2&4&7/2&1/2&6&1&10&5&4&4\\
9/2&-15/2&1/2&4&7/2&1/2&-6&-1&-10&-5&-4&-4\\ 
65637&65637&0&-4641&-264&-3&0&0&0&0&0&0\\ 729&729&0&-224&22&5&0&0&0&0&0&0\\
12393&12393&0&476&232&-7&0&0&0&0&0&0\\ 
104247&104247&0&12376&-1001&3&0&0&0&0&0&0\\ 
1716&-1716&0&0&0&0&-252&-5&-176&176&134&100\\
22451&-22451&0&0&0&0&-980&8&16556&2464&-96&-385\\ 
104&-104&0&0&0&0&0&5&44&-56&-15&20\\ 910&-910&0&0&0&0&77&8&-847&-32&120&-34\\
3003&-3003&0&0&0&0&330&0&-1540&704&-216&55\\
5200&-5200&0&0&0&0&616&-14&3388&-595&132&-20\end{array}\right)}\,.\nonumber\end{align}

\newcommand\biba[7]   {\bibitem{#1} {#2:} {\sl #3.} {\rm #4} {\bf #5,}
                    {#6 } {#7}}
                    \newcommand\bibx[4]   {\bibitem{#1} {#2:} {\sl #3} {\rm #4}}

\def\ASENS            {Ann. Sci. \'Ec. Norm. Sup.}
\def\AM   {Acta Math.}
   \def\AnM              {Ann. Math.}
   \def\CMP              {Commun.\ Math.\ Phys.}
   \def\IJM              {Internat.\ J. Math.}
   \def\JAMS             {J. Amer. Math. Soc.}
\def\JFA              {J.\ Funct.\ Anal.}
\def\JMP              {J.\ Math.\ Phys.}
\def\JRA              {J. Reine Angew. Math.}
\def\JPAA             {J.\ Pure Appl.\ Algebra.}
\def\JSP              {J.\ Stat.\ Physics}
\def\LMP              {Lett.\ Math.\ Phys.}
\def\RMP              {Rev.\ Math.\ Phys.}
\def\RNM              {Res.\ Notes\ Math.}
\def\RIMS             {Publ.\ RIMS.\ Kyoto Univ.}
\def\Inv              {Invent.\ Math.}
\def\npbp             {Nucl.\ Phys.\ {\bf B} (Proc.\ Suppl.)}
\def\nupb             {Nucl.\ Phys.\ {\bf B}}
\def\nup              {Nucl.\ Phys. }
\def\nupp             {Nucl.\ Phys.\ (Proc.\ Suppl.) }
\def\adma             {Adv.\ Math.}
\def\coma             {Con\-temp.\ Math.}
\def\PAMS             {Proc. Amer. Math. Soc.}
\def\PJM              {Pacific J. Math.}
\def\ijmp             {Int.\ J.\ Mod.\ Phys.\ {\bf A}}
\def\jpa              {J.\ Phys.\ {\bf A}}
\def\PLB              {Phys.\ Lett.\ {\bf B}}
\def\RIMS             {Publ.\ RIMS, Kyoto Univ.}
\def\Top               {Topology}
\def\TAMS             {Trans.\ Amer.\ Math.\ Soc.}
\def\Duke              {Duke Math.\ J.}
\def\K                 {K-theory}
\def\JOP               {J.\ Oper.\ Theory}
\def\JKT               {J.\ Knot Theory and its Ramifications}

\vspace{0.2cm}\addtolength{\baselineskip}{-2pt}
\begin{footnotesize}
\noindent{\it Acknowledgement.}

The authors thank the University of Alberta Mathematics Dept, Cardiff School of Mathematics, 
Swansea University Dept of Computer Science, 
and Universit\"at  W\"urzburg Institut f\"ur Mathematik for generous
hospitality while researching this
paper. Their research was supported in part by  
EU-NCG Research Training Network: MRTN-CT-2006 031962, DAAD (Prodi Chair),
and NSERC. We  thank Matthias Gaberdiel, Pinhas Grossman, Paulo Pinto, Ingo Runkel, and 
Feng Xu for discussions.

\end{footnotesize}

\end{document}